\documentclass{article}
\usepackage{animate}
\usepackage[utf8]{inputenc}
\usepackage{amssymb}
\usepackage{graphicx}
\usepackage{tikz}
\usepackage[margin=1in]{geometry}
\usepackage[absolute,overlay]{textpos}
\usepackage{graphicx}
\usepackage{hyperref} 
\usepackage[numbers,sort]{natbib}

\title{Eigenvalue collisions for
periodic matrix families associated with Ginibre matrices}
\author{Carlos Vargas\thanks{homepage: \url{https://carlos-vargas-math.github.io/}, e-mail: \texttt{carlos.vargas.freemath@gmail.com}}}
\date{\today}

\begin{document}
	\maketitle	
	\begin{abstract}
		We study, count and locate the exceptional points where eigenvalues collide 
		for certain families of matrices 
		$$R(s,t) = \cos(s \pi / 2)C + \sin(s \pi / 2)U(t), \quad  s,t \in [0,1]$$ 
		where $C$ is a realization of a Ginibre random matrix, or a closely related matrix, 
		and $U(t)$ is a $t$-periodic diagonal matrix whose eigenvalues move along 
		the unit circle or other parameterized curves (simple or not).
	
		To do this, we track eigenvalues continuously 
		along simple curves, searching for eigenvalue discrepancies after looping, 
		indicating such collision-events inside said curves. 
		
		Although stochastic in nature, the collision count for the circle cases repeatedly
		yielded $N(N-1)$ under fairly general circumstances, where $N$ is the matrix dimension. 
		More collisions were observed for the elliptic generalization and for other curves.
	
		We include a package to compute and visualize these results.	
	\end{abstract} 

	\section{Introduction} \label{section:introduction}
	We consider the two-parameter family of matrices
	\[
	R(s,t) = \cos(s \pi / 2)C + \sin(s \pi / 2)U(t), \quad  s,t \in [0,1]
	\]
	where $C \in M_N(\mathbb{C})$ is a realization of a Ginibre matrix, 
	and $U(t)$ is a $t$-periodic diagonal matrix with eigenvalues 
	flowing along the unit circle, or more generally, along a curve. 
	In the circular case, we set
	\[
	U(t) = \mathrm{diag}\left(\omega^{tN}, \omega^{tN+1}, \omega^{tN+2}, \dots, \omega^{tN+N-1}\right),
	\]
	where $\omega$ is the first non-trivial counter-clockwise $N$-th root of unity. 
	In this case the model may be re-written as
	$$R(s,t) = \cos(s \pi / 2)C + \sin(s \pi / 2)\exp(2i \pi t)U(0).$$

	The random matrix $R(s,0)$ has been studied in the large-$N$ limit. 
	In particular, Zhong \cite{Zhong-2022} computed its asymptotic spectral distribution, using ideas from Haagerup and Larsen \cite{Haagerup-Larsen-2010} 
	on the Brown measure of $R$-diagonal elements. 

	In this work, we focus on how the eigenvalues evolve as $s$ and $t$ vary, 
	with particular attention to locating the exceptional points $(s,t)$ where eigenvalues collide. 
	Since the dependence on $(s,t)$ is analytic (see Section~\ref{subsection:techincal-remarks}), 
	eigenvalues can be continuously tracked.

	We are especially interested in the \emph{monodromy} of the eigenvalues 
	arising from the evolution in $t$: 
	For fixed $s$, as $t$ runs from $0$ to $1$, the eigenvalues of $R(s,t)$ 
	rotate counter-clockwise at speeds inversely correlated to their moduli.
	After one full period in $t$, the eigenvalues must return to the 
	original configuration, but they don't generally end up at their 
	same starting positions, 
	thus defining a permutation $\sigma(s)$, for all $s\in [0,1]$.

	At $s=0$ and $s=1$, this permutation is trivial, but for intermediate values of $s$, 
	$\sigma(s)$ is a piece-wise constant process, 
	with different values explained by composing consecutive transpositions 
	associated to eigenvalue collisions.
	By tracking the evolution of $\sigma(s)$ finely across $s \in [0,1]$, 
	we locate the points where $\sigma(s)$ jumps, corresponding to collisions.
	
	We repeatedly observed a count of $N(N-1)$ 
	of such collision-events for the circular case, 
	an increased (stochastic) number of collisions for other curves, 
	and also some variations in the elliptic case. 
	The count for the circular case, although surprisingly constant, is indeed stochastic: 
	the commutative version of the model 
	(included in the non-commutative model with overwhelmingly small probability) 
	yields $N(N-1)/2$ collisions (see Section~\ref{subsection:commutative-case}). 
	
	Eigenvalue collisions, more precisely, correspond to points $(s,t,\lambda)$ 
	where $R(s,t)$ has a repeated eigenvalue $\lambda$, of order $2$ or greater.
	We detect these using the monodromy action around small square loops in the $(s,t)$-space: 
	if the eigenvalues do not return to their original ordering after traversing the square loop, 
	then a collision must exist inside the square.
	
	The dynamics of eigenvalues as $s$ and $t$ vary are visually striking: 
	swirling trajectories, local repulsion, and a discrete set of cycle-merging/separating-events.
	When $U(t)$ flows along the circle, the eigenvalues initially revolve 
	around their own centers, before teaming-up to rotate around the origin,
	until finally each eigenvalue loops by itself, as the value of $s$ increases.
	
	We elaborate on the model and its variants in Section~\ref{section:model-description}. 
	We use permutations to detect and count eigenvalue collisions 
	(Section~\ref{section:eigenvalue-permutations-and-collisions}). 
	In Section~\ref{section:collision-statistics} we discuss statistical properties of these collisions 
	and examine the influence of alternative curves,
	initial eigenvalue distributions and eccentricity on the collision counts. 
	Section~\ref{section:about-the-algorithm} briefly describes implementation details.

	Appendix~\ref{appendix:several-curves} 
	provides high-resolution illustrations of the eigenvalue-permutation evolution 
	for $N=10$ under four different scenarios and three different curves.
	Appendix~\ref{appendix:elliptic-cases} presents analogous results for elliptic deformations.
	Appendix~\ref{appendix:repeated-eigenvalues} includes eigenvalue tracks for traceless Bernoulli matrices, 
	highlighting specific cases with repeated eigenvalues.
	Appendix~\ref{appendix:N=100} displays eigenvalue tracks for the circular case with $N=100$.

	We provide a software package for computing and visualizing these phenomena
	\url{https://github.com/carlos-vargas-math/eigenvalue-collisions}.
	To give the reader a general feeling for these eigenvalue dynamics, 
	we prepared some videos showing the evolution and first collisions 
	(w.r.t. the parameter $s$)
	for a case with $N=20$ (available at 
	\url{https://carlos-vargas-math.github.io/eigenvalue-collisions/animation}).

	\section{Acknowledgements}

	The author is grateful to O. Arizmendi, A. Beshenov, D. L\'opez and J. Santos, 
	for constant feedback and moral support to conclude this project.	
	The author acknowledges P. Zhong, for sharing formula (\cite{Zhong-2022}, Th. 8.8)
	which was used in our model, C. D\'iaz-Aguilera and N. Sakuma, 
	for suggesting interesting examples, 
	and R. Speicher, for general remarks that improved this paper. 
	
	The author declares that there is no conflict of interest regarding the publication of this work.

	\section{Model Description} \label{section:model-description}
  
	\subsection{Ginibre Matrix} 

	\begin{figure}[htbp]
		\centering
		\includegraphics[width=0.3\textwidth]{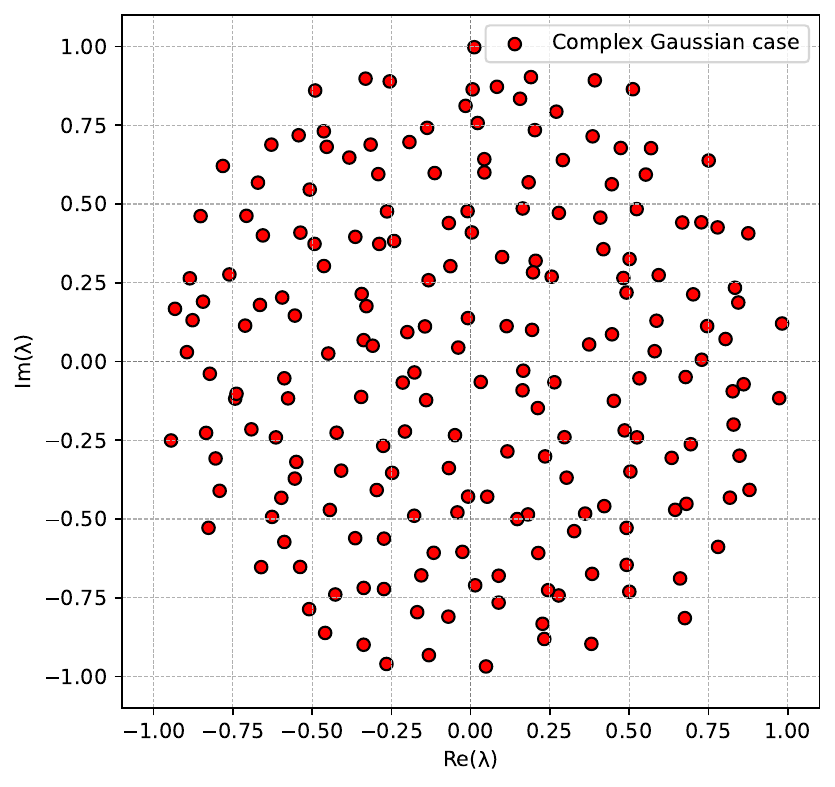}	
		\includegraphics[width=0.3\textwidth]{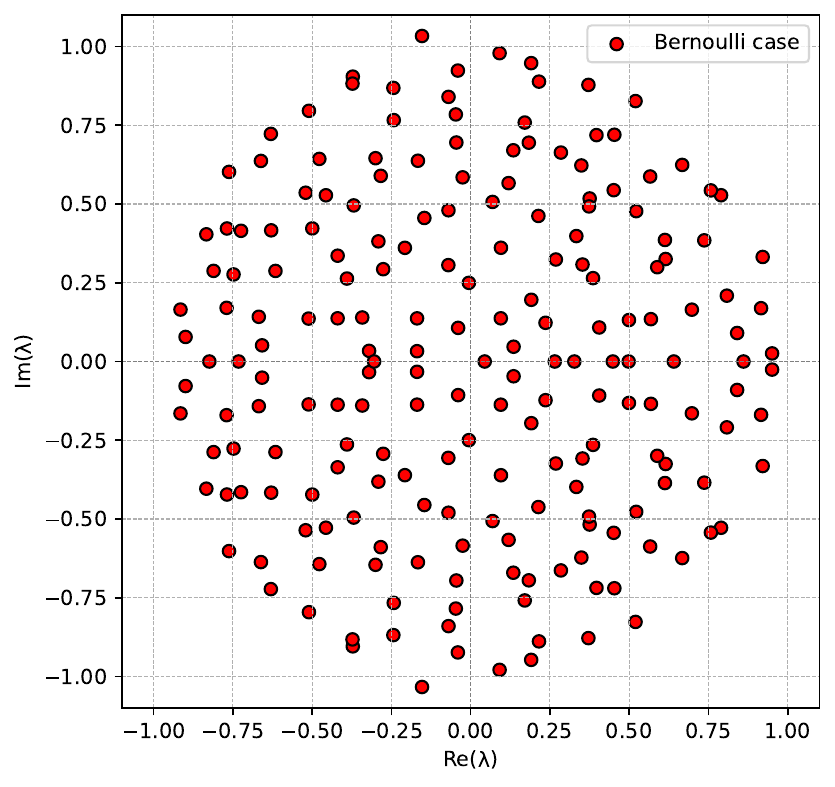}
		\includegraphics[width=0.3\textwidth]{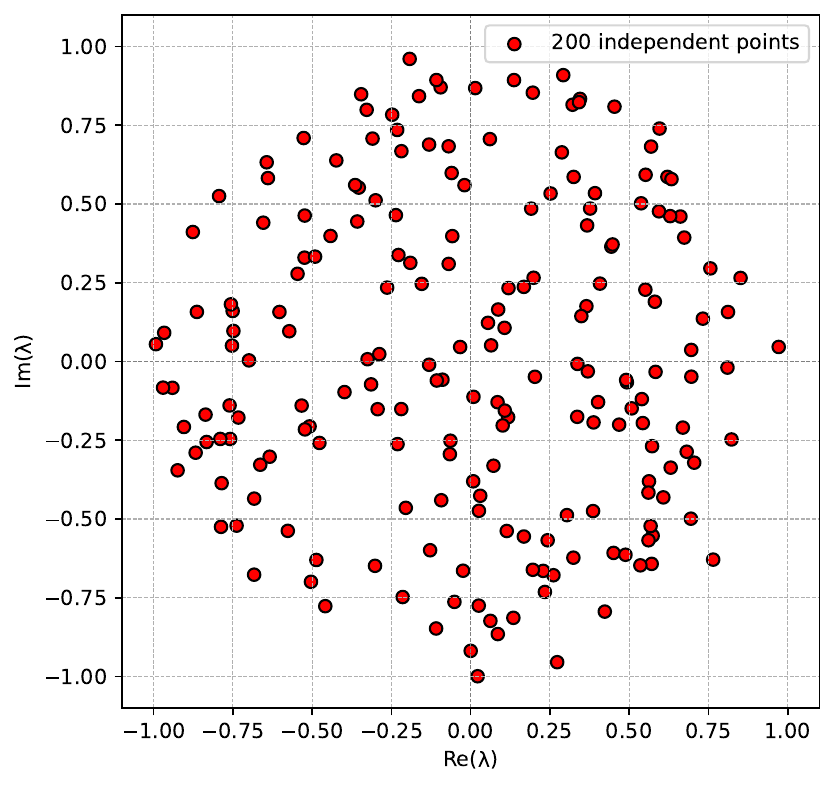}
		\caption{ Eigenvalues of complex Gaussian Ginibre matrix (left), symmetric Bernoulli Ginibre matrix (center) and 
		200 independent points in the unit disc (right)}
		\label{fig:eigenvalues}
	\end{figure}

	A Ginibre matrix is an $N\times N$ random matrix, 
	consisting of independent, identically distributed entries 
	with zero-mean and the same finite variance $\sigma^2$. 
	The well-known \emph{circular law}
	\cite{Ginibre-1965, Girko-1984,Bai-1997,Gotze-Tikhomirov-2010,Tao-Vu-2010}  
	states that if $\sigma^2 = 1/N$, the empirical distribution of the eigenvalues converges almost surely,
	as $N \to \infty$, to the uniform distribution on the centered unit disk.

	Similar to the central limit theorem, the circular law is universal:  
	it holds for matrices with independent, centered, finite-variance entries,  
	regardless of the underlying distribution (e.g., complex Gaussian, real Gaussian, symmetric Bernoulli; 
	see Figure~\ref{fig:eigenvalues}).

	In Figure~\ref{fig:eigenvalues}, we observe that 
	the real-valued matrix has conjugate eigenvalues (as expected),  
	but otherwise, the two eigenvalue distributions (left and center) 
	appear fairly similar and spread uniformly over the unit disk.

	For the complex Gaussian case, Ginibre computed the explicit joint density of eigenvalues in 1965 \cite{Ginibre-1965}:
	$$\rho (\lambda_1, \dots , \lambda_N) = 
	\frac{1}{\pi^n \prod_{k=1}^N k!}
	\exp(-\sum_{k=1}^N |\lambda_{k}|^{2}) 
	\prod_{1\leq j < k \leq N} |\lambda_k - \lambda_j|^2$$

\begin{figure}[htbp]
		\centering

		\includegraphics[width=0.45\textwidth]{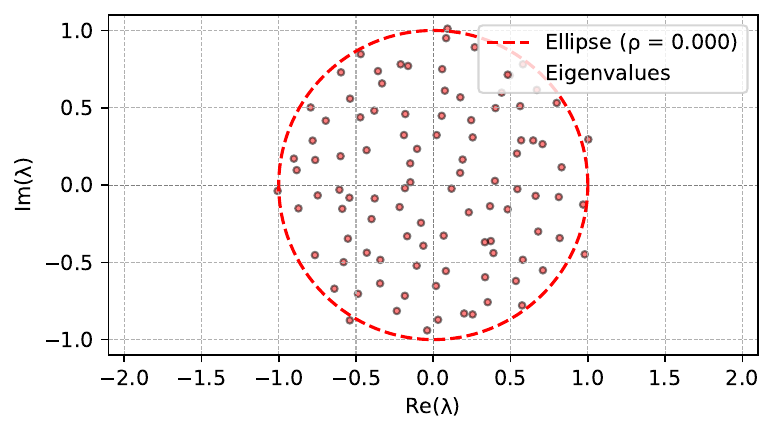}	
		\includegraphics[width=0.45\textwidth]{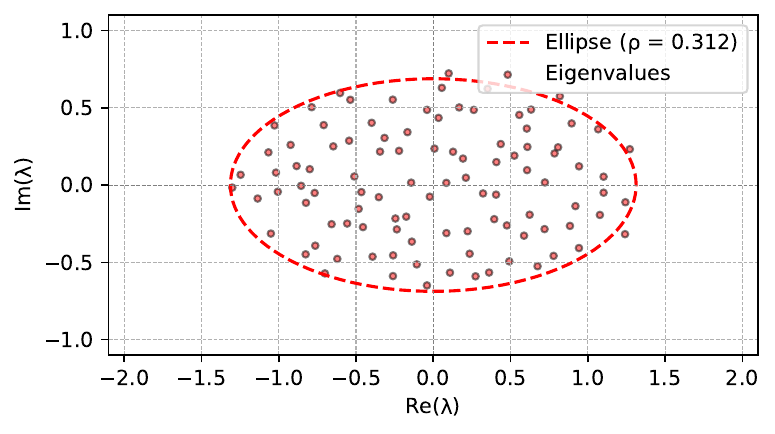}

		\includegraphics[width=0.45\textwidth]{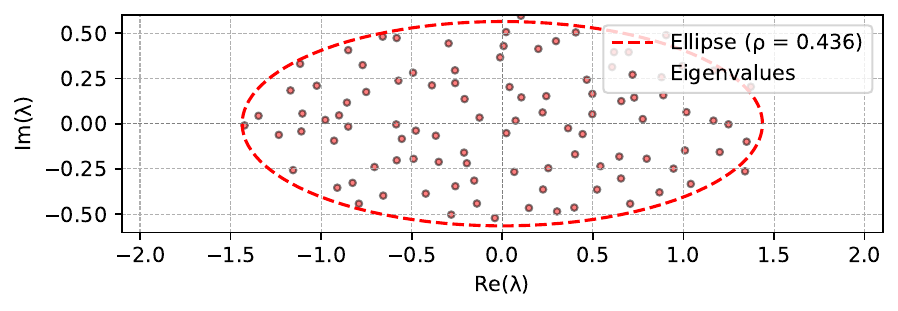}	
		\includegraphics[width=0.45\textwidth]{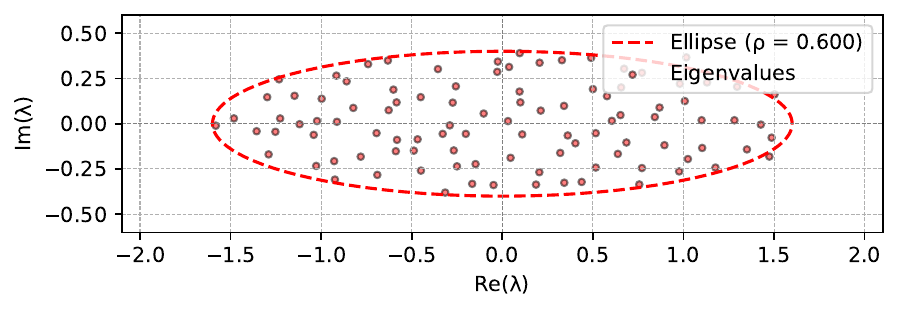}

		\includegraphics[width=0.45\textwidth]{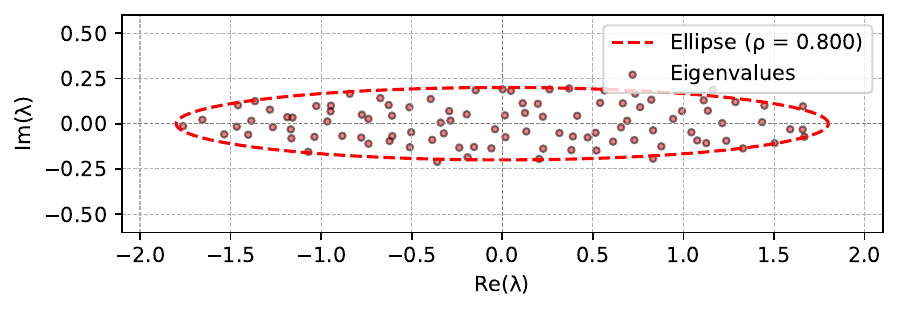}	
		\includegraphics[width=0.45\textwidth]{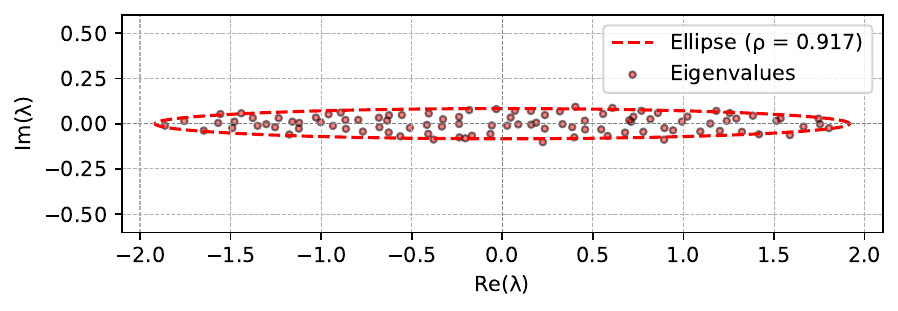}
		\caption{Eigenvalues of $E_x$ for some values of $x\in [0, \pi/4]$ (respectively $\rho =0, 0.312, 0.436, 0.6, 0.8, 0.917$)}
		\label{fig:elliptic}
	\end{figure}

	The eigenvalues of the conjugate transpose $C^*$ of a Ginibre matrix are simply the conjugates of the eigenvalues of $C$.
	More generally, the eigenvalues of the matrix
	\[
	E_x = \cos(x) C + \sin(x) C^*
	\]
	follow the \emph{elliptic law} (see \cite{Naumov-2012,Nguyen-ORourke-2015}):  
	asymptotically, they are uniformly distributed on an ellipse with semi-axes 
	$1+\rho$ and $1-\rho$, where $\rho = \sin(2x)$.  
	Figure~\ref{fig:elliptic} illustrates this behavior for several values of $x$.

	The eigenvalues of Ginibre matrices exhibit \emph{repulsion}  
	(as seen from the Vandermonde determinant factor),  
	leading to highly uniform distributions and faster, 
	stronger convergence to the asymptotic limit,  
	especially when compared to independent random points.  
	In our case, this repulsion is helpful for continuously 
	tracking eigenvalues along $(s,t)$-curves.

	\subsection{Motivation and model description}

	Consider a realization $C \in M_N(\mathbb{C})$ of an $N\times N$ Ginibre matrix and 
	let $\omega$ be the first counter-clockwise non-trivial complex $N$-th root of $1$. 
	For $t\in [0,1]$ define
	$$U(t) = \mathrm{diag}(\omega^{tN}, \omega^{tN+1}, \omega^{tN+2}, \dots   , \omega^{tN+N-1}),$$ 
	be a diagonal matrix with equidistant points along the circle.

	We study the model	
	$$R(s,t) = \alpha(s)C + \beta(s)U(t),\quad  \alpha(s)= \cos((s\pi)/2),\quad \beta(s)= \sin((s\pi) /2),$$ 
	which is periodic in $t$, since $U(0) = U(1)$ implies $R(s,0) = R(s,1)$.

\begin{figure}[htbp]
		\centering
		\includegraphics[width=0.3\textwidth]{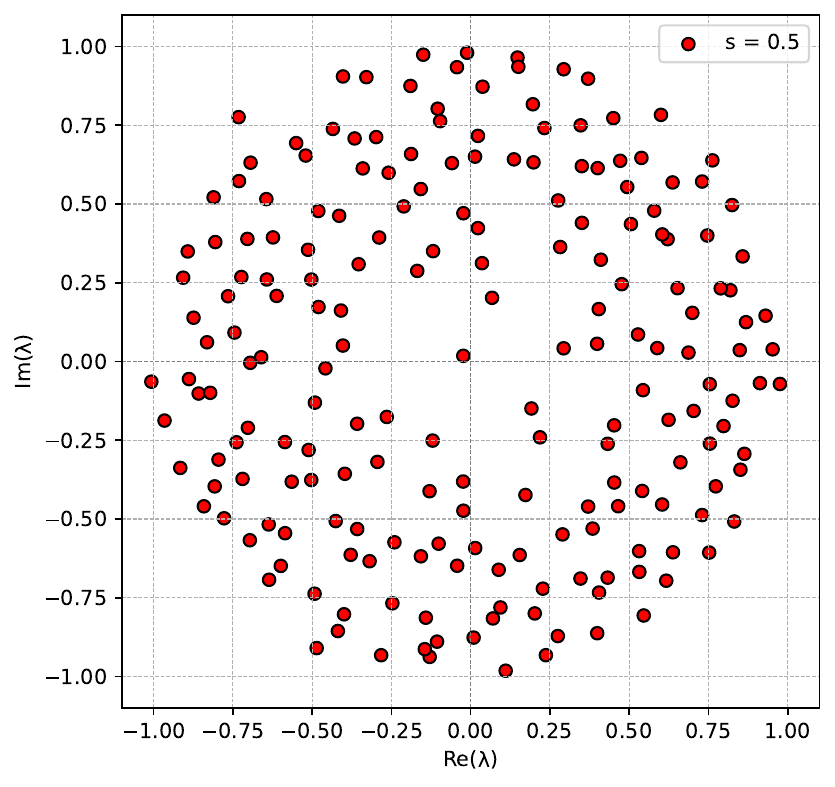}	
		\includegraphics[width=0.3\textwidth]{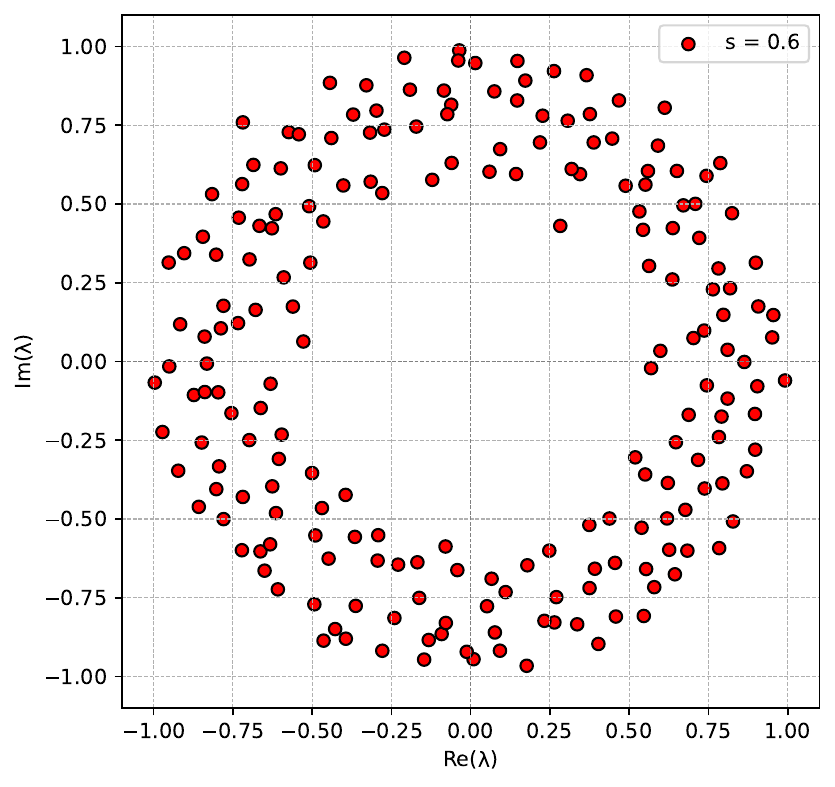}
		\includegraphics[width=0.3\textwidth]{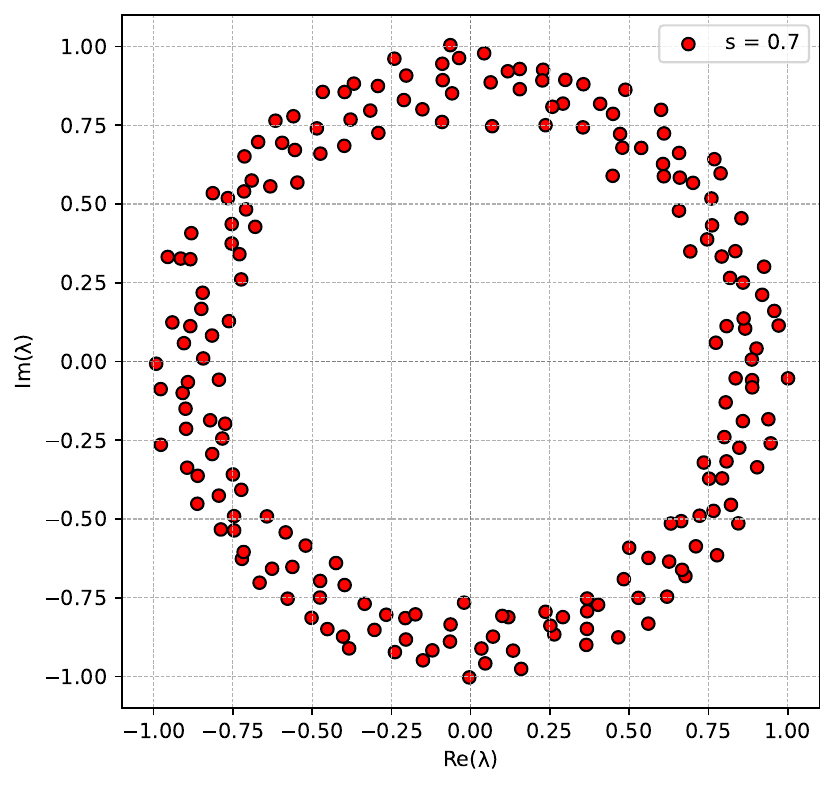}
		\caption{Eigenvalues of $R(s,t)$ for $s = 0.5, 0.6, 0.7$}
		\label{fig:annulus}
	\end{figure}

	The starting point of this investigation is the model $R(s,0)$. 
	For large $N$, it approximates a weighted sum of a circular operator and a free Haar-unitary operator.  
	The asymptotic distribution depends only on $s$ and is independent of $t$,  
	with $s$ controlling the width of the annulus where the eigenvalues are supported 
	(see Figure~\ref{fig:annulus}).

	We noticed that if we replace $U(0)$ 
	with a matrix $V(0)$ having eigenvalues flowing along a different curve,  
	the empirical distributions behave similarly:  
	repelling eigenvalues supported on said curve, 
	thickened by some  $s$-dependent width.
	
	The effect of varying $t$ becomes noticeably interesting 
	when $s$ is large enough and the support is an annulus (not a disk).  
	As $t$ increases, eigenvalues rotate around the curve, 
	repelling each other neatly without collisions.
	
	The main observation behind this investigation is the following:  
	The eigenvalues rotate at different speeds, inversely correlated with their norms.  
	Since the model is periodic, the eigenvalues must return to their original configuration;  
	however, most (especially those near the outer edge) do not return to their initial locations,  
	leading to a nontrivial permutation.

	Computing these permutations allows us to color eigenvalue tracks by permutation cycle size,  
	producing animations tracking eigenvalues continuously, 
	allowing us to visually separate different cycles, yet maintaining the periodicity of the animation,  
	as each eigenvalue ends-up in a same-colored position 
	(see figures in the appendices, or the above-mentioned animations).  

	In general, computing eigenvalue distributions of non-self-adjoint random matrices is difficult.  
	However, matrices in this model approximate R-diagonal elements,  
	whose Brown measures and asymptotic distributions have been computed explicitly  
	(see \cite{Haagerup-Larsen-2010}).

	Depending on $\alpha$ and $\beta$, the asymptotic distribution (as $N \to \infty$)  
	is supported on a centered disk or annulus,  
	with greater density toward the outer edge.  
	Following Zhong \cite{Zhong-2022} (Theorem 8.8),  
	we choose $\alpha^2 + \beta^2 = 1$ to ensure that the (asymptotic) outer radius remains 1.
	
	Thus, we consider the interpolation
	\[
	\alpha(s) = \cos\left( \frac{s\pi}{2} \right), \quad \beta(s) = \sin\left( \frac{s\pi}{2} \right), \quad s \in [0,1],
	\]
	from the circular element to a rotating Haar-unitary element.  
	Eigenvalue behaviors for other positive pairs $(\alpha, \beta)$  
	can be obtained by simple rescaling.	

	\subsection{Increasing $t$}

	In this work, rather than studying asymptotic distributions, 
	we want to draw attention to the more empirical question of detecting eigenvalue collisions 
	as we vary the parameters $s$ and $t$, even for small $N$.

	\begin{figure}[htbp]
		\centering
		\includegraphics[width=0.4\textwidth]{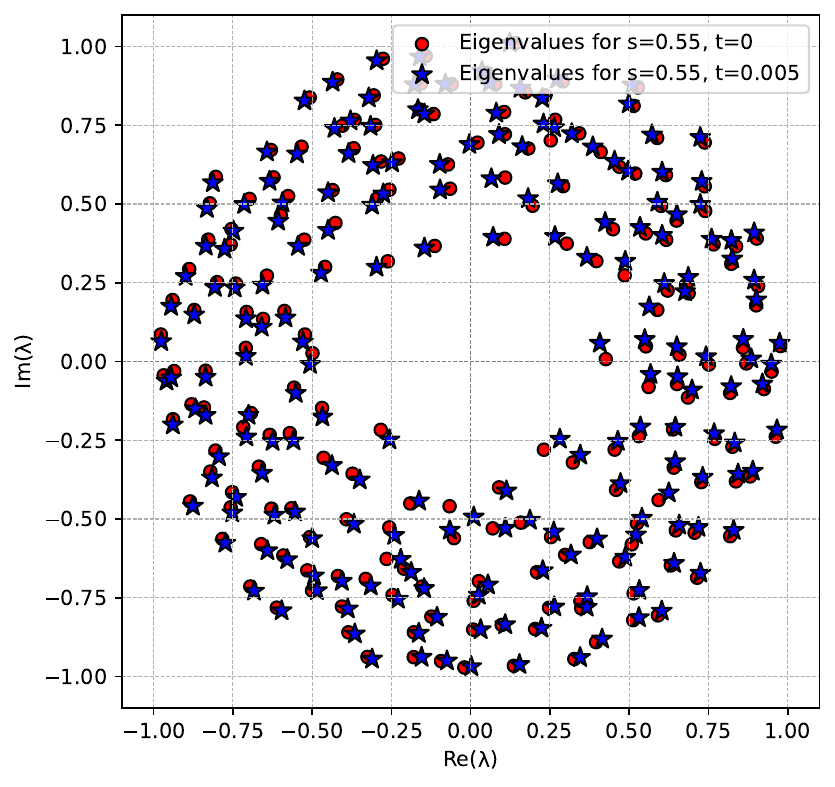}
		\includegraphics[width=0.4\textwidth]{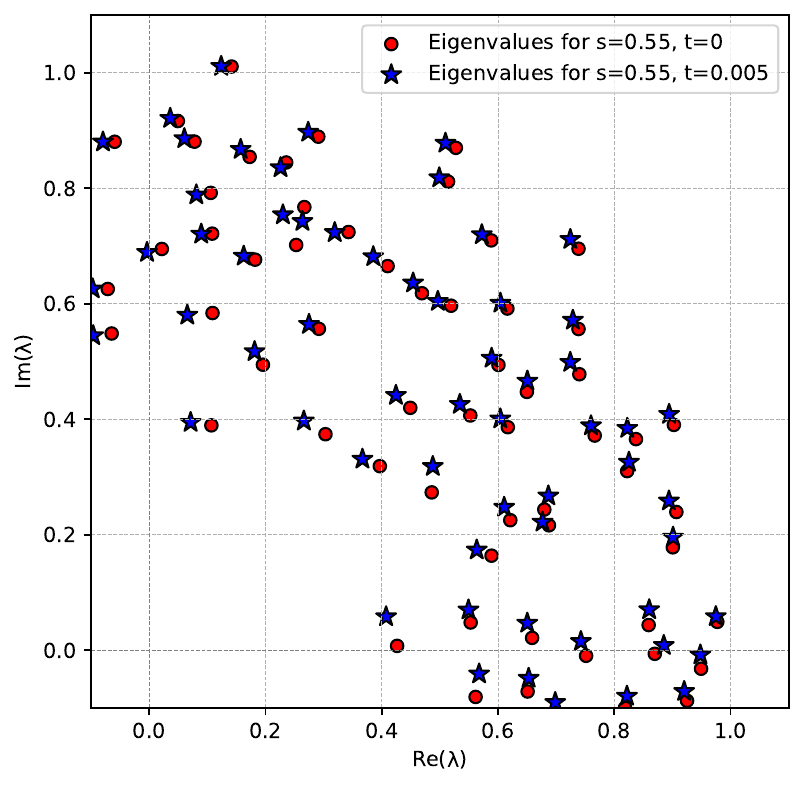}
		   \caption{Eigenvalues of $R(s,t)$ for $t=0$ (red dots) and $t=0.005$ (blue stars), $s=0.55$}
		\label{fig:rotating_circle}
	\end{figure}

	Increasing $t$ induces a counter-clockwise rotation of the eigenvalues. 
	Figure~\ref{fig:rotating_circle} displays the eigenvalues at a fixed $s$ 
	for two nearby values of $t$: the red dots correspond to $t=0$, 
	and the blue stars to $t=0.005$. 
	The eigenvalues near the boundary of the domain move more slowly than those near the center.

	Since $U(0) = U(1)$, as $t$ evolves from $0$ to $1$, 
	each eigenvalue must eventually return to the original position of some eigenvalue, 
	although most will not trace a full cycle individually. 
	This yields a non-trivial permutation $\sigma(s)$ 
	associated with the matrix process at each $s \in [0,1]$.
	When $\sigma(s_0) \neq \sigma(s_1)$ for nearby values $s_0 < s_1$,
	(or more concretely, if their conjugacy classes don't match)  
	a collision is expected to occur for some $(s,t) \in [s_0,s_1] \times [0,1]$.
	To locate such collisions more efficiently, 
	we refine the $(s,t)$-domain by partitioning $[0,1]^2$ into smaller sub-squares 
	rather than stripes (see Section~\ref{subsection:grid-search}).

\begin{figure}[htbp]
		\centering
		\includegraphics[width=0.4\textwidth]{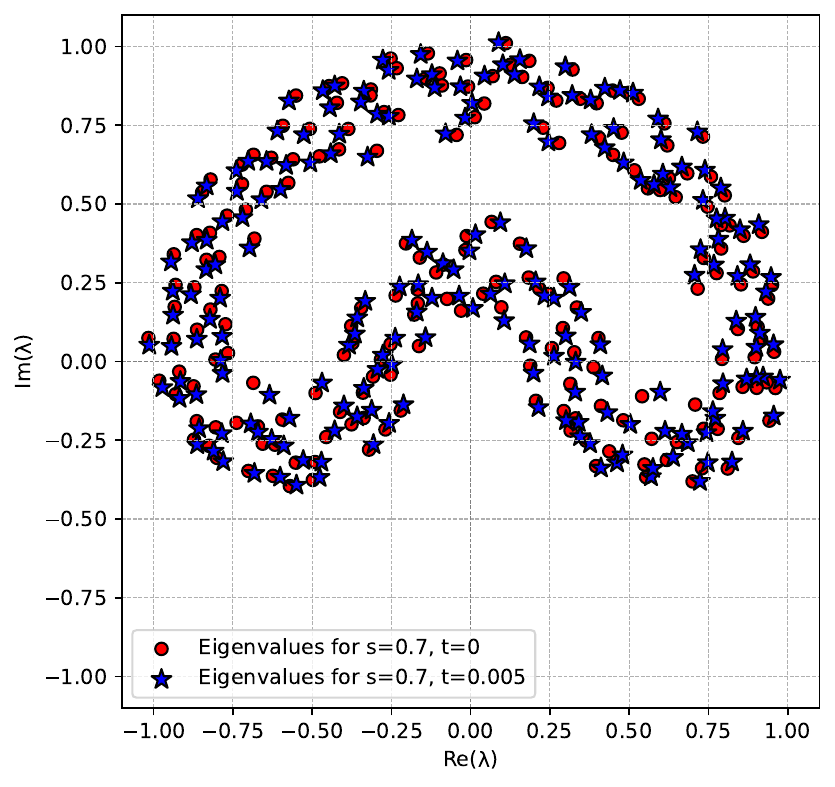}
		\includegraphics[width=0.4\textwidth]{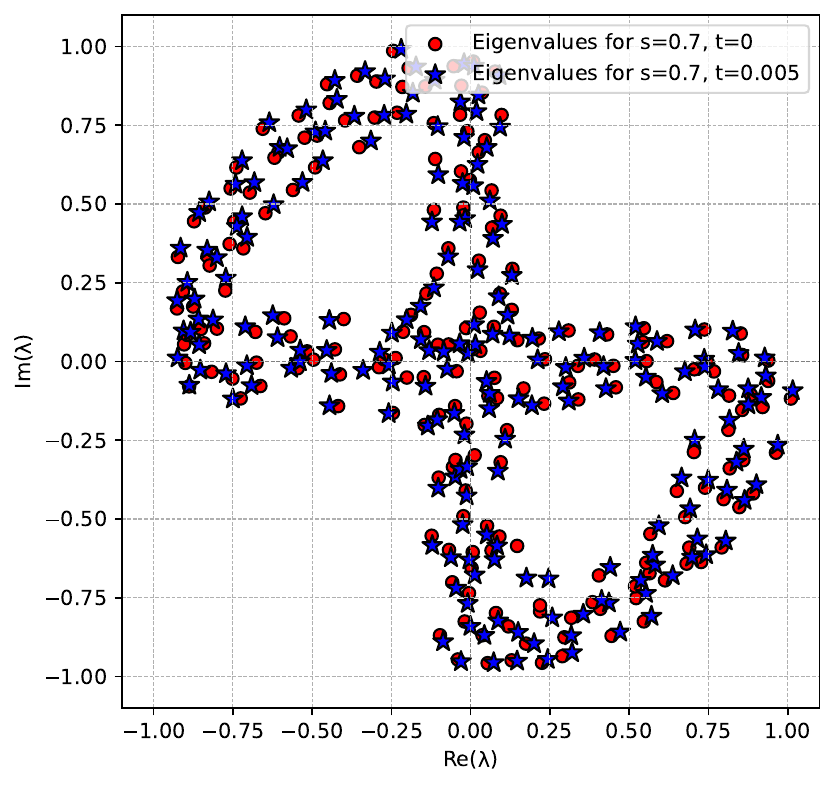}
			   \caption{Circuit curve (left) and crossing curve (right)}
		\label{fig:curves}
	\end{figure}

	The model remains well-posed if the unit circle is replaced
	by any diagonal matrix whose entries lie along a different parametrized curve. 
	Due to the bi-unitary invariance of the Ginibre matrix, 
	the specific shape of the matrix $U$ should not substantially affect the model.
	The diagonal structure provides a convenient way to prescribe eigenvalues 
	flowing along a curve. 
	In this work, we consider three specific examples: 
	the unit circle; 
	a \emph{circuit} curve formed by concatenating the upper unit semicircle with three smaller semicircles 
	of radius $1/3$ (Figure~\ref{fig:curves}, left); 
	and a \emph{crossing} curve flowing along the edges of two opposing sectors (Figure~\ref{fig:curves}, right).

	We chose these curves for their simple parametrizations 
	and for enclosing rational portions of the unit disk's area ($5/9$ and $1/2$, respectively). 
	We also briefly address interesting elliptic deformations of the circle case.

	\subsection{Commutative case} \label{subsection:commutative-case}

	A 'classical' or 'commutative' analog of this model
	$$(1-s)C + sU(t), \quad s,t \in [0,1],$$
	where $C$ is \emph{also a diagonal matrix} (and hence commutes with $U(t)$) is much simpler 
	but still may be worth a few lines. 
	In this case the eigenvalue tracks for fixed $s$ are just circles.
	As $s$ increases, the length of the common radious increases and the centers 
	of the circular tracks get closer to the origin (homothetically).
	In this case there is no repulsion and the eigenvalues move independently, 
	so it suffices to solve the problem for a given pair of eigenvalues.
	
	In general position (i.e. the centers are not the same), 
	two counter-clockwise curves with a fixed non-zero shift difference 
	clash exactly \emph{once} for $(s,t) \in[0,1]^2$.	
	
	Indeed, if we write the difference of two eigenvalues $(j,k)$ with centers in $c_j$ and $c_k$ for $s,t \in [0,1]$ as
	$$
	\lambda_j(s, t) - \lambda_k(s, t) = (1 - s)(c_j - c_k) + s\left( \omega_j(t) - \omega_k(t) \right),
	$$
	the equation for a collision between the eigenvalues \( j \) and \( k \) becomes
	$$
	\frac{1 - s}{s}(c_j - c_k) = -\left( \omega_j(t) - \omega_k(t) \right).
	$$

	Hence, we are comparing a fixed vector $(c_j - c_k)$ 
	(scaled by a decreasing positive number in $[0,\infty]$) 
	to a difference of two points in the circle with a fixed arc length distance
	(which describes itself a circle as $t$ increases).
	Thus, only one value of $s$ gives the right modulus, 
	and only one value of $t$ the right angle.
	
	The same proof works for any flow along a smooth, strictly convex curve instead of $U(t)$. 
	Indeed, the difference of two points at a fixed arc-length 
	distance also attains each possible angle exactly once, for fixed $s$, and $t$ in $[0,1]$.
	Since the modulus of the left-hand side is strictly decreasing on $s$, 
	there will be a single solution for the collision problem.
	Hence, the classical situation yields $N(N-1)/2$ collisions.

	For the free situation (where $C$ is an actual Ginibre matrix), in the circle case,
	we obtained $N(N-1)$ collisions for all of our trials (excluding singular cases).

	This commutative example shows that the number of collisions in the free case is stochastic,
	although we found it to be remarkably stable for the circular case.
	Indeed, although overwhelmingly small, there is a positive probability 
	that the Ginibre matrix is close to diagonal and thus reproduces collision-counts 
	similar to those of the commutative situation.

	\subsection{Technical remarks} \label{subsection:techincal-remarks}

	In the unitary case, the model can be expressed as	
	$$R(z,w) = \cos(z \pi / 2)C + \sin(z \pi / 2)\exp(2i \pi w)U(0),$$
	which is analytic in $(z,w)$. 
	Consequently, the eigenvalues can be tracked continuously, 
	except at a finite number of branching points 
	(see Kato~\cite{Kato-1995}, Section~II-1). 
	Similar behavior is expected when the unit circle is replaced by other smooth curves.

	Our methods rely heavily on tracking eigenvalues from an initial matrix $R(s_0,t_0)$ 
	to a nearby matrix $R(s_1,t_1)$, and, more generally, 
	along sequences in $(s,t)$-space, inserting intermediate points as necessary.
	
	Standard algorithms for computing eigenvalues of non-self-adjoint matrices 
	typically produce eigenvalue lists ordered either by norm 
	or according to internal conventions tied to the underlying numerical methods. 
	Thus, we must perform a matching between two unordered lists of eigenvalues. 
	Our approach uses a greedy matching scheme: 
	each eigenvalue from one list is matched to its closest neighbor in the other. 
	If multiple eigenvalues compete for the same match, the greedy scheme fails. 
	When failure occurs, the $(s,t)$-step is refined by inserting intermediate points.

	Eigenvalue repulsion facilitates that greedy matching generally succeeds 
	for sufficiently small (but not too small), step sizes, 
	with few refinement steps. 
	For smooth eigenvalue trajectories in our visualizations, 
	we used at least $1000$ $t$-steps.
	
	The algorithm encounters difficulties when elliptic deformations 
	are too close to Hermitian, 
	in which real eigenvalues oscillate with $t$. 
	Nevertheless, it remains robust for moderately non-Hermitian cases (e.g., $\rho = 0.866$). 
	Our monodromy-based tracking method 
	fails when the system lacks sufficient dimensionality.

	\section{Eigenvalue permutations and collisions} \label{section:eigenvalue-permutations-and-collisions}	

	We now describe how eigenvalue collisions can be detected 
	using the monodromy action of eigenvalue permutations, 
	focusing on the case where the curve is the unit circle. 
	The same method applies to other curves.

	\begin{figure}[htbp]
		\centering
		\includegraphics[width=0.3\textwidth]{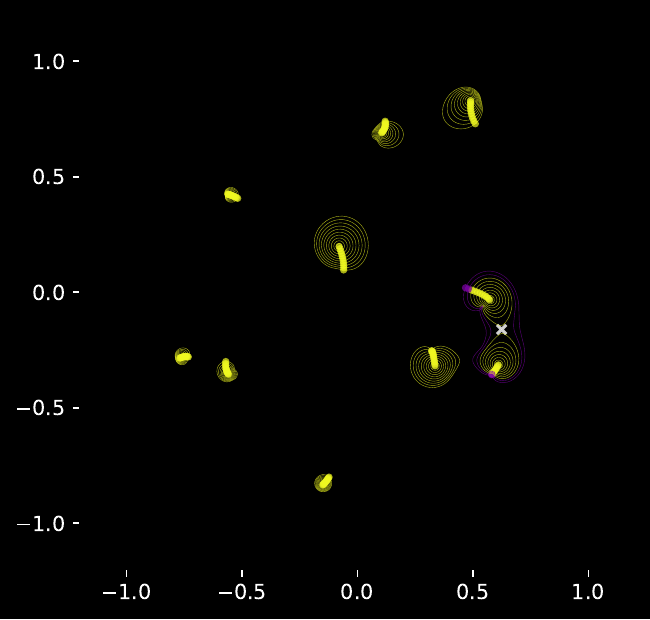}
		\includegraphics[width=0.3\textwidth]{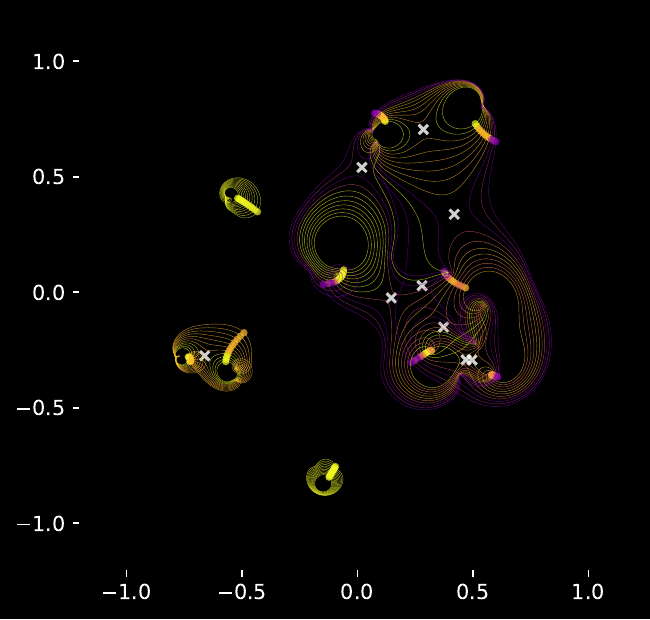}
		\includegraphics[width=0.3\textwidth]{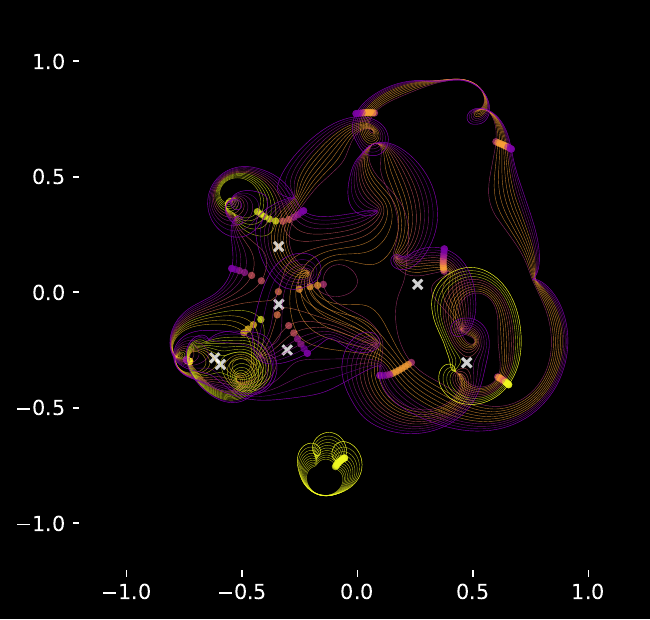}
		\caption{Eigenvalue tracks for three $s$-windows: 
		$[0, 0.01, \dots, 0.10]$ (left), 
		$[0.10, 0.11, \dots, 0.20]$ (center), 		
		$[0.20, 0.21, \dots, 0.30]$ (right)}
		\label{fig:3 s-windows}
	\end{figure}

	Unless $s$ is very close to $0$ (the static case) or near $1$ 
	(perfectly equidistant flow along the curve), 
	increasing the periodic parameter $t$ from $0$ to $1$ 
	induces a nontrivial rotation of the eigenvalues of $R(s,t)$. 
	Outer and inner eigenvalues rotate at different speeds, 
	yielding a non-identity permutation $\sigma(s)$ 
	that relates the eigenvalues of $R(s,0)$ and $R(s,1)$. 
	We study collisions using this permutation.

	Even at small values of $s$, where the rotating annular support is not yet apparent, 
	numerous eigenvalue collisions already occur (see Fig \ref{fig:3 s-windows}). 
	At $s = 0$, the deformation is inactive (as $U$ has no effect), 
	and the eigenvalue paths remain the static eigenvalues of a Ginibre matrix. 
	As $s$ increases, the eigenvalues begin to trace longer paths. 
	In Figure~\ref{fig:3 s-windows}, 
	the color of each path reflects the combinatorial length of the corresponding permutation cycle, 
	ranging from yellow (for the shortest cycles, often singletons) to dark purple (for the longest cycles).

	Recall that composing, or pre-composing a permutation $\sigma$ with a transposition $\tau = (i,j)$ 
	affects the cycle structure in one of two ways:
	\begin{enumerate}
		\item If $i$ and $j$ are in the same cycle of $\sigma$, the transposition splits 
		the cycle into two disjoint sub-cycles: one containing $i$ and one containing $j$, 
		each with consecutive elements of the original cycle.
		\item If $i$ and $j$ belong to different cycles, 
		the transposition merges the cycles into a single one.
	\end{enumerate}

	Strictly speaking, what we are visualizing are not the permutations $\sigma(s)$ themselves, 
	but rather the eigenvalue trajectories colored according to the cycle length. 
	The color changes only when the transposition alters a cycle, 
	so a change in color indicates a collision involving that eigenvalue.

	To speak more precisely about permutations rather than just conjugacy classes, 
	we would need a reference labeling of the eigenvalues. 
	One natural scheme is to label the eigenvalues of $R(s,t)$ at $(s,t) = (0,0)$ (e.g. by norm), 
	then track the labeled eigenvalues along $t=0$ as $s$ increases, 
	and finally compute $\sigma(s)$ by comparing the labeling at $R(s,0)$ and $R(s,1)$. 
	This yields a piece-wise constant permutation process $\sigma(s)$. 
	The outcome depends on the initial labeling: 
	if we label at $t=0.5$ instead, it would result in a different process $\sigma'(s)$, 
	which shares the same conjugacy class with $\sigma(s)$, but will most likely 
	differ as a permutation.

	In this paper, we restrict ourselves to local labeling schemes. 
	One might wonder whether a global labeling exists in which each pair of eigenvalues collides 
	exactly twice (to match the count of $N(N-1)$ collisions; 
	see Section~\ref{section:collision-statistics}). 
	However, we observe that this is not the case for any immediately available global labeling. 
	Whether such a labeling exists remains unclear.

	An interesting direction for future work is to study the statistical properties 
	of the conjugacy classes of the permutations $\sigma(s)$ across different intervals 
	$[s_0, s_1]$. For example, there are stages for the parameter $s$ where long cycles are prevalent
	vs stages with many singletons. 
	Our small software package 
	(available at \url{https://github.com/carlos-vargas-math/eigenvalue-collisions}) 
	allows for such an investigation, 
	though we will not proceed further in this direction.

	\begin{figure}[htbp]
		\centering
		\includegraphics[width=0.3\textwidth]{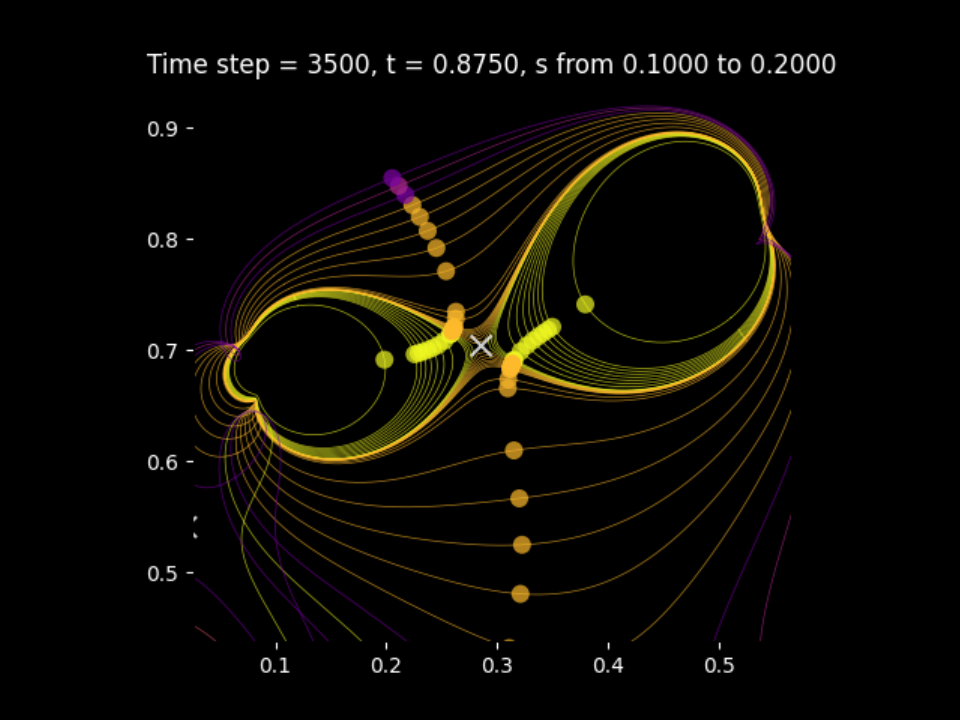}
		\includegraphics[width=0.3\textwidth]{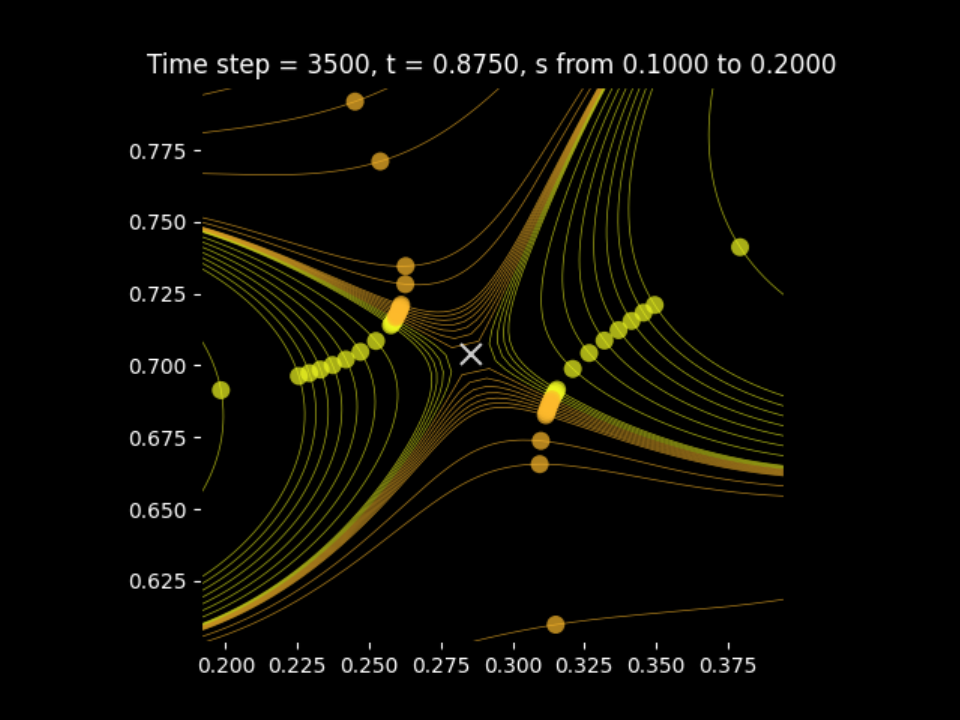}
		\includegraphics[width=0.3\textwidth]{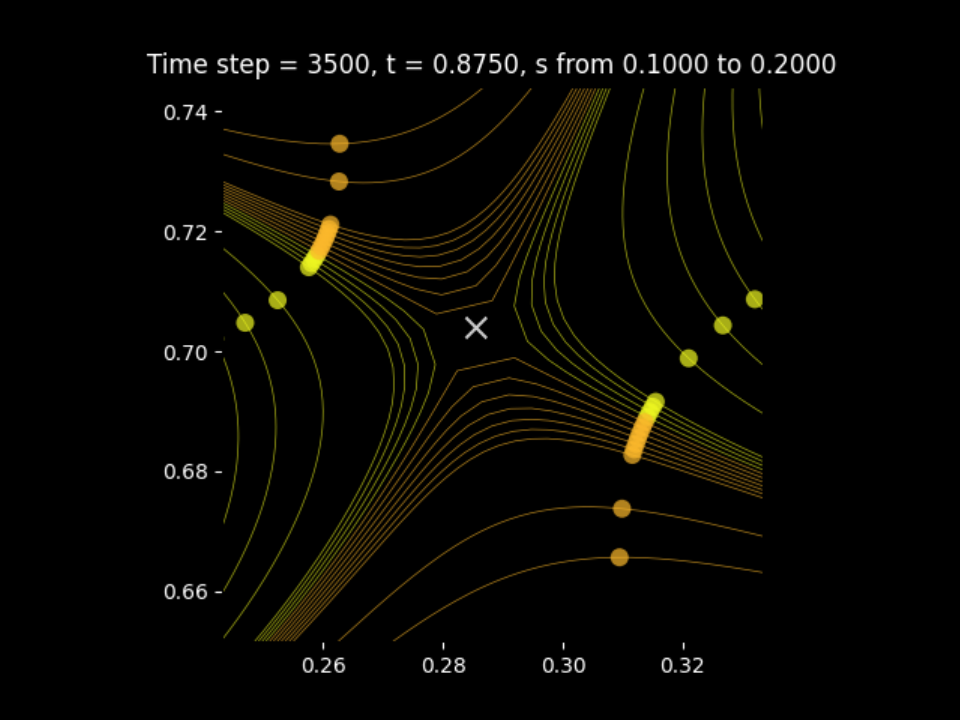}
		\caption{Zoomed-in figures}
		\label{fig:zoomed}
	\end{figure}

	To generate the trajectory plots, we sampled at least $1000$ $t$-steps for each $s$. 
	The animations 
	show the eigenvalue bundles evolving continuously in $t$. 
	Collisions are visible when two bundles of parallel paths appear to exchange eigenvalues, 
	and the cycle length (and hence the color) of the involved tracks change. 
	
	In Figure~\ref{fig:3 s-windows} (left), all but two eigenvalues remain singleton cycles 
	throughout the $s$-window (please zoom-in generously). A clear change in cycle structure occurs between 
	$s = 0.08$ and $s = 0.09$: the permutation $\sigma(s)$ becomes a transposition. 
	This is indicated by the corresponding tracks changing color, 
	and a collision is marked at the estimated $s$-value using our grid-based method.
	
	In the center and right panels of Figure~\ref{fig:3 s-windows}, 
	which show $s$-windows $[0.10, 0.20]$ and $[0.20, 0.30]$, 
	we observe nine and seven collisions, respectively. 
	To analyze one such collision in more detail, 
	we zoom-in using smaller $s$-increments.

	In Figure~\ref{fig:zoomed}, we interpolate additional $s$ values between $0.111$ and $0.119$, 
	and even further refine the interval $[0.117, 0.118]$ to steps of $10^{-4}$. 
	Between $s = 0.1173$ and $s = 0.1174$, two eigenvalue tracks visibly swap, 
	confirming the collision. Away from the collision, 
	the tracks remain nearly unchanged at this resolution.

	\begin{figure}[htbp]
		\centering
		\includegraphics[width=0.3\textwidth]{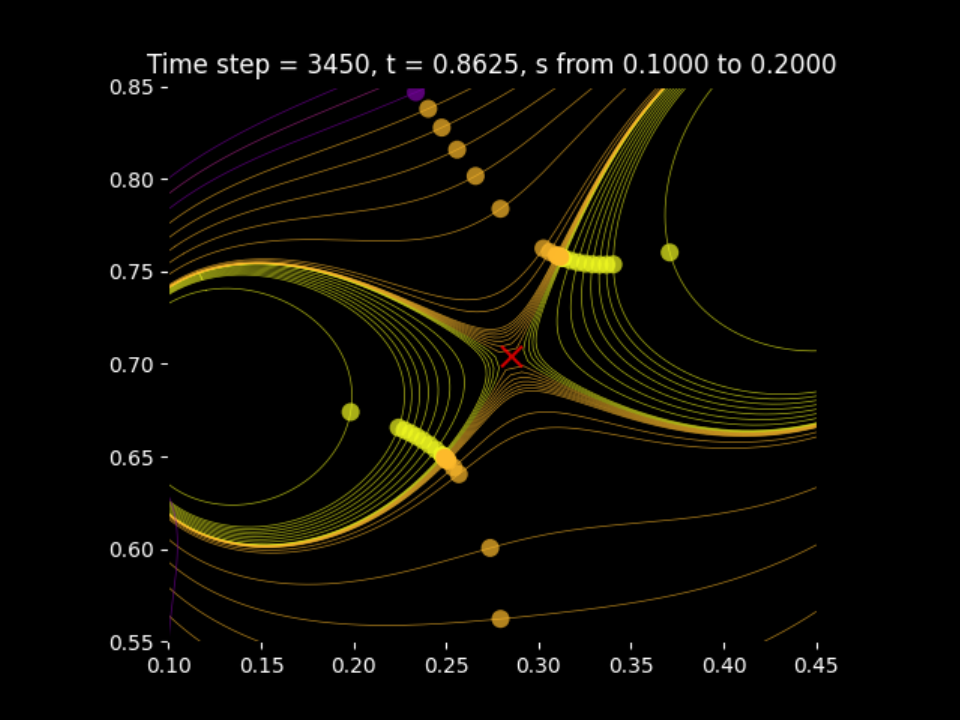}
		\includegraphics[width=0.3\textwidth]{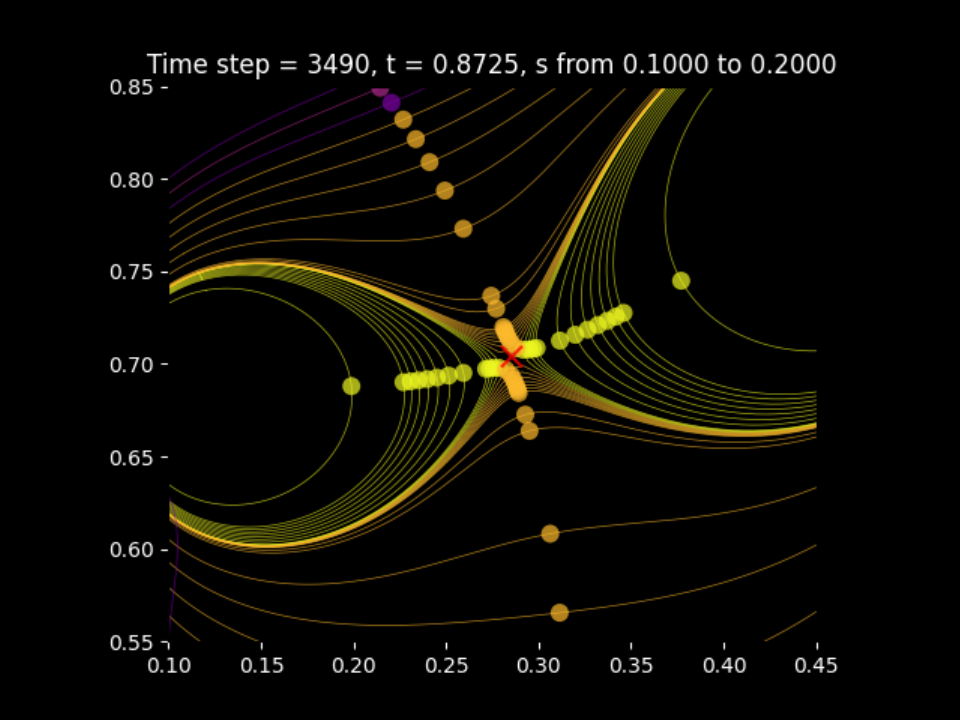}
		\includegraphics[width=0.3\textwidth]{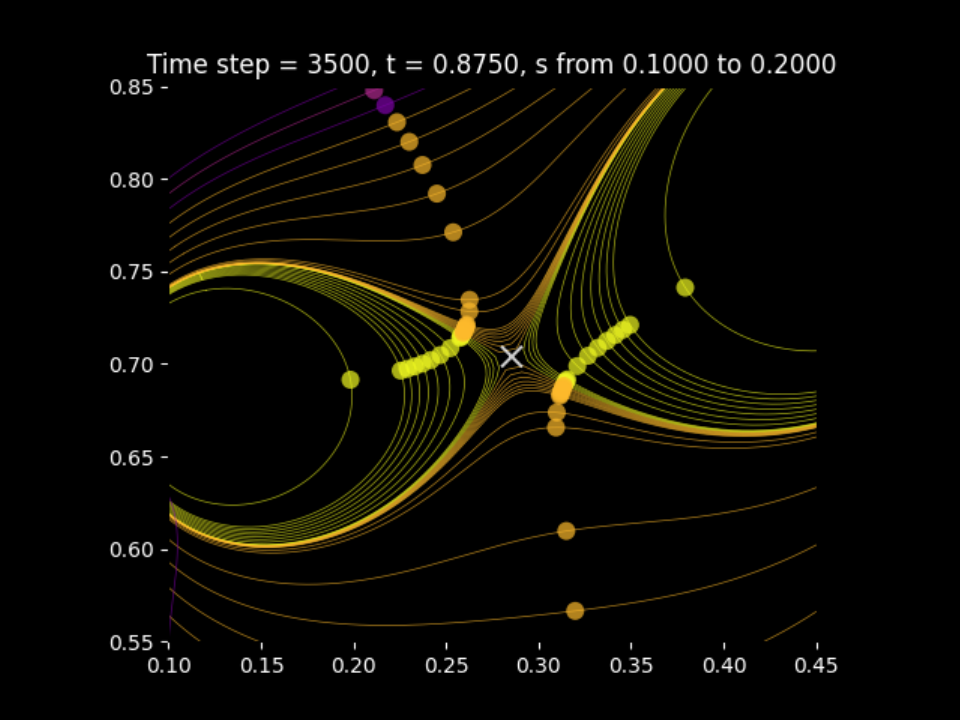}
		\caption{Collision sequence}
		\label{fig:collision sequence}
	\end{figure}

	By computing transpositions across different seeds for the Ginibre matrix 
	and using fine $s$-steps, we observed that the total number of collisions 
	is approximately $N(N-1)$ for the circular curve, and slightly larger for other curves. 
	Because this count grows rapidly with $N$, care must be taken when estimating 
	the permutations $\sigma(s_0)$ and $\sigma(s_1)$ over small intervals $[s_0,s_1]$. 
	Let $\Delta s = s_1 - s_0$ and $\Delta \sigma$ be the permutation such that 
	$\sigma(s_1) = \Delta \sigma \circ \sigma(s_0)$. Then:
	\begin{enumerate}
		\item If $\Delta s$ is small enough that each eigenvalue collision 
		in the interval involves disjoint pairs, $\Delta \sigma$ is a product of commuting transpositions.
		
		If we increase $\Delta s$ slightly and any of these pairs got repeated,
		the identical transpositions would cancel out in the difference permutation, 
		and our method would fail to detect them. So we should keep $\Delta s$ small.	

		\item Suppose that $\Delta s$ is slightly larger,
		and the collisions
		of the eigenvalues in the strip involve
		one of the eigenvalues more than once,
		but not so large that the collision multi-graph $G$
		(i.e., draw a non-directed edge connecting $(i,j)$
		for each collision between $(i, j)$)
		is a forest, but not a (sub-)matching. 
		
		We will not be able to track-down 
		the specific collisions 
		from the permutation discrepancy $\Delta\sigma$ 
		associated to each tree component of $G$ of size greater that $2$.

		For example, if there are two collisions 
		$(1,2)$, and $(2,3)$ within $[s_0,s_1]$,
		we will only see the permutation discrepancy 
		$(2,3)(1,2)=(1,3,2)$.
	
		Unless we compute the permutation at some convenient 
		intermediate values in $[s_0,s_1]$, 
		we will not be able to determine whether 
		$(1,2)$ and $(2,3)$ occurred in that order,
		or if it was $(1,3)$ followed by $(2,1)$, 
		or $(2,3)$ followed by $(3,1)$.
	
		We can, nevertheless, still 
		tell the exact number of collisions from this situation.

		\item If the collision multi-graph $G$ has a non-tree component, 
		then we will actually be missing two or more collision points, 
		and we would need to do a refinement to count properly. 

	\end{enumerate}	

	Even if we could find a collision via a full $t$-stripe domain's permutation discrepancy, 
	we would still need to use an additional method to approximate the value of $t$
	at the collision. Switching to more regular domains allow us to overcome these problems.

	\subsection{Grid search} \label{subsection:grid-search}
	
	Since the order of collisions is quadratic, 
	the full $t$-interval stripes turn out to be rather inconvenient 
	for the task of counting and locating eigenvalue collisions. 
	Instead, we choose some $m > 0$ and split the $(s,t)$ values into a grid of $mN \times mN$ sub-squares
	$$[s_i, s_{i+1}] \times [t_j, t_{j+1}], \quad i, j \leq mN.$$ 

	In this way we have much more control on the expected number of collisions 
	per sub-square. Similar to the strip case, we track the eigenvalues along the curve $\gamma_{i,j}$ 
	that loops along the edge of each square, starting at $(s_i, t_j)$, 
	continuing to $(s_{i+1}, t_j)$, $(s_{i+1}, t_{j+1})$, $(s_i, t_{j+1})$, 
	and finally back to $(s_i, t_j)$.
	
	\begin{figure}[htbp]
		\centering
		\includegraphics[width=0.6\textwidth]{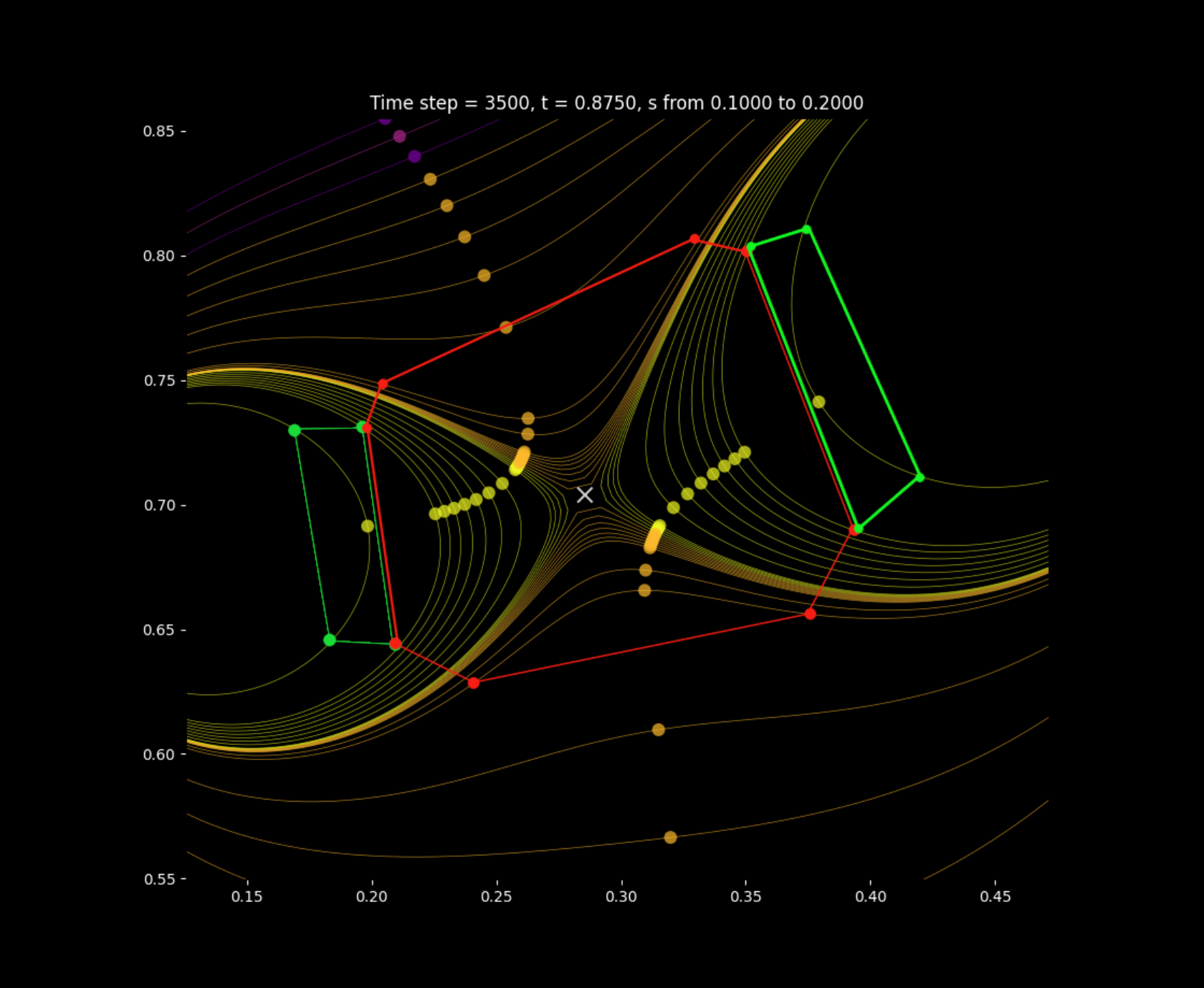}
		\caption{Moving along the edges of two $(s,t)$-squares}
		\label{fig:subsquare}
	\end{figure}

	Figure \ref{fig:subsquare} illustrates this principle: the green points mark the positions 
	of two eigenvalues at the four corners of a certain square $(s_i,t_j)$, 
	while the red points show the same two eigenvalues at the corners of the adjacent square 
	$(s_{i+1},t_j)$, where the increment in $s$ is $0.01$. 
	Notice that for the red points, instead of returning to their original positions, 
	each eigenvalue ends up where the other started, 
	after cycling around the square. 
	
	Hence, a permutation discrepancy in the ordered lists of eigenvalues
	after looping around the square, indicates a collision inside the square.
	Sometimes a collision occurs elsewhere, within the same $(s,t)$-square involving 
	some other pair of eigenvalues.
	This is not a problem, as we are keeping track of \emph{all} the eigenvalues, 
	as long as the associated transpositions are commutative.
	One can always refine so that this is the case.

	The fact that the eigenvalues crash rather violently prevents them from crashing again 
	against each other within a 
	sufficiently small $(s,t)$-time-window, in contrast with the full $t$-strip case (where two eigenvalues 
	would crash again more often, at rather distant values of $t$). 
	Hence, these square grids are much more efficient 
	for detecting eigenvalue collisions. The same rules (i) - (iii) 
	apply for the square contours  
	but situations (ii) and (iii) are much less likely to occur.
	
	\section{Collision statistics} \label{section:collision-statistics}

	For $N=10$ we computed the eigenvalue collisions for 100 trials, in several cases.
	We report that we used the seed values $1000$ to $1099$ (for reproductibility).

	\subsection {$N=10$}

	\begin{enumerate}
		\item Ginibre-to-Circle. For the complex Gaussian case, 
		all 100 trials resulted in exactly $90$ collisions. 
		At $m=5$, the grids missed 2 collisions for 4 out of the 100 seeds. 
		All $100$ found $90$ collisions for $m=6,7,8,9,10$.    
		We originally expected some randomness in the number of collisions and are surprised 
		to find this not to be the case. 
		
		Removing the trace or switching from complex Gaussian to Bernoulli Ginibre matrices 
		seems to have little effect on the collision count.
	
		Out of the 100 (traceless) Bernoulli trials, our algorithm did throw some fifteen seeds 
		where it seemed that fewer collisions were occurring. 
		After we investigated each of these cases, we noticed that they coincided exactly with the 
		initial values of Ginibre matrices with repeated eigenvalues. We even got one case
		with a triple eigenvalue and two cases with two pairs of repeated eigenvalues.
				
		We did find some interesting different behaviors for initial collisions 
		(see Appendix~\ref{appendix:repeated-eigenvalues}).
		
		Repeated eigenvalues for traceless Bernoulli Ginibre matrices quickly become unlikely 
		at slightly larger values of $N$.

		\item Ginibre-to-Circuit. Going from the circular law to the circuit curve results in
		more collisions. This time there is some variance on the results, as seen in the first histogram. 

		\begin{figure}[htbp]
			\centering
			\includegraphics[width=0.45\textwidth]{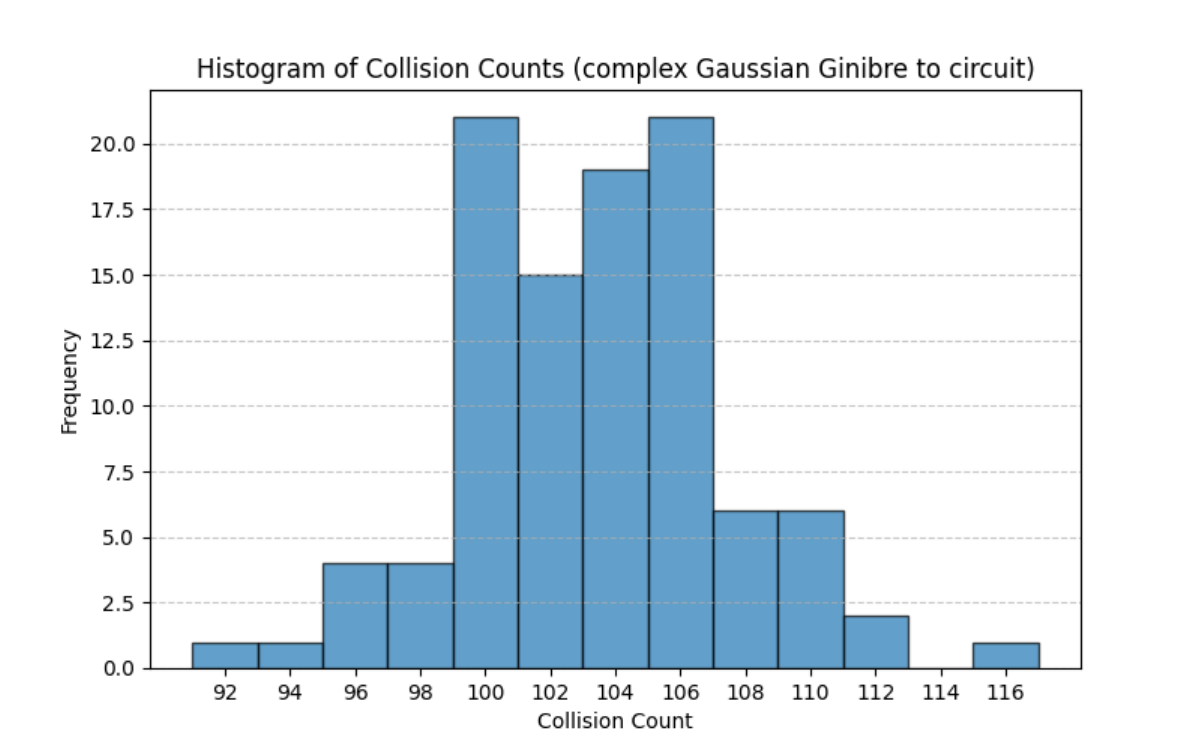}
			\includegraphics[width=0.45\textwidth]{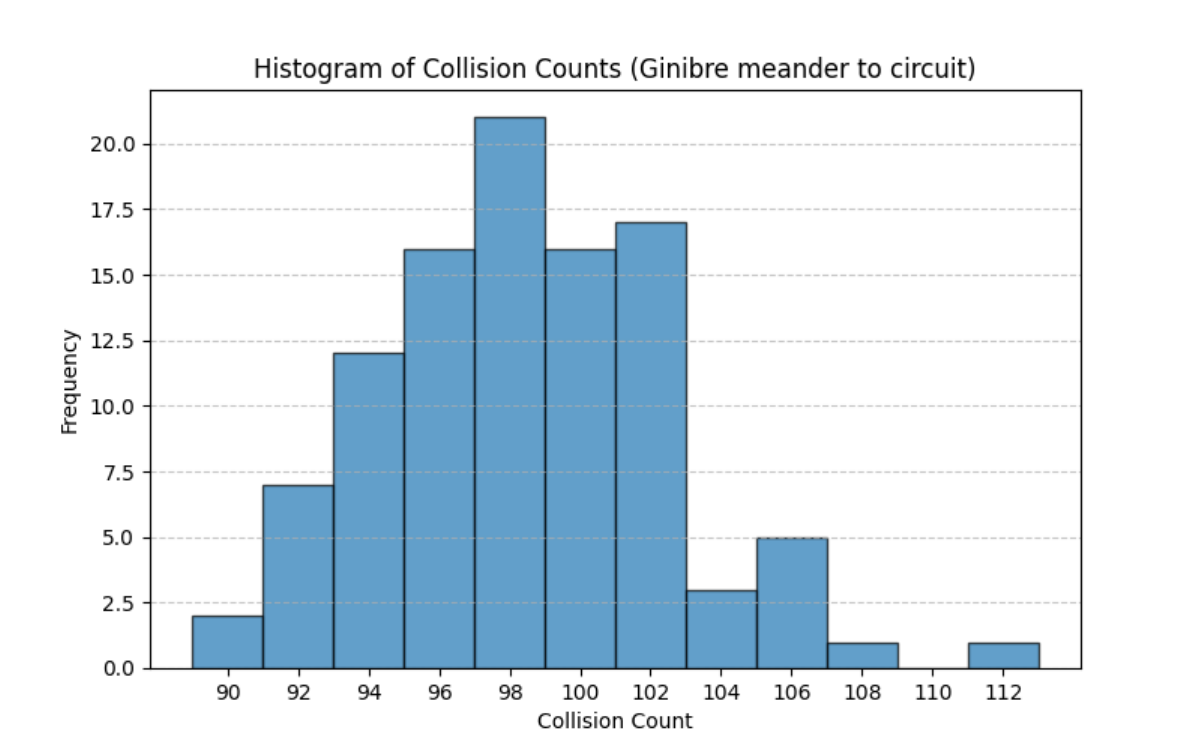}	
			\caption{Histograms of collision counts for cases: complex Gaussian Ginibre 
			to circuit curve (left) and 
			Ginibre restricted to meander to circuit curve (right)}
			\label{fig:pdf_image}
		\end{figure}

		\item Meander-to-Circuit. We tried to improve the collision count for the circuit case 
		using some free-probabilistic intuition.

		Since the circuit curve is completely contained in the unit disk
		and has $5/9$ the area of the unit disk, we may consider a Ginibre matrix of dimension $9N/5$ and
		keep only the eigenvalues that lie in the region enclosed by the circuit curve.
		After a few trials, thanks to the eigenvalue repulsion, one gets exactly $N$ points. 
		
		We put these eigenvalues back in a diagonal matrix $D$ and consider its rotation 
		by an independent, random unitary matrix $VDV^*$. 
		This puts the eigenspaces of $VDV^*$ again in general positions w.r.t. 
		the eigenspaces of the curve matrix $U$ while keeping the evenly spaced eigenvalues 
		of the Ginibre matrix restrincted to our meander shaped domain. 
	
		As seen in the second histogram, the number of collisions gets reduced by an 
		average of 5 collisions. 
		
		We wonder if there is another way (perhaps taking the curvature into consideration, 
		instead of using equidistant points) to further reduce collisions in this case.

		\item For the crossing curve we actually used $N=11$. Our current implementation of the algorithm
		has some issues with the even cases for values of $s$ close to $1$ 
		(as the points in the curve will crash at the crossing). 
		For odd $N$ this is no problem and by taking sufficiently small steps the points will cross
		alternately from both directions.

		Similar to the circuit cases, the number of collisions has non-zero variance, 
		with mean around $134$ when starting from a complex Gaussian Ginibre matrix, 
		and $125$ when starting from randomly rotated eigenvalues of a Ginibre matrix
		restricted to the opposing sectors enclosed by the crossing curve.

		\begin{figure}[htbp]
			\centering
			\includegraphics[width=0.45\textwidth]{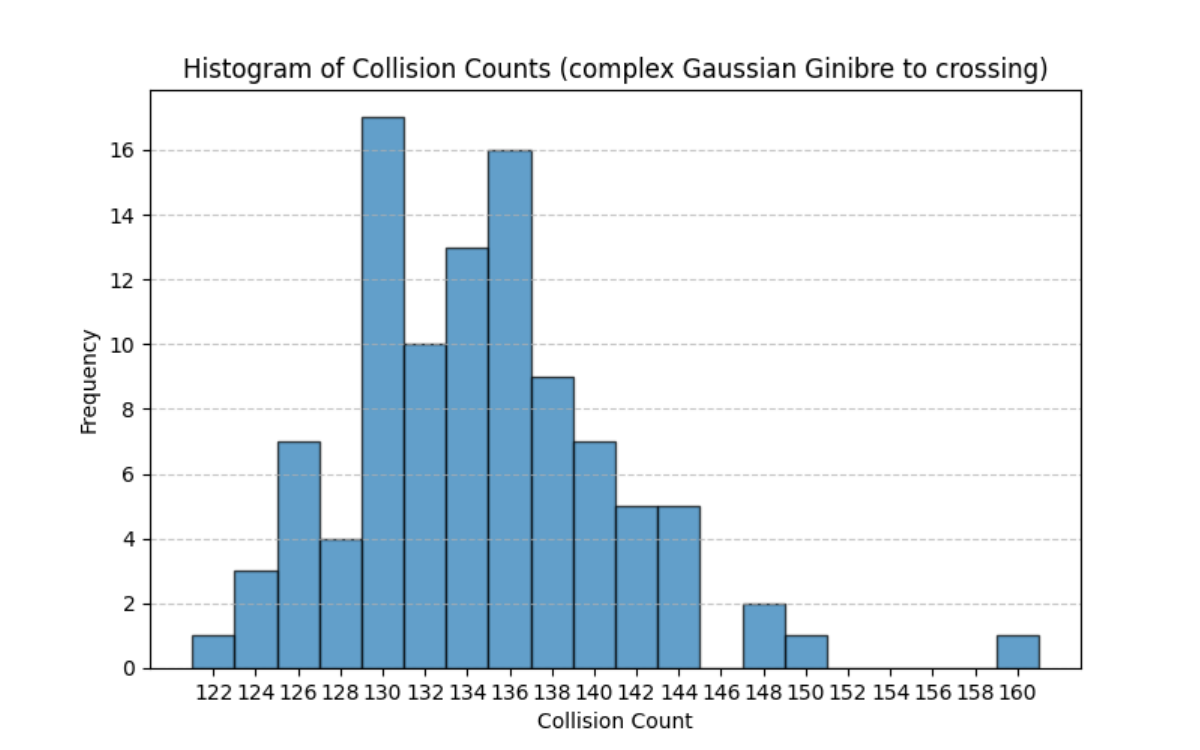}
			\includegraphics[width=0.45\textwidth]{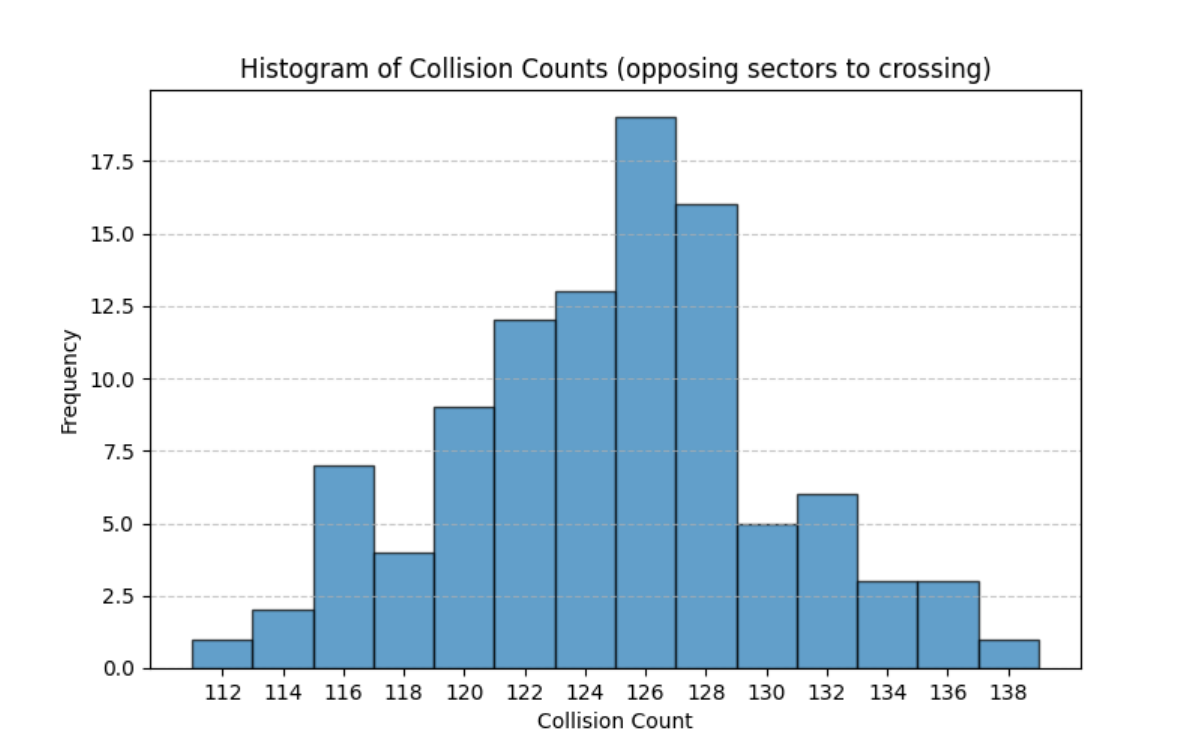}	
			\caption{Histograms of collision counts for cases: complex Gaussian Ginibre 
			to crossing curve (left) and 
			Ginibre (restricted to opposing sectors) to crossing curve (right)}
			\label{fig:pdf_image}
		\end{figure}

	\end{enumerate}		

	\subsection {Elliptic case}

	We considered a single seed (with value 2007, for reproductibility) and computed the eigenvalue collisions 
	at different values of $x$,
	replacing $C$ by the elliptic law $E_x = \cos(x) C + \sin(x) C^*$ (see Fig. \ref{fig:elliptic}) 
	and $U(t)$ by the corresponding 
	arc-length equidistant flow along the ellipse. 
	The model is analytic in the variable $x$ and hence the eigenvalues can be tracked too.
	
	As we increased $x$ from $0$ to $\pi/4$ (and $\rho$ from $0$ to $1$), 
	the number of collisions seems piece-wise constant with only a few jumps, 
	as seen in the table below:
	
	\vspace{1em}
	\begin{center}
	\begin{tabular}{l | c c c c c c c c c c c}
		$\rho$ & $0$ & $0.436$ & $0.527$ & $0.572$ & $0.6$ & $0.662$ & $0.714$ & $0.760$  & $0.8$ & $0.835$ & $0.866$ \\
		\hline
		collisions & 90 & 90 & 90 & 94 & 96 & 98 & 104 & 108 & 108 & 106 & 102 \\
	\end{tabular}
	\end{center}
	\vspace{1em}

	The collision points should evolve continuously on $x$. 
	We expect some ramification phenomena of the fibers of exceptional points 
	explaining these discrepancies on the observed collision counts.	
	
\begin{figure}[htbp]
		\centering
		\includegraphics[width=0.45\textwidth]{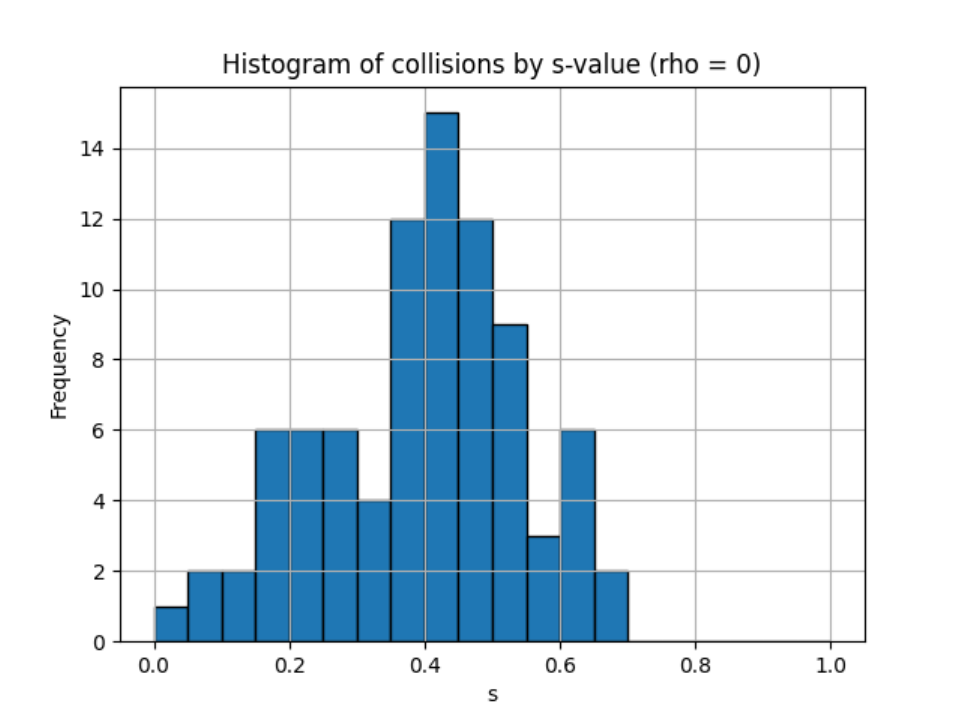}
		\includegraphics[width=0.45\textwidth]{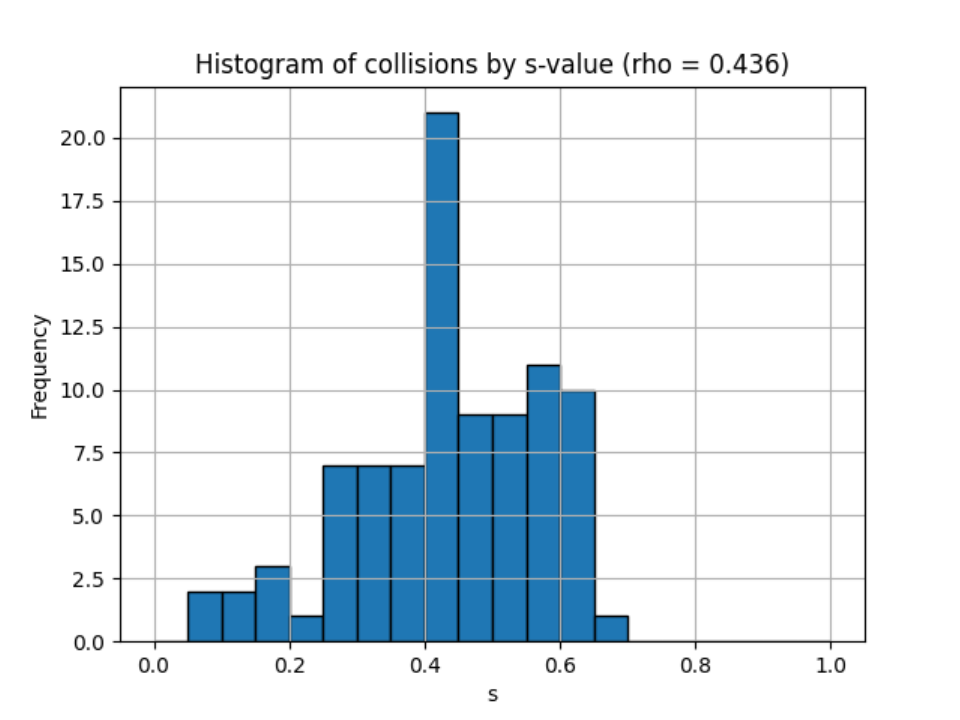}

		\includegraphics[width=0.45\textwidth]{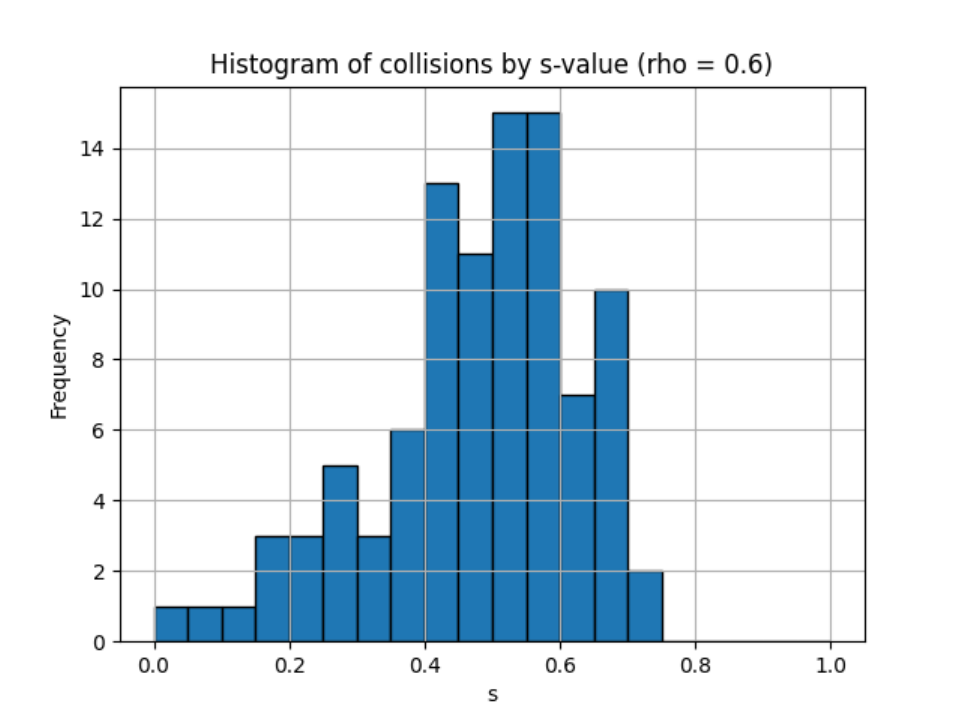}
		\includegraphics[width=0.45\textwidth]{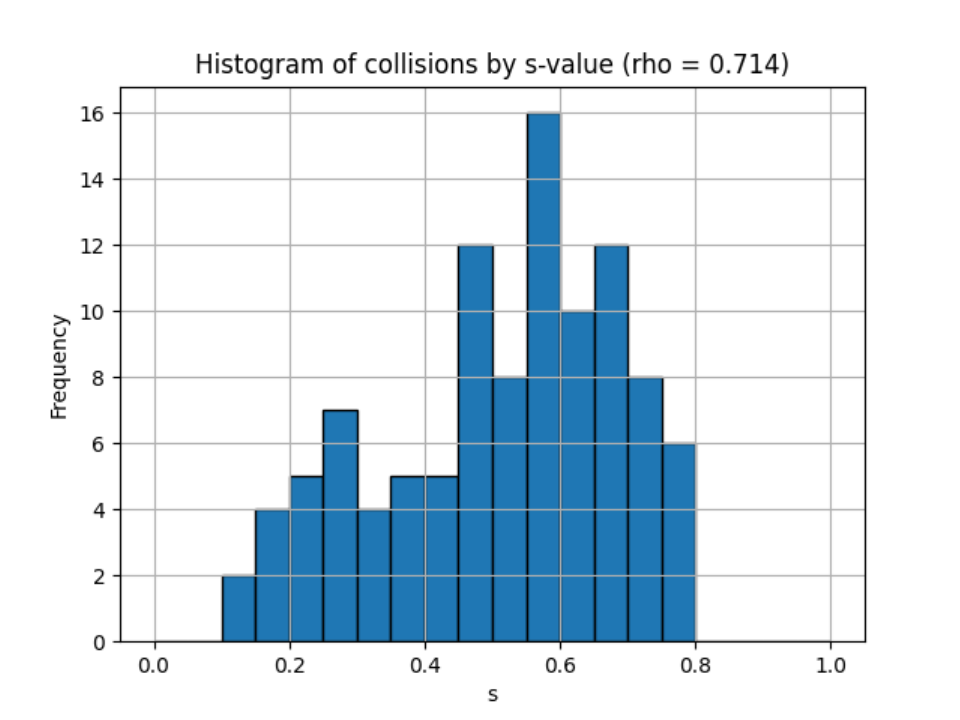}

		\includegraphics[width=0.45\textwidth]{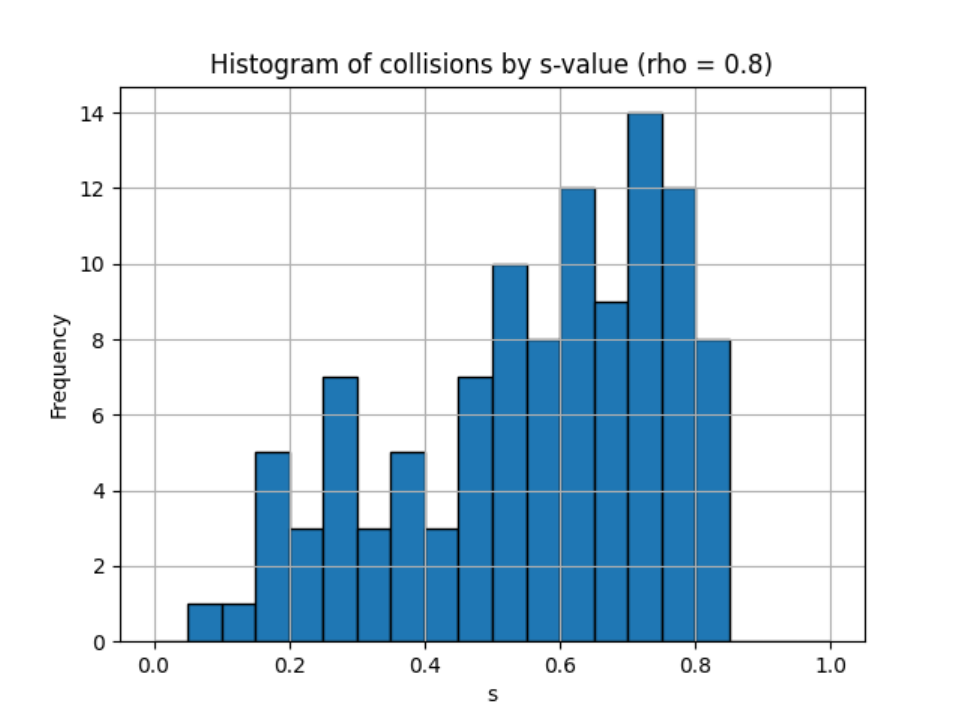}
		\includegraphics[width=0.45\textwidth]{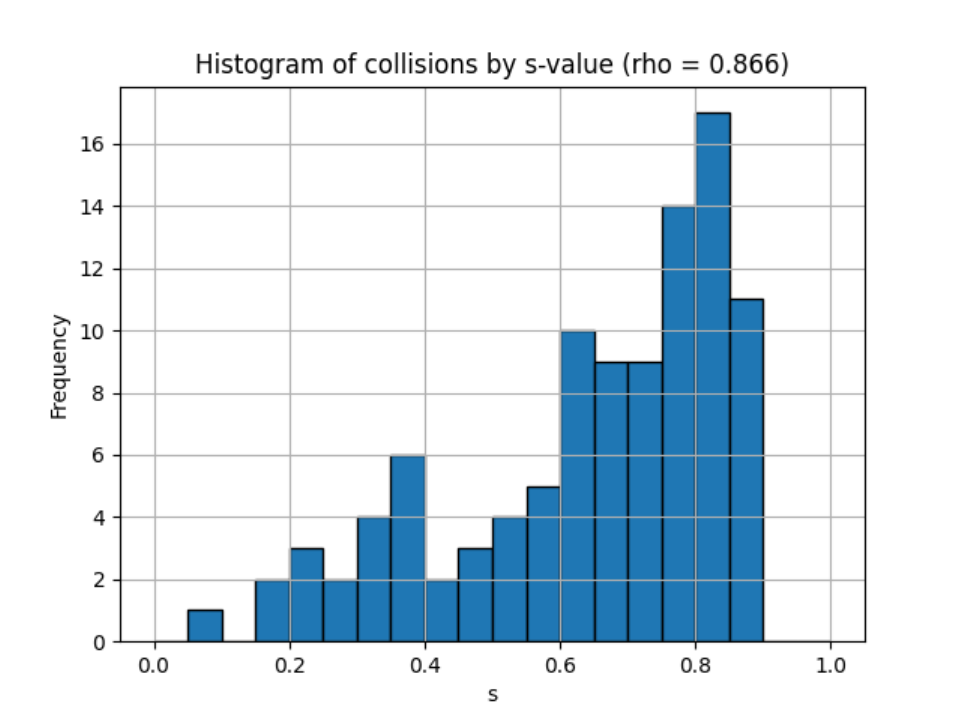}

		\caption{Histograms of collisions for different values of $\rho$}
		\label{fig:elliptic_histograms}
	\end{figure}

	We had to stop at $\rho = 0.866$. 
	For close-to-Hermitian situations, the collisions move too far to the edge
	(see Figure \ref{fig:elliptic_histograms}), 
	and our method did not seem reliable there.

	It would be interesting to implement a method to find the collisions 
	for close to Hermitian cases (and similarly, close to purely imaginary cases). 
	We would then be able to track the collision fibers, 
	while varying $x$ along the full interval $[0, 2\pi]$.

	\subsection {Larger $N$}

	We computed fewer trials for larger $N=20, 40, 100$. The general results from $N=10$ persist 
	(the circle case gives $N(N-1)$, the circuit case gives more collisions, 
	the meander trick reduces collisions).
	
	The algorithm works for large $N$, but it is time-consuming in its current state.
	An animation showing the first $60$ out of $380$ detected collisions for $N=20$ 
	can be viewed at \url{https://carlos-vargas-math.github.io/eigenvalue-collisions/animation}.	

	\section{About the algorithm} \label{section:about-the-algorithm}

	The package allows the user to select the different settings from the model
	(e.g. $N$, seed, matrix-valued distribution, curve, 
	elliptic interpolation parameter (in case curve=circle)) 
	and the algorithm (e.g. number of base s-steps, t-steps, grid resolution).

	For $m > 1$ the current algorithm divides the square $[0,1]^2$
	into $(N m)^2$ squares. 
	
	For each square of the grid the algorithm 
	takes $m_{cont}$ steps to walk along each edge of the square.
	At each step it checks if a greedy matching is possible between the current list of
	eigenvalues and those at the next step. 
	If true, it orders the eigenvalues of the next step according to the greedy matching.
	If false, it inserts an intermediate step and tries again, until it finally loops around the square.
	We use Delaunay triangulations to avoid comparing all eigenvalues against each other 
	at each step. If the tracking fails the square is reported as 'unprocessed'.
	All our statistics from Section~\ref{section:collision-statistics} were computed
	at grid resolutions with no 'unprocessed' squares.

	After going around the square the algorithm compares the lists of eigenvalues
	and computes the permutation $\sigma$ that relates them.
	If $\sigma$ is the identity the algorithm proceeds to the next square of the grid.
	Otherwise, the algorithm will process this square to approximate the values 
	of $(s,t, \lambda)$ at the collisions:

	(i) if $\sigma$ is a product of commutative transpositions, 
	the algorithm checks each pair of eigenvalues $(i,j)$ along the square
	using $m_{fine}$ steps
	and finds an approximate for $(s,t)$ and $\lambda$  
	using weighted convex combinations of the minima for their distance 
	at the square sides.

	(ii) if $\sigma$ is not a product of commutative transpositions, 
	our current algorithm increases the collision count appropriately, 
	but does not subdivide the square inductively
	to proceed as in (i) and collect the collisions.  
	We currently simply use a sufficiently large $m_{grid}$ so that this case does not occur 
	(e.g. $m_{grid}= 20$ for $N=10$).
	We will try (harder) to implement a more efficient divide-and-conquer method to improve
	this part of the algorithm. 
	Our first efforts on this only arrived at slower algorithms, 
	due to not sufficiently careful implementations.
	
	One important current bottleneck of the algorithm 
	is the linear algebraic calculation of the eigenvalues.
	Since we are dealing with matrices that are very close from each other, 
	a more specialized, explicitly pivoted method could 
	help reduce computing time at this step.
	In general, it is hard to deal with much larger matrices with the current implementation.
	The algorithm works, but it quickly becomes slow as $N$ increases. 
	
	\section{Concluding remarks and future work} \label{section:concluding-remarks}

	We are quite surprised by the invariance of the number of collisions in the process.
	There could be probably be a reasonable explanation for this by means of topological-algebraic 
	and/or linear-algebraic arguments.
	In particular, we are completely ignoring the eigenvectors in our discussions. 

	Along this note, we have pointed out a few aspects of the implementation 
	that we will gradually improve. 
	The project is publicly available and will gladly accept contributions from external users. 
	It would be interesting to find a fixed-point method type of solution
	(in particular, to replace our current convex approximation, when a collision is detected). 
	For the visualization of the eigenvalue tracks, however, 
	there seems to be no escape from our approach, 
	which requires to compute many lists of eigenvalues. 
	The methods are parallelizable.

	\bibliographystyle{plainnat}	
	\bibliography{refs}

	\appendix

	\newpage
	
	\begin{figure}[htbp]
		\centering
		\includegraphics[width=0.45\textwidth]{figures/N=10Circle00to01.pdf}
		\includegraphics[width=0.45\textwidth]{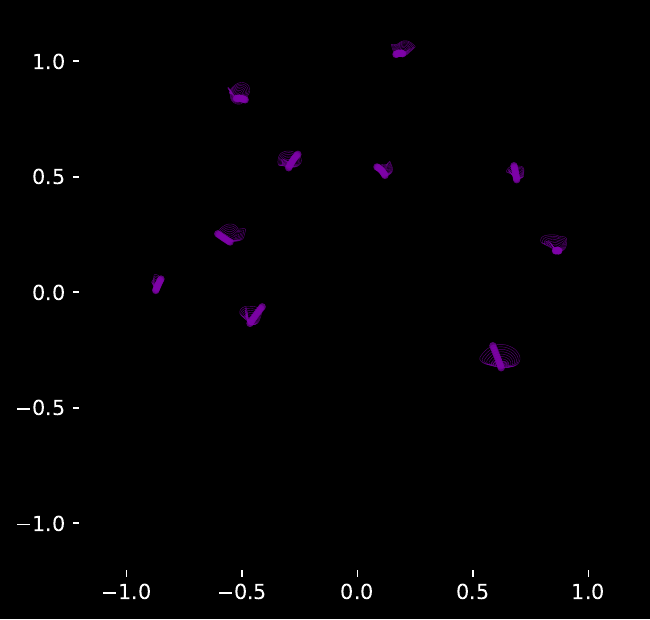}

		\includegraphics[width=0.45\textwidth]{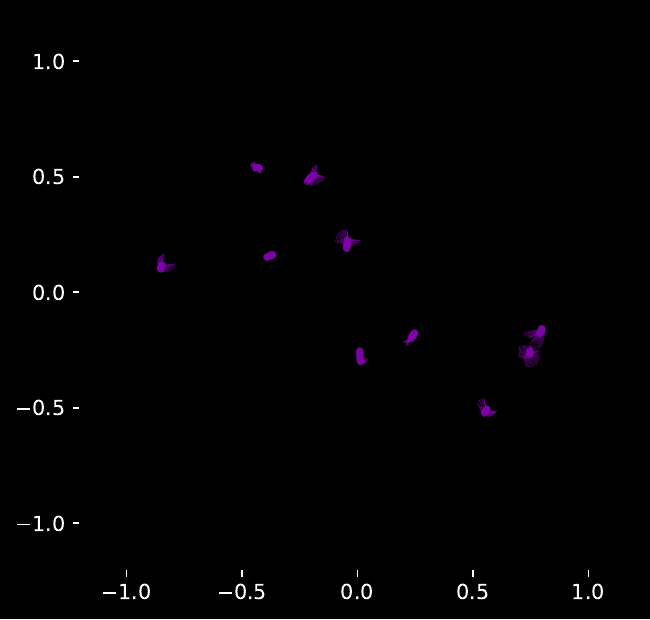}
		\includegraphics[width=0.45\textwidth]{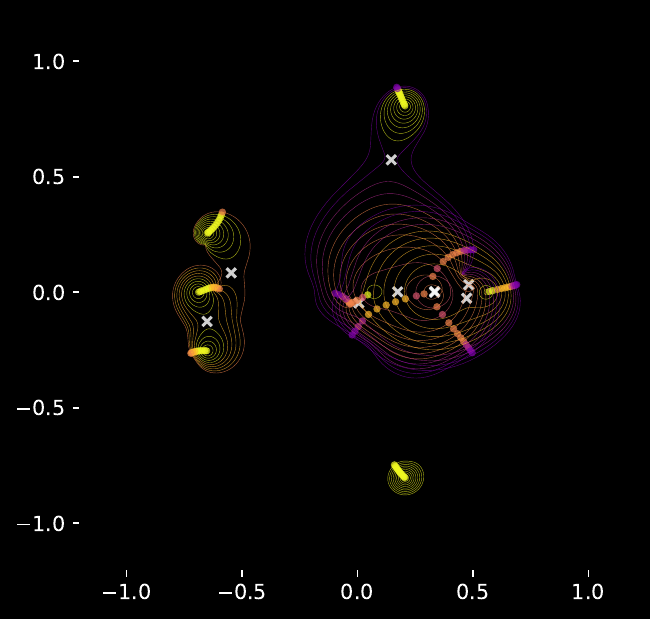}

		\caption{Eigenvalue trajectories for $s= 0.0, 0.01, \dots , 0.09, 0.10$ }
		\label{fig:pdf_image}
	\end{figure}

	\section{Collisions for $N=10$} \label{appendix:several-curves}

	We include some high resolution pictures for $N=10$ of the eigenvalue tracks and collisions for the 
	cases: 
	
	\begin{itemize}
		\item Complex-Gaussian to circle (upper-left)

		\item Ginibre Meander to circuit (upper-right),

		\item Opposing sectors to crossing (lower-left)

		\item 	Traceless Bernoulli-Ginibre to circle, starting from a triple eigenvalue (lower-right). 
		Notice the strong shock-wave/triple-helix produced by the initial triple collision
	
		In Appendix C we include pictures of the initial $s$-windows for other examples 
		of traceless Bernoulli Ginibre matrices with repeated eigenvalues.

	\end{itemize}
	 
	\newpage

	\begin{figure}[htbp]
		\centering
		\includegraphics[width=0.45\textwidth]{figures/N=10Circle01to02.pdf}
		\includegraphics[width=0.45\textwidth]{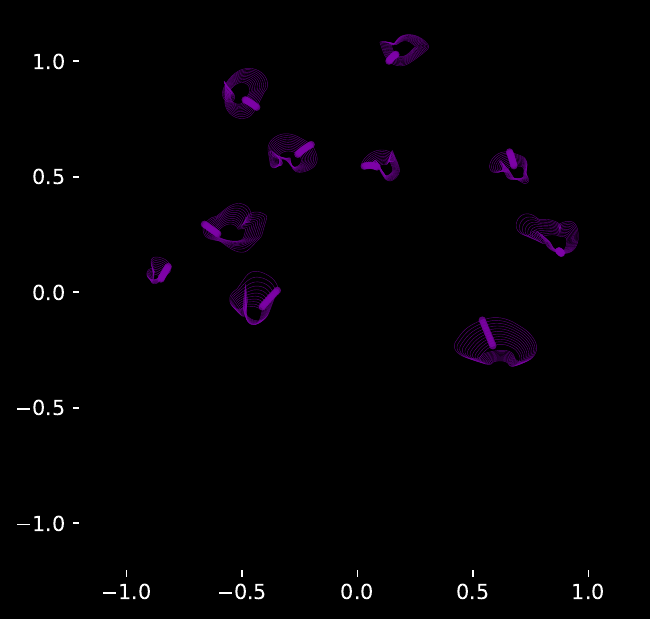}

		\includegraphics[width=0.45\textwidth]{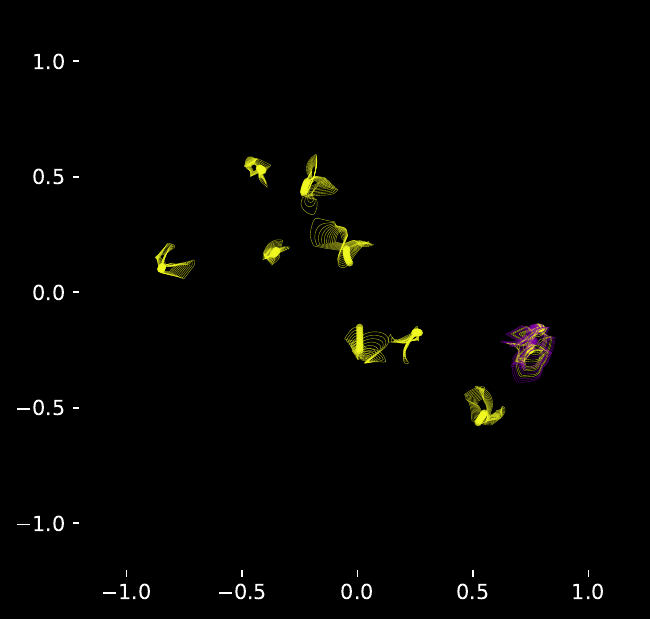}
		\includegraphics[width=0.45\textwidth]{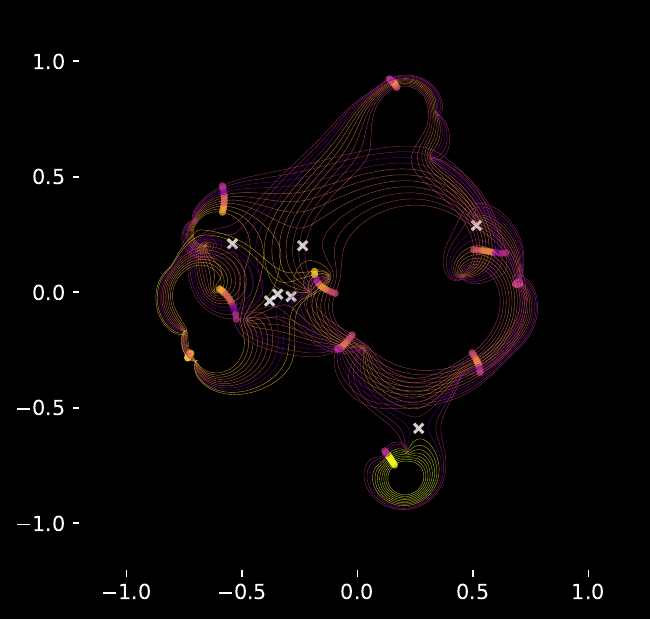}
		\caption{Eigenvalue trajectories for $s= 0.10, 0.11, \dots , 0.19, 0.20$ }
		\label{fig:pdf_image}
	\end{figure}

	\newpage

	\begin{figure}[htbp]
		\centering
		\includegraphics[width=0.45\textwidth]{figures/N=10Circle02to03.pdf}
		\includegraphics[width=0.45\textwidth]{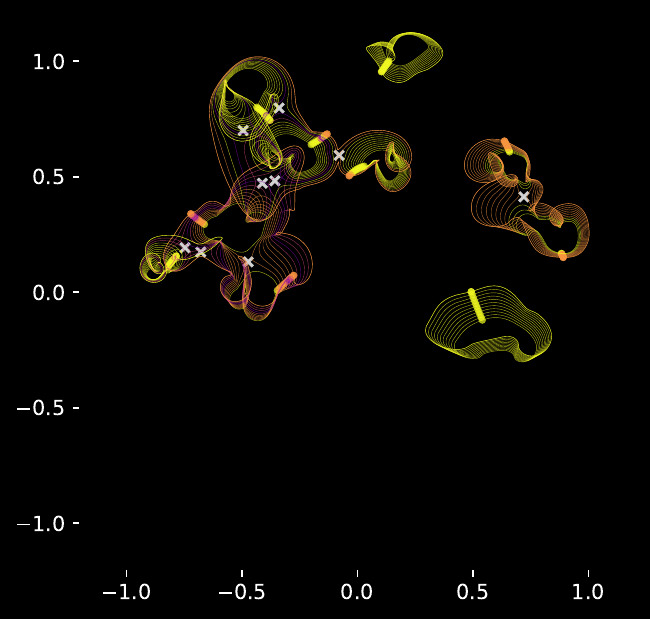}

		\includegraphics[width=0.45\textwidth]{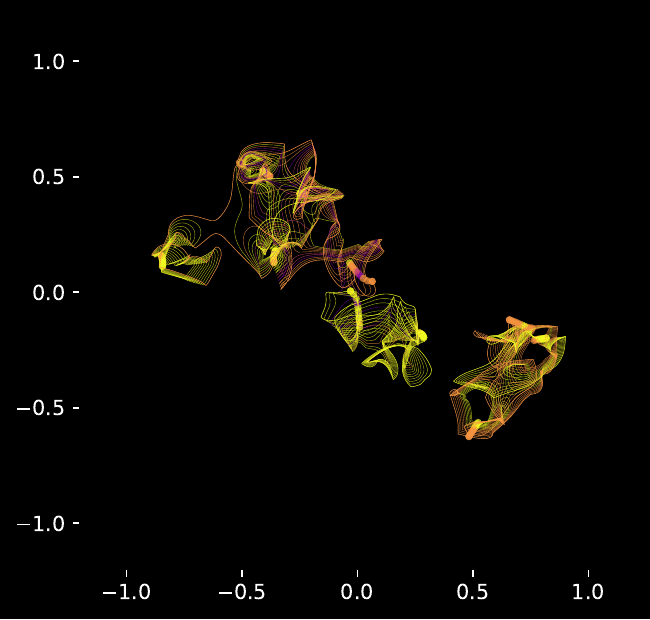}
		\includegraphics[width=0.45\textwidth]{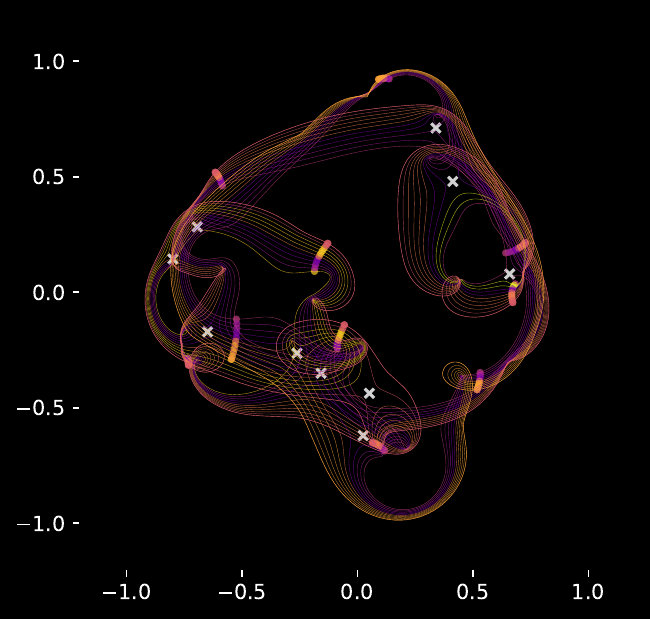}
		\caption{Eigenvalue trajectories for $s= 0.20, 0.21, \dots , 0.29, 0.30$ }
		\label{fig:pdf_image}
	\end{figure}

	\newpage

	\begin{figure}[htbp]
		\centering
		\includegraphics[width=0.45\textwidth]{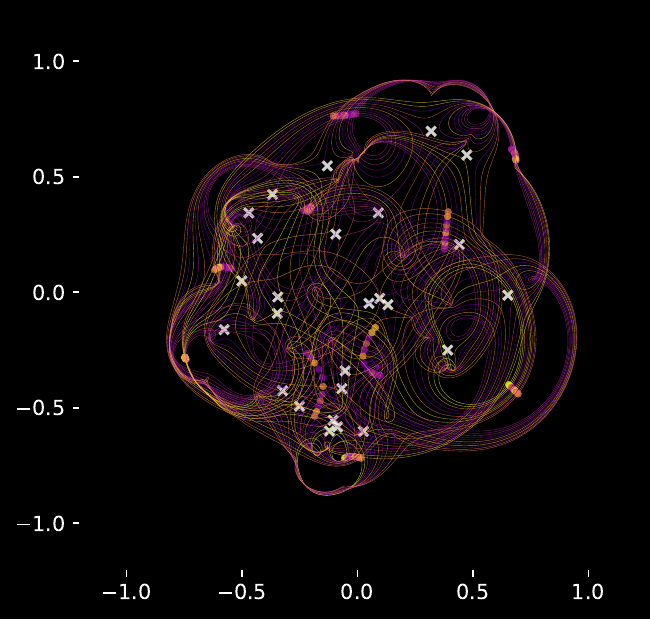}
		\includegraphics[width=0.45\textwidth]{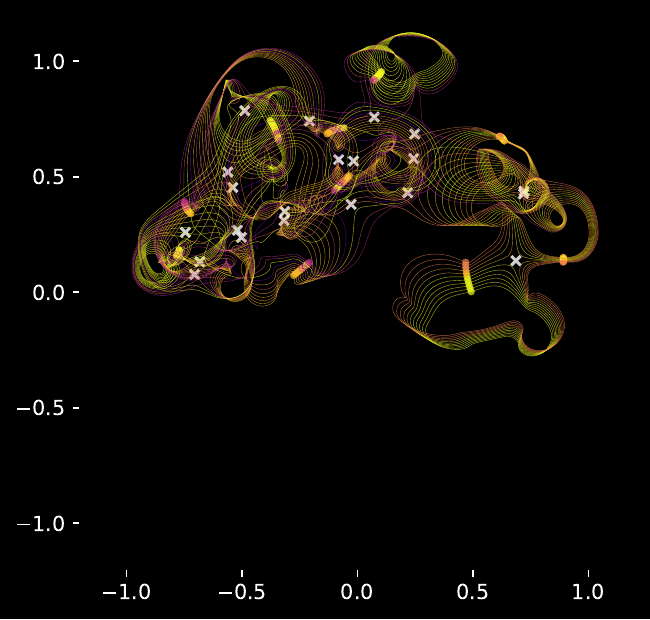}

		\includegraphics[width=0.45\textwidth]{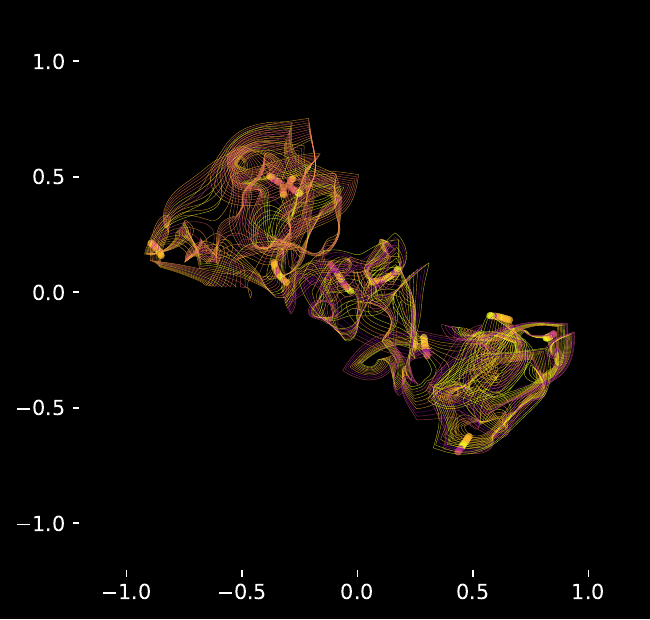}
		\includegraphics[width=0.45\textwidth]{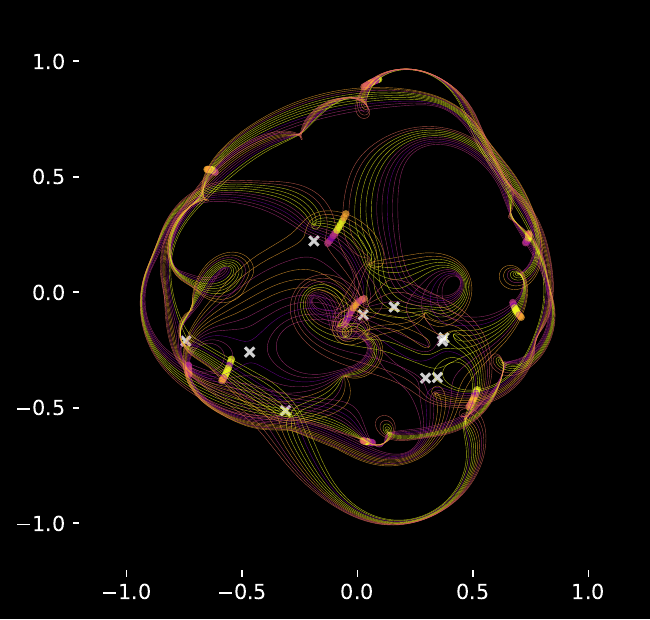}
		\caption{Eigenvalue trajectories for $s= 0.30, 0.31, \dots , 0.39, 0.40$ }
		\label{fig:pdf_image}
	\end{figure}

	\newpage

	\begin{figure}[htbp]
		\centering
		\includegraphics[width=0.45\textwidth]{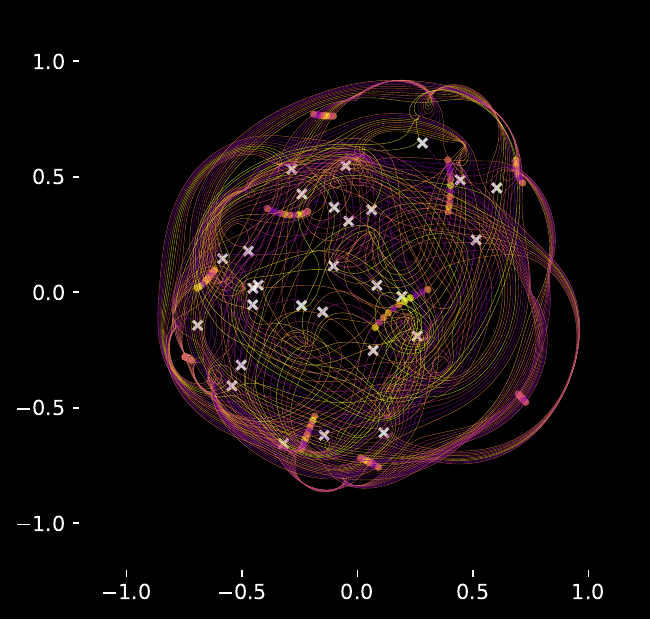}
		\includegraphics[width=0.45\textwidth]{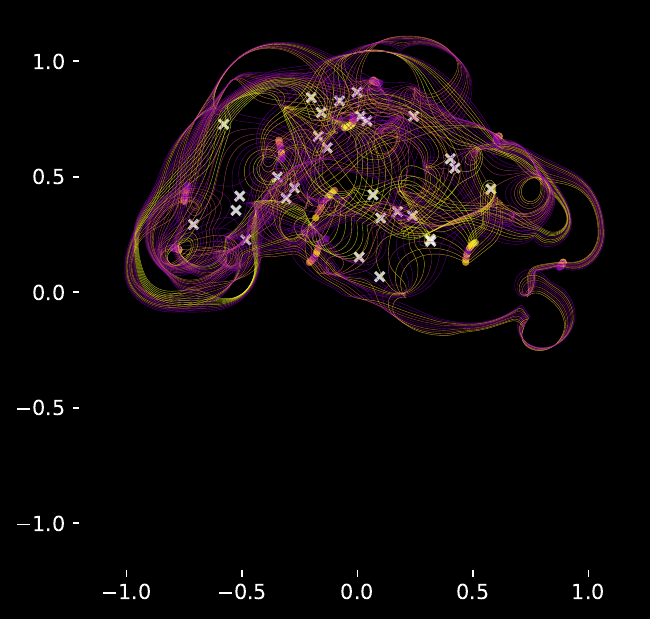}

		\includegraphics[width=0.45\textwidth]{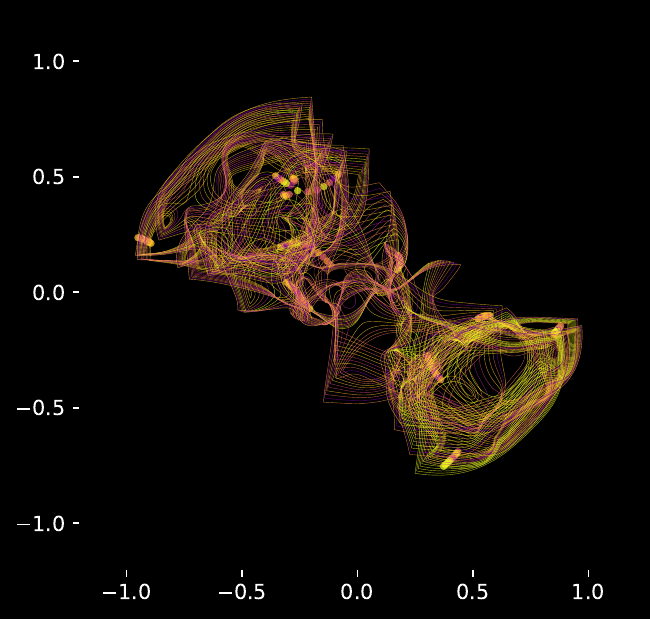}
		\includegraphics[width=0.45\textwidth]{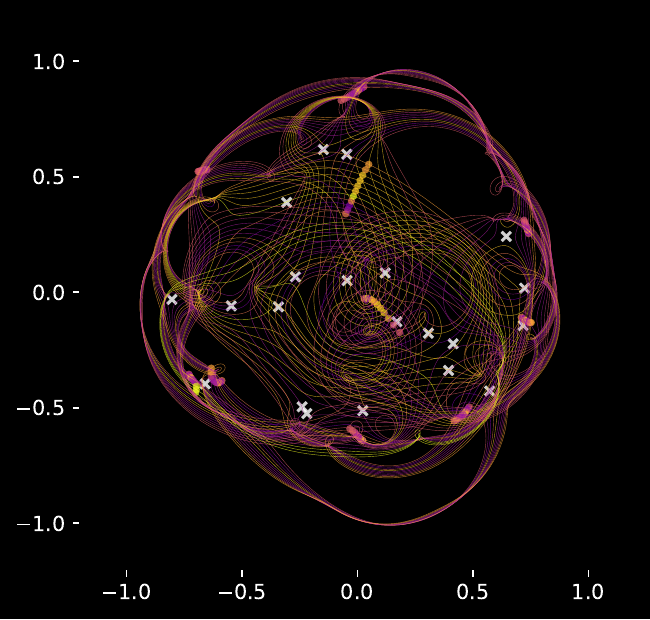}
		\caption{Eigenvalue trajectories for $s= 0.40, 0.41, \dots , 0.49, 0.50$ }
		\label{fig:pdf_image}
	\end{figure}

	\newpage

	\begin{figure}[htbp]
		\centering
		\includegraphics[width=0.45\textwidth]{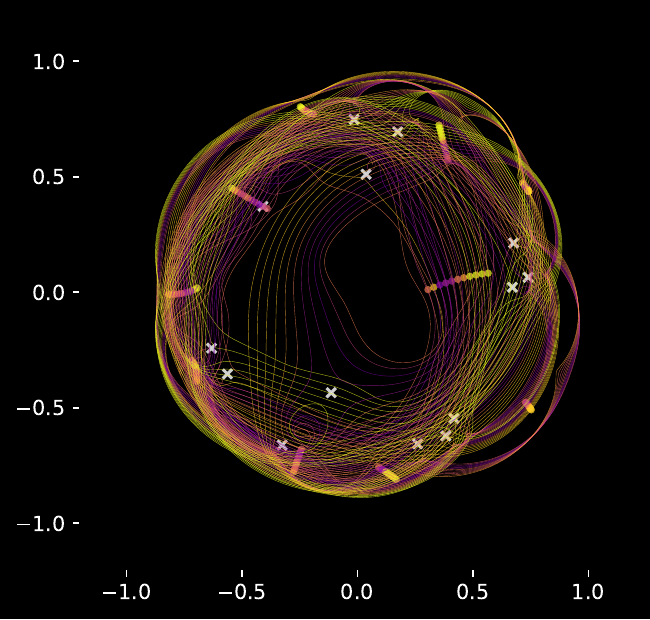}
		\includegraphics[width=0.45\textwidth]{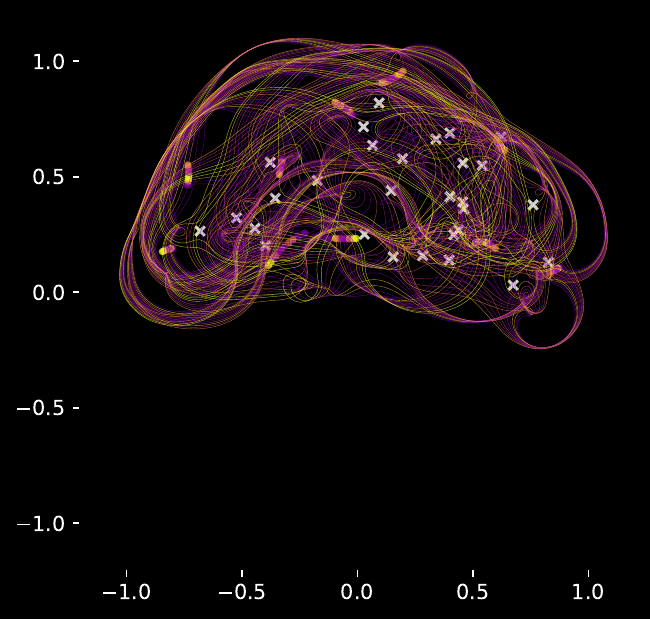}

		\includegraphics[width=0.45\textwidth]{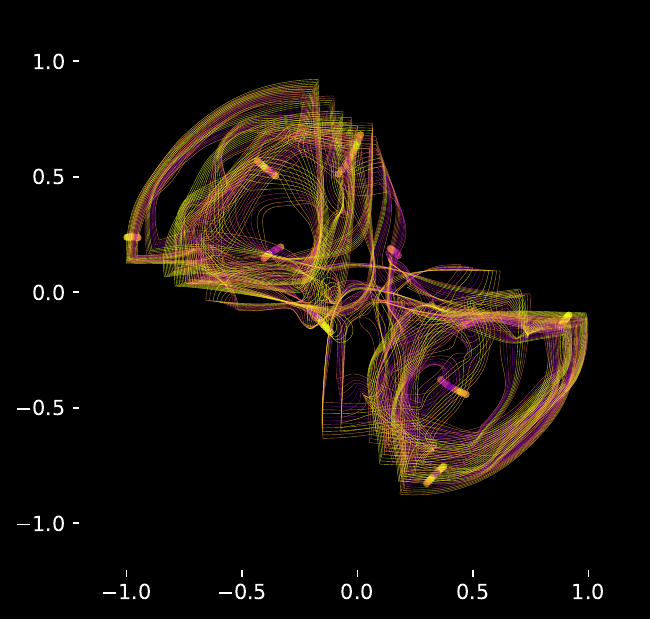}
		\includegraphics[width=0.45\textwidth]{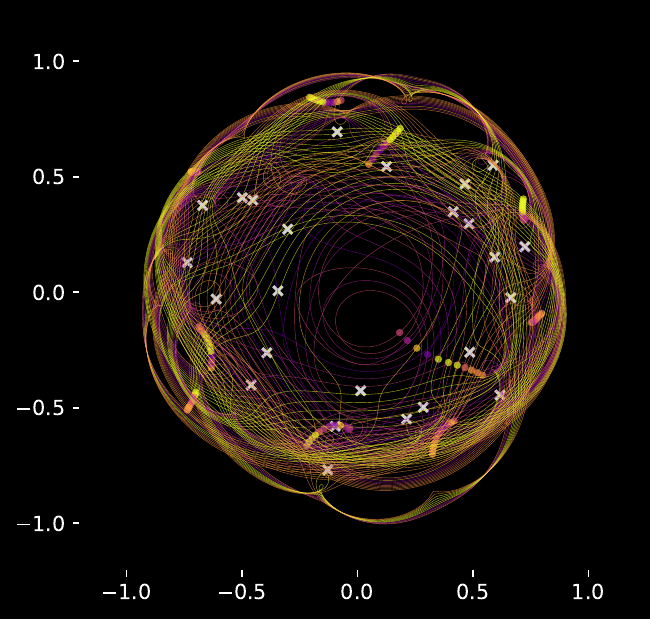}
		\caption{Eigenvalue trajectories for $s= 0.50, 0.51, \dots , 0.59, 0.60$ }
		\label{fig:pdf_image}
	\end{figure}

	\newpage

	\begin{figure}[htbp]
		\centering
		\includegraphics[width=0.45\textwidth]{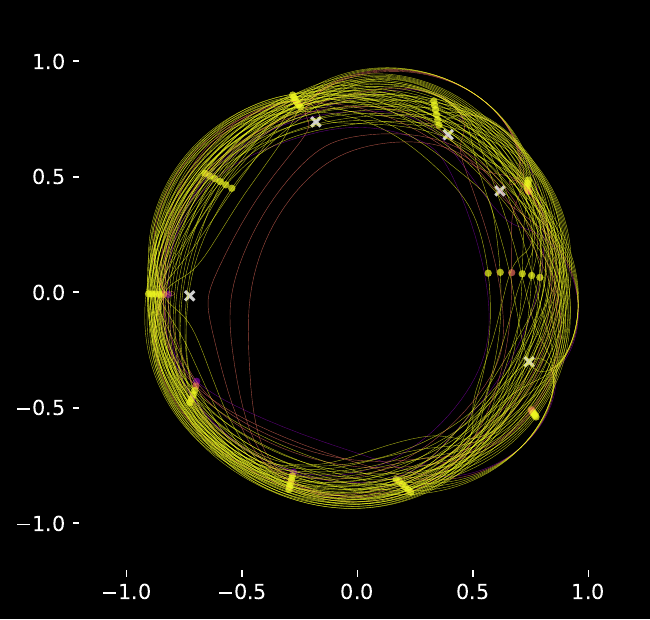}
		\includegraphics[width=0.45\textwidth]{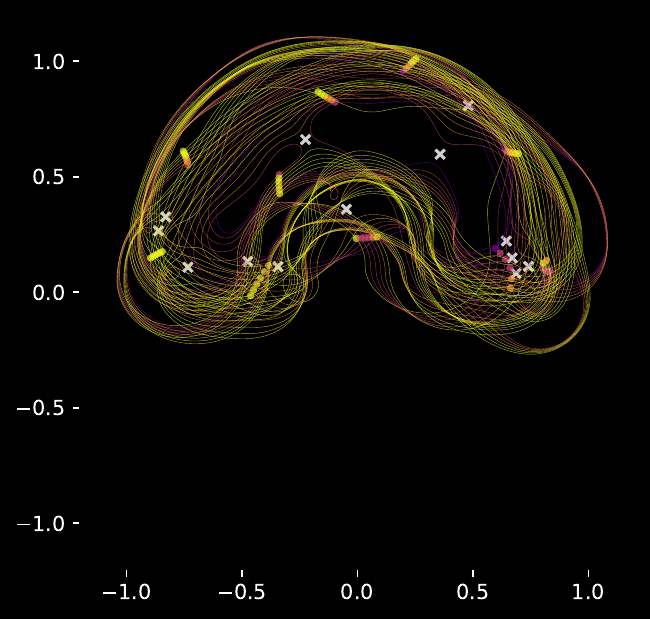}

		\includegraphics[width=0.45\textwidth]{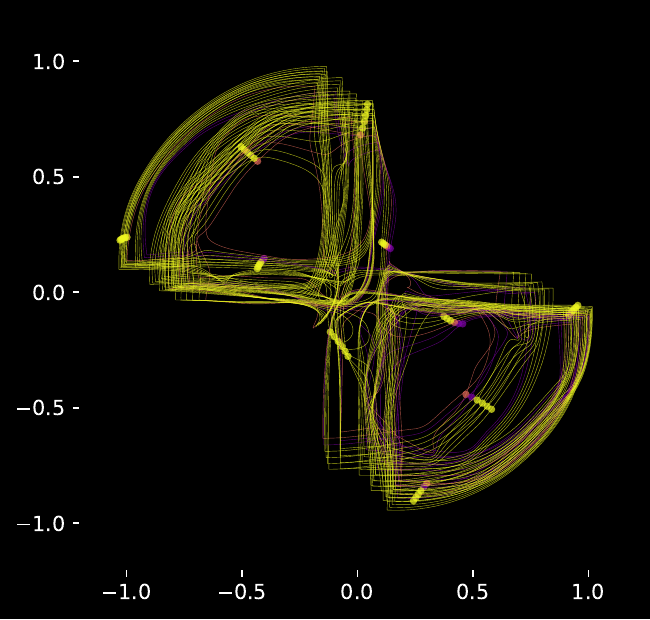}
		\includegraphics[width=0.45\textwidth]{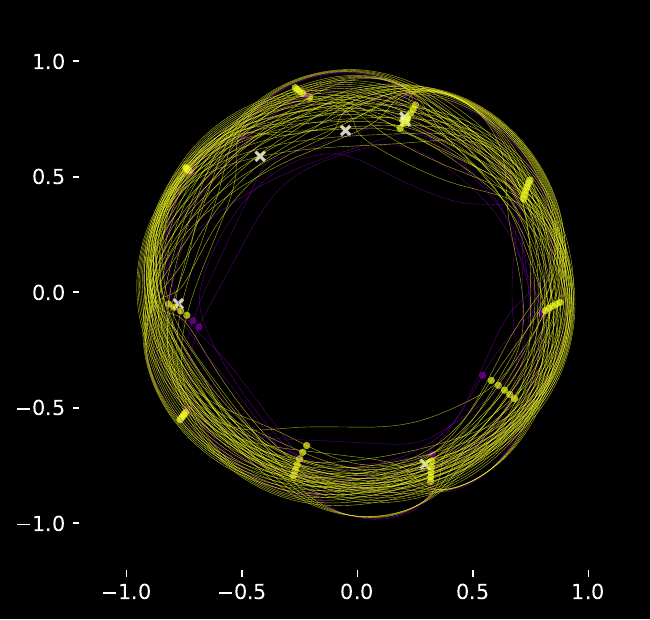}
		\caption{Eigenvalue trajectories for $s= 0.60, 0.62, \dots , 0.68, 0.70$ }
		\label{fig:pdf_image}
	\end{figure}

	\newpage
	\begin{figure}[htbp]
		\centering
		\includegraphics[width=0.45\textwidth]{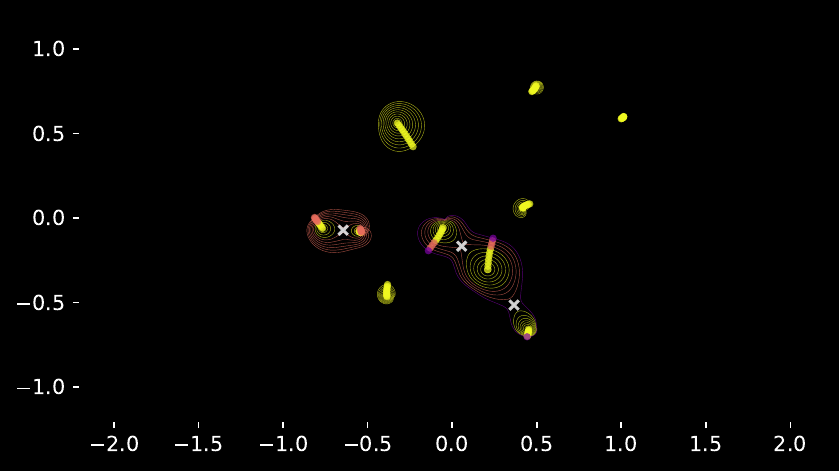}
		\includegraphics[width=0.45\textwidth]{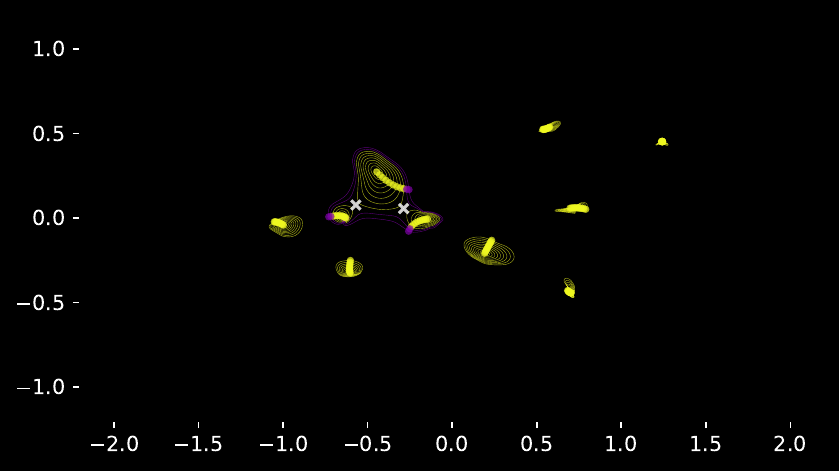}

		\includegraphics[width=0.45\textwidth]{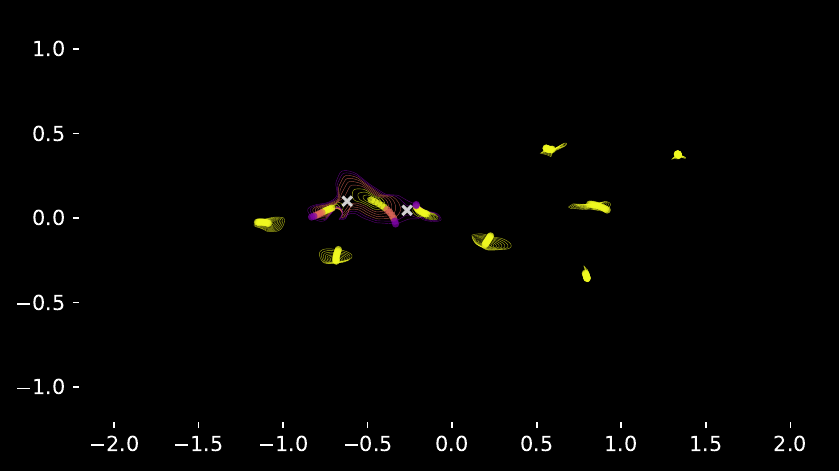}
		\includegraphics[width=0.45\textwidth]{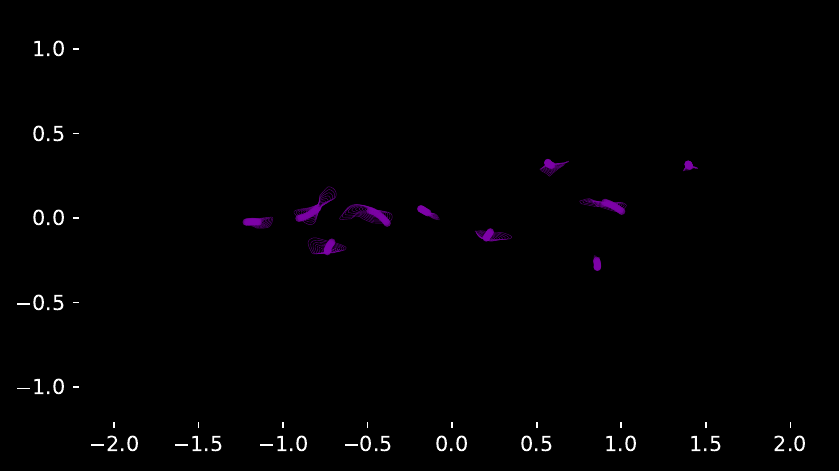}

		\includegraphics[width=0.45\textwidth]{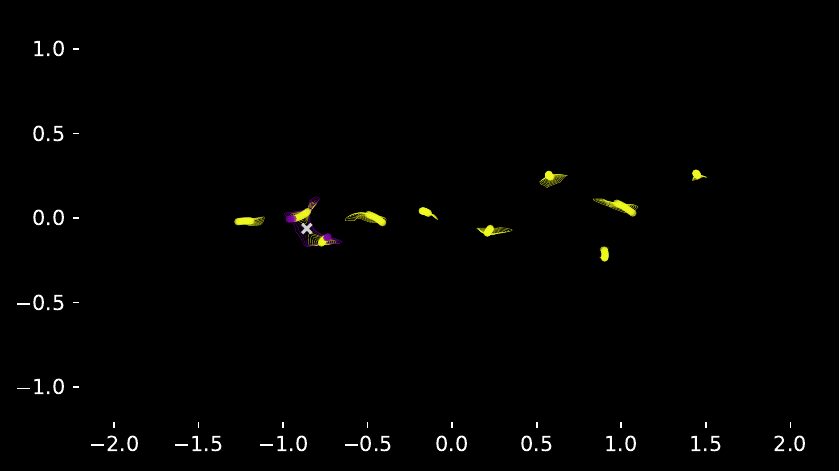}
		\includegraphics[width=0.45\textwidth]{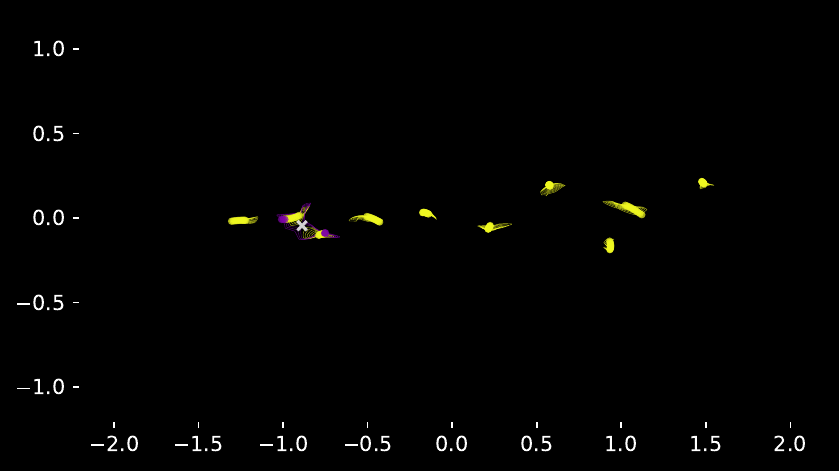}

		\caption{Eigenvalue trajectories for $s= 0.0, 0.01, \dots , 0.09, 0.10$ }
		\label{fig:pdf_image}
	\end{figure}

	\section{Elliptic cases} \label{appendix:elliptic-cases}

	We present the evolution of the eigenvalue tracks using the same seed (2007), with
	longer semi-axis $1+\rho$, for $\rho = 0, 0.436, 0.6, 0.714, 0.8, 0.866.$

	\newpage

	\begin{figure}[htbp]
		\centering
		\includegraphics[width=0.45\textwidth]{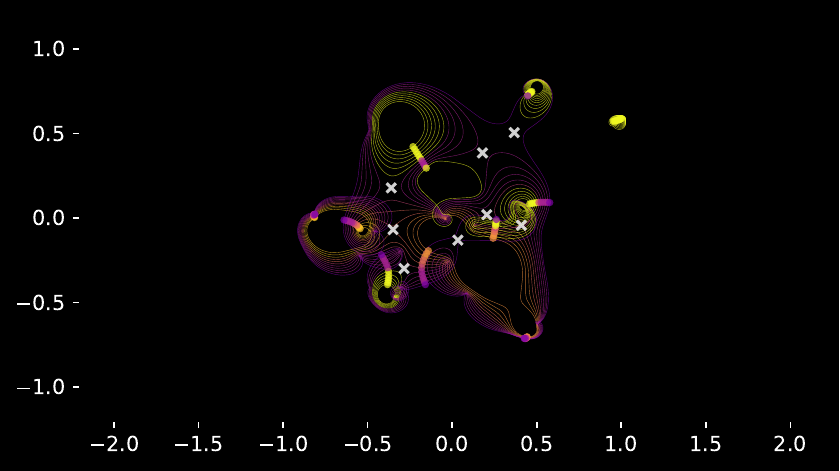}
		\includegraphics[width=0.45\textwidth]{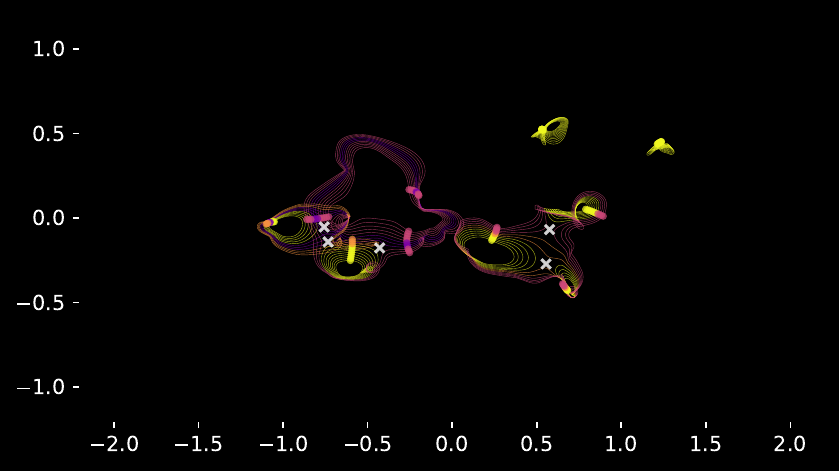}

		\includegraphics[width=0.45\textwidth]{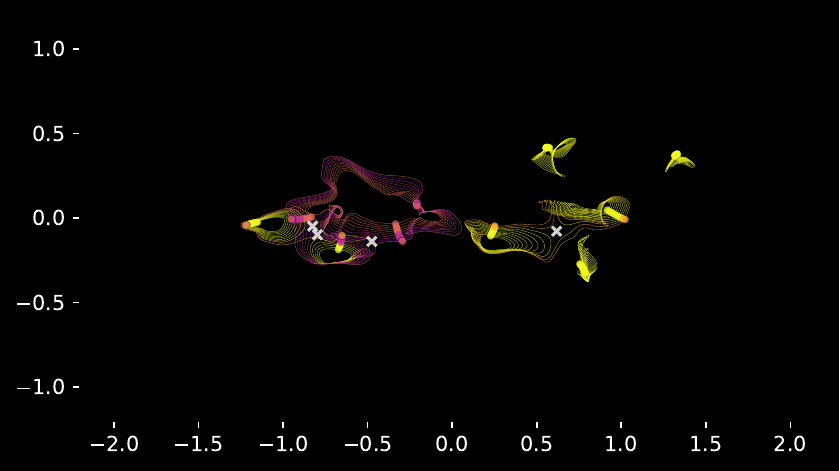}
		\includegraphics[width=0.45\textwidth]{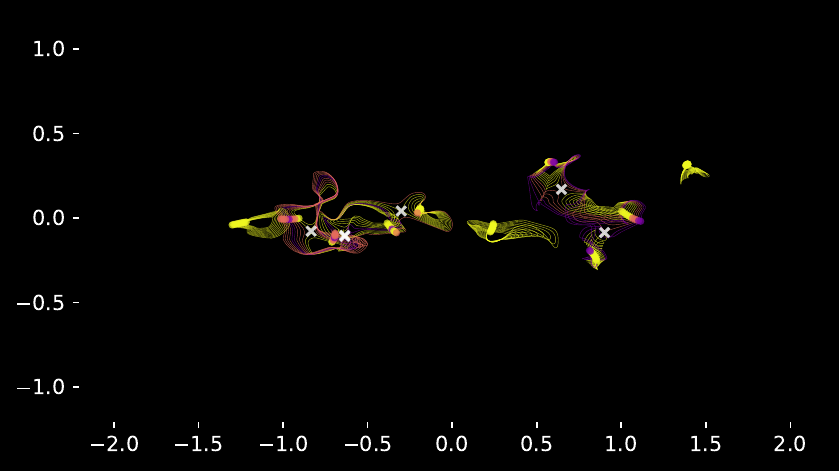}

		\includegraphics[width=0.45\textwidth]{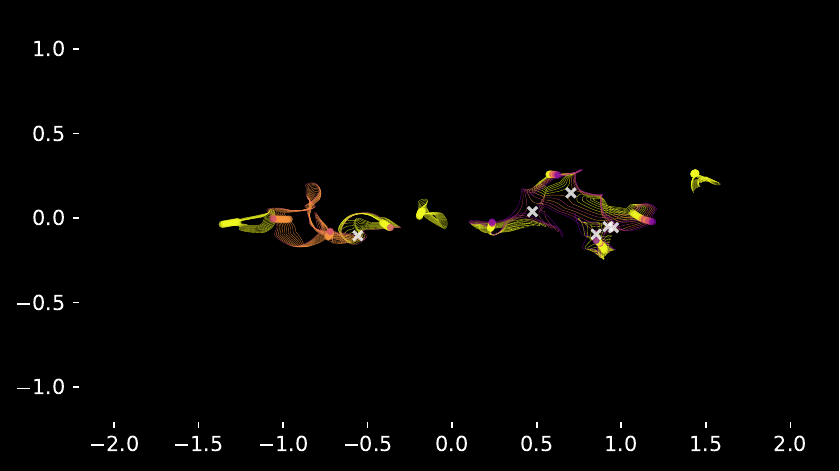}
		\includegraphics[width=0.45\textwidth]{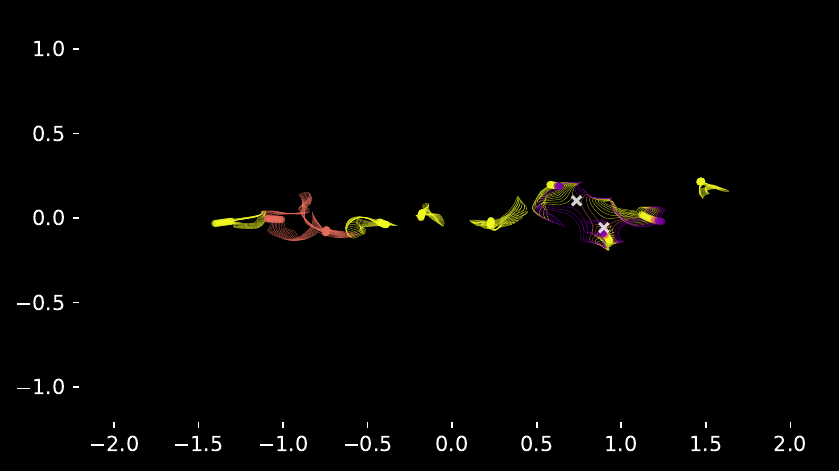}

		\caption{Eigenvalue trajectories for $s= 0.10, 0.11, \dots , 0.19, 0.20$ }
		\label{fig:pdf_image}
	\end{figure}

	\newpage

	\begin{figure}[htbp]
		\centering
		\includegraphics[width=0.45\textwidth]{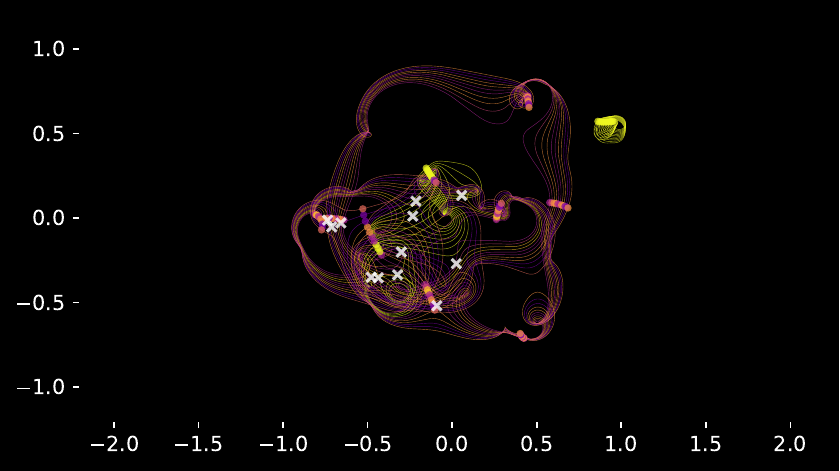}
		\includegraphics[width=0.45\textwidth]{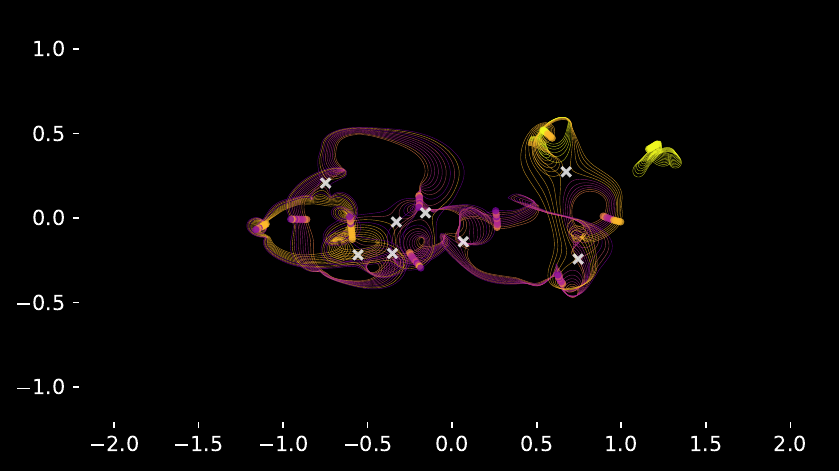}

		\includegraphics[width=0.45\textwidth]{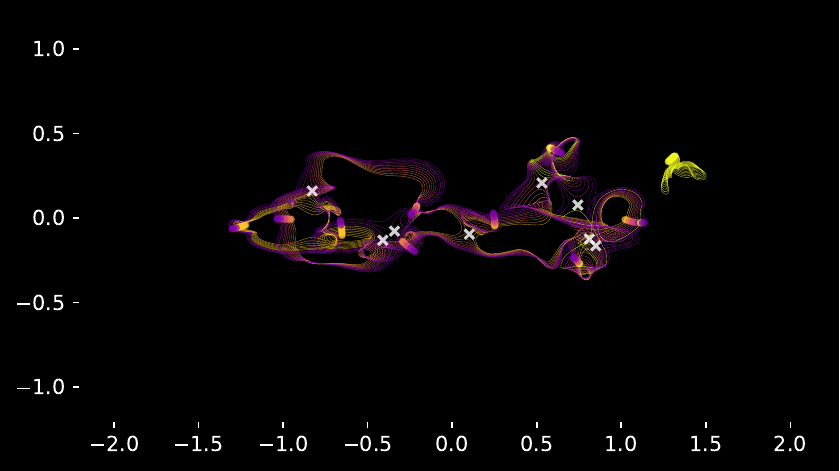}
		\includegraphics[width=0.45\textwidth]{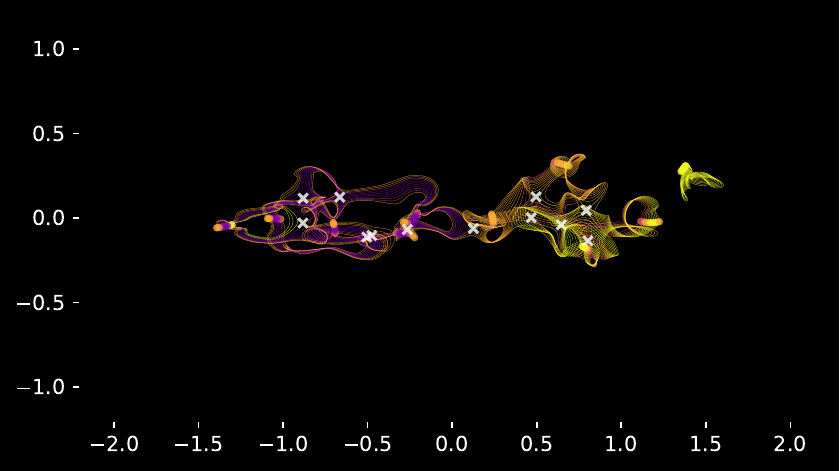}

		\includegraphics[width=0.45\textwidth]{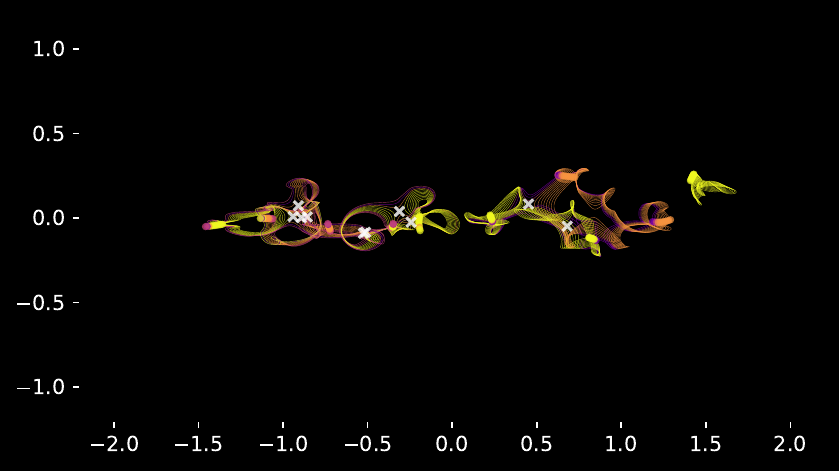}
		\includegraphics[width=0.45\textwidth]{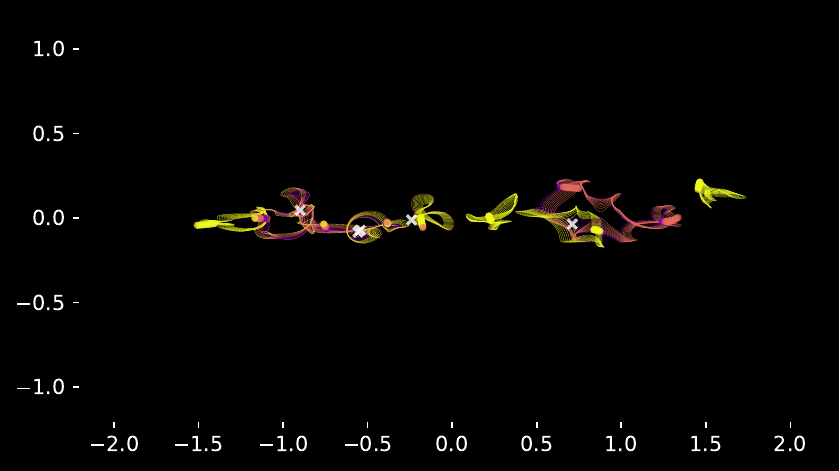}

		\caption{Eigenvalue trajectories for $s= 0.20, 0.21, \dots , 0.29, 0.30$ }
		\label{fig:pdf_image}
	\end{figure}

	\newpage

	\begin{figure}[htbp]
		\centering
		\includegraphics[width=0.45\textwidth]{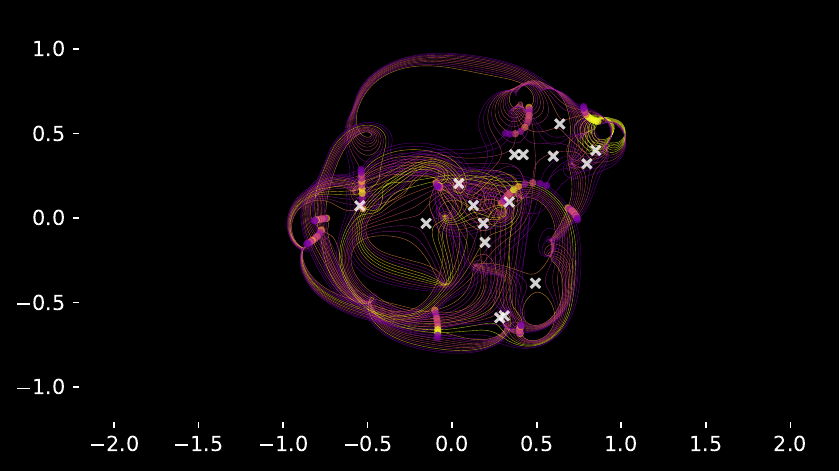}
		\includegraphics[width=0.45\textwidth]{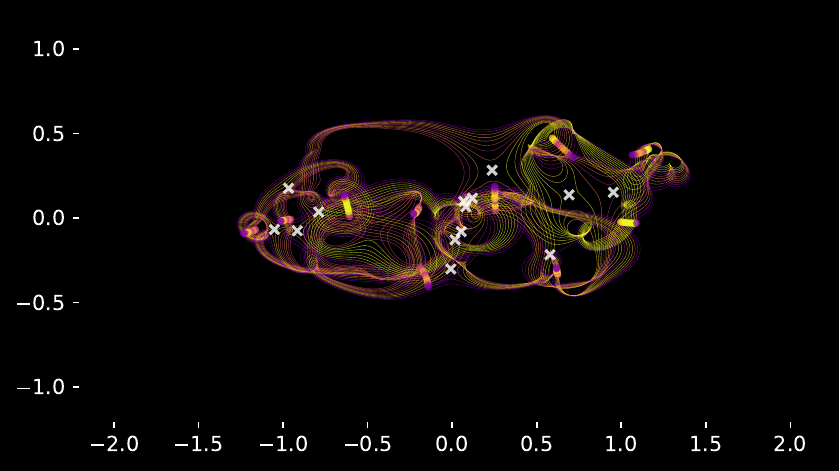}

		\includegraphics[width=0.45\textwidth]{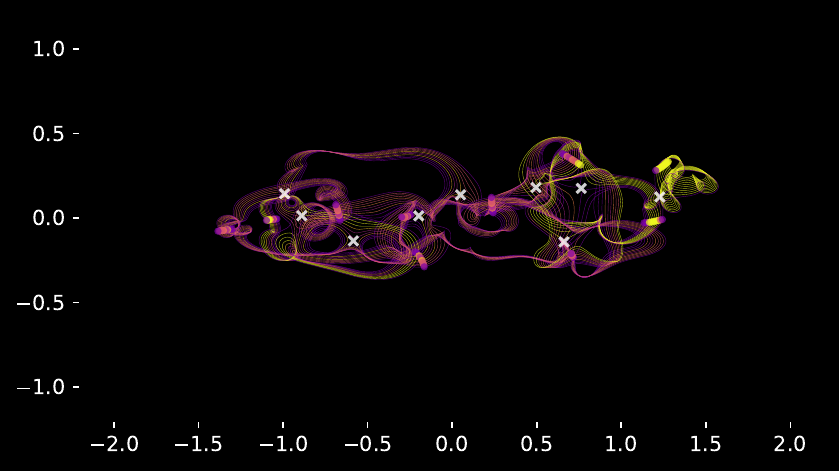}
		\includegraphics[width=0.45\textwidth]{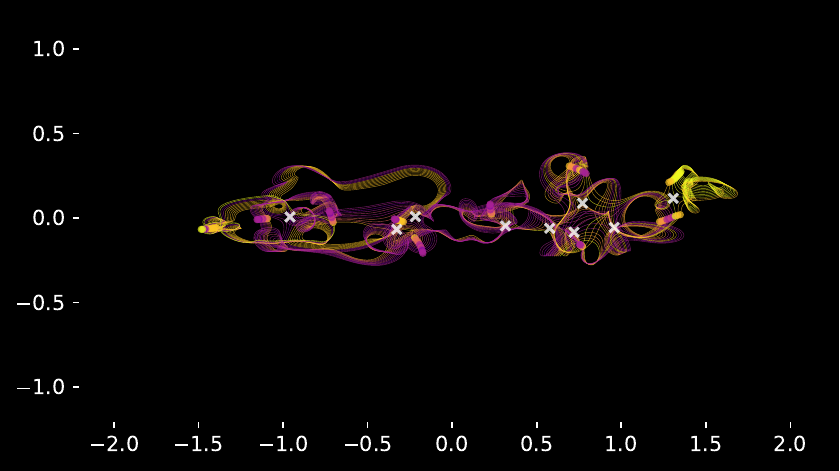}

		\includegraphics[width=0.45\textwidth]{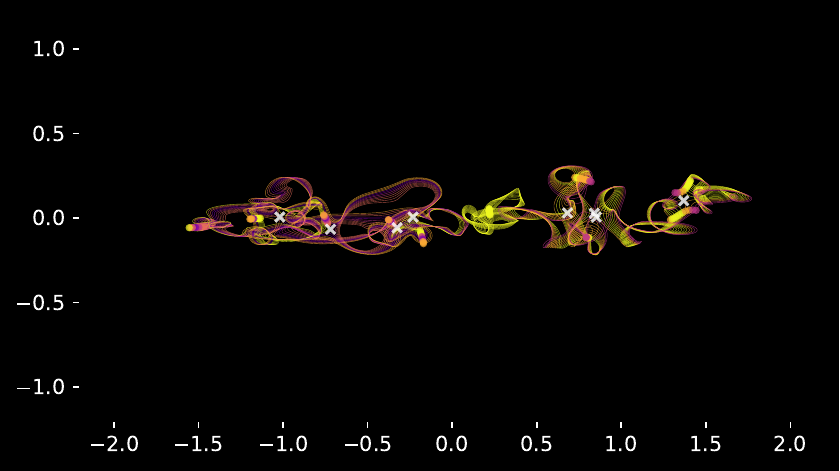}
		\includegraphics[width=0.45\textwidth]{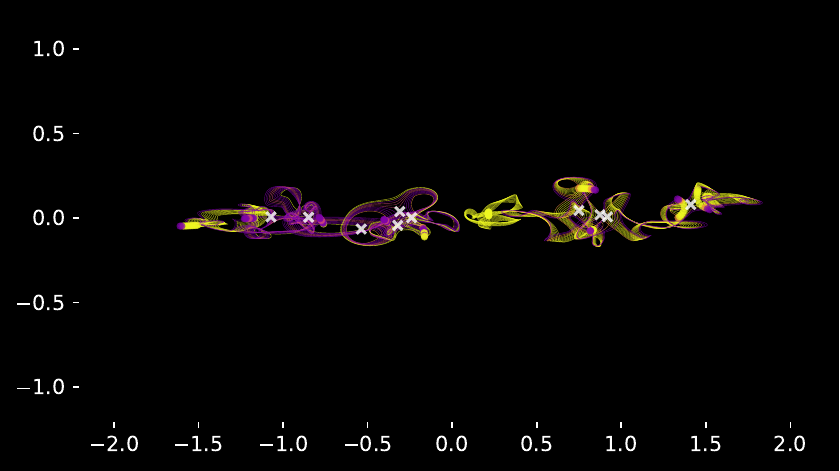}

		\caption{Eigenvalue trajectories for $s= 0.30, 0.31, \dots , 0.39, 0.40$ }
		\label{fig:pdf_image}
	\end{figure}

	\newpage

	\begin{figure}[htbp]
		\centering
		\includegraphics[width=0.45\textwidth]{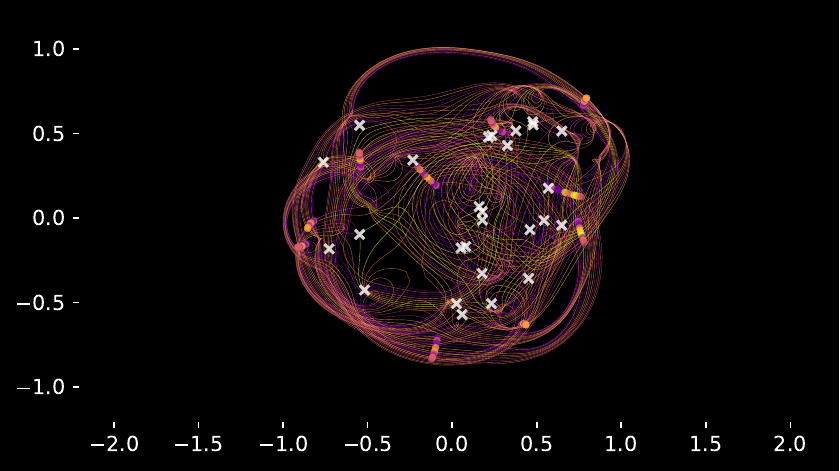}
		\includegraphics[width=0.45\textwidth]{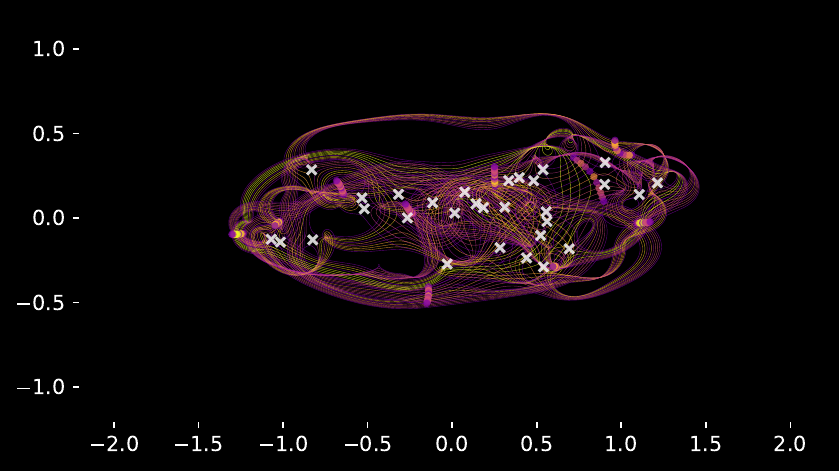}

		\includegraphics[width=0.45\textwidth]{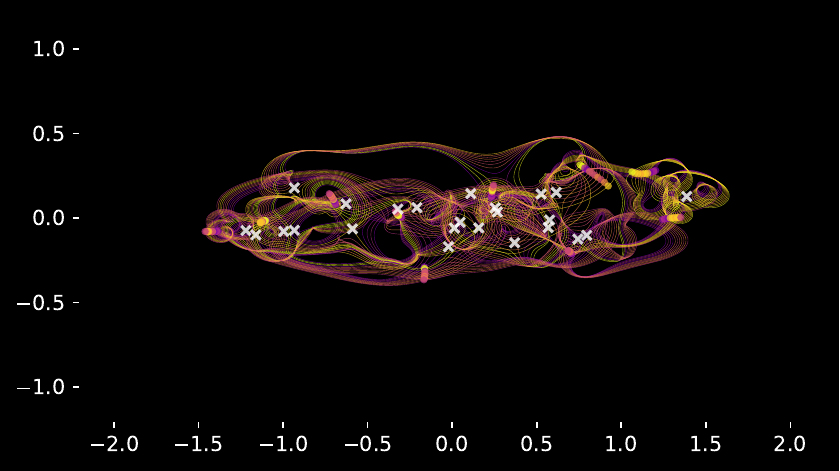}
		\includegraphics[width=0.45\textwidth]{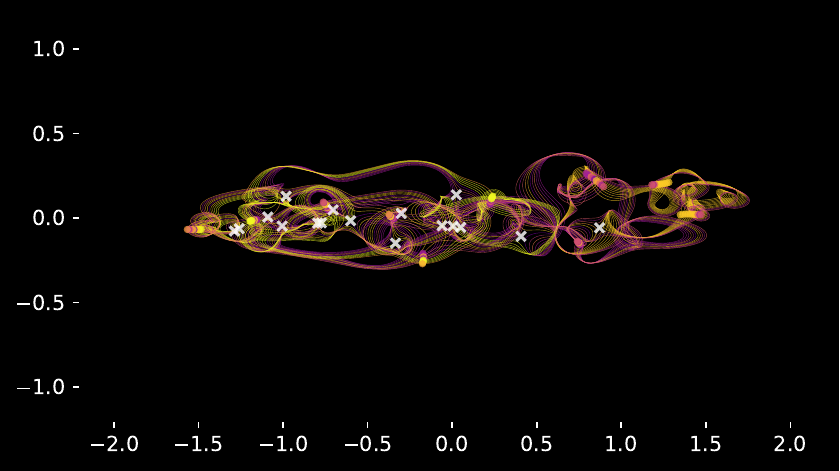}

		\includegraphics[width=0.45\textwidth]{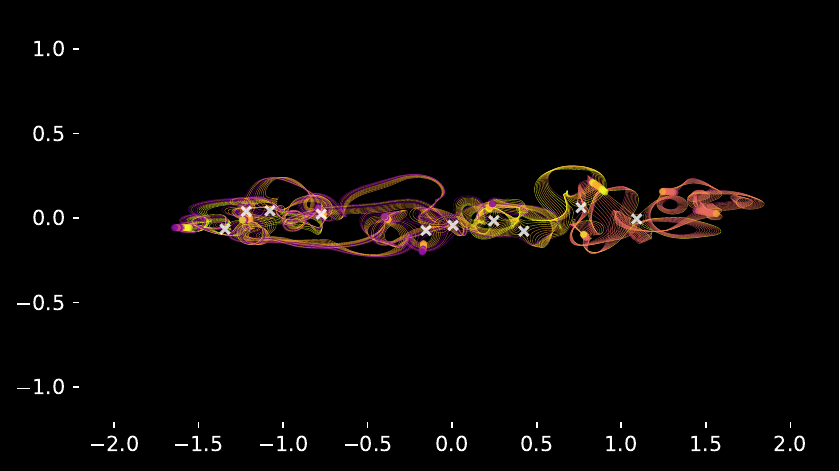}
		\includegraphics[width=0.45\textwidth]{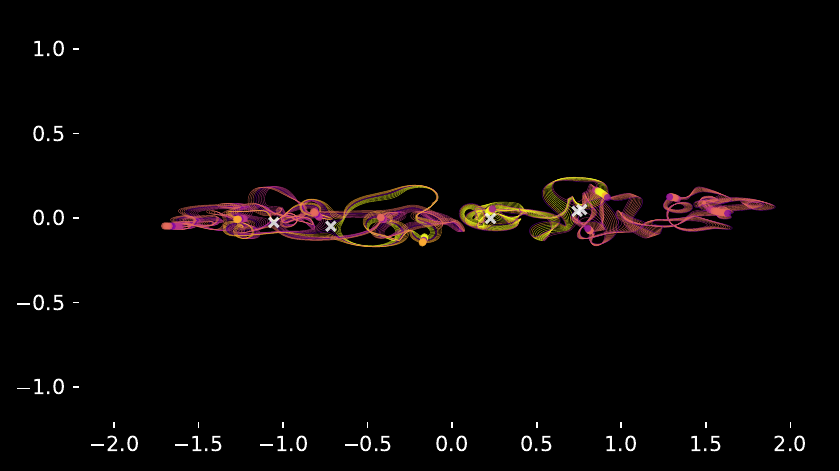}

		\caption{Eigenvalue trajectories for $s= 0.40, 0.41, \dots , 0.49, 0.50$ }
		\label{fig:pdf_image}
	\end{figure}

	\newpage

	\begin{figure}[htbp]
		\centering
		\includegraphics[width=0.45\textwidth]{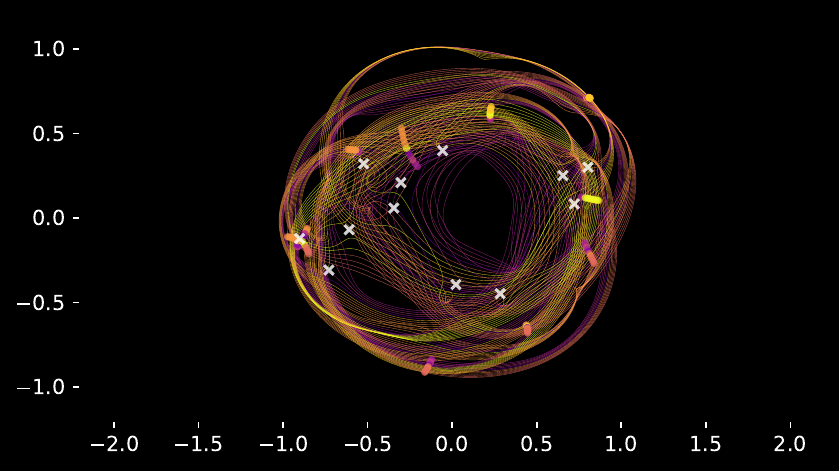}
		\includegraphics[width=0.45\textwidth]{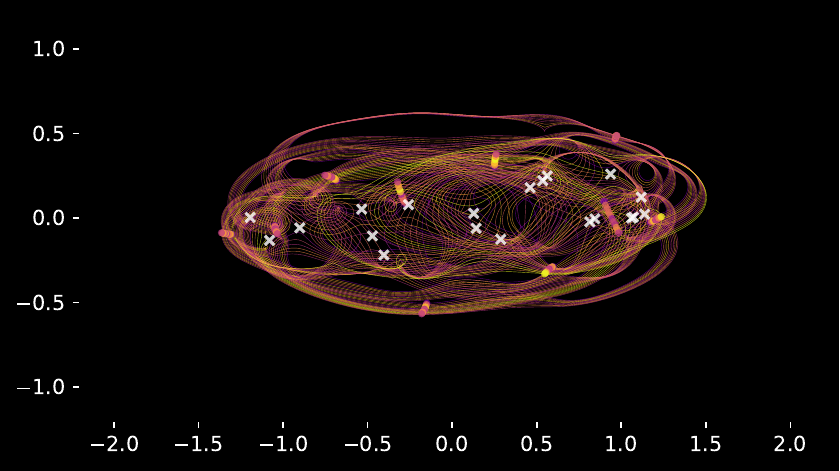}

		\includegraphics[width=0.45\textwidth]{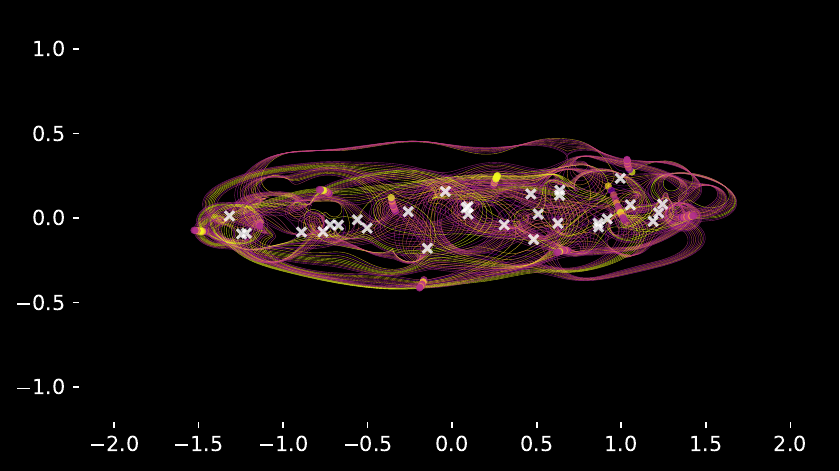}
		\includegraphics[width=0.45\textwidth]{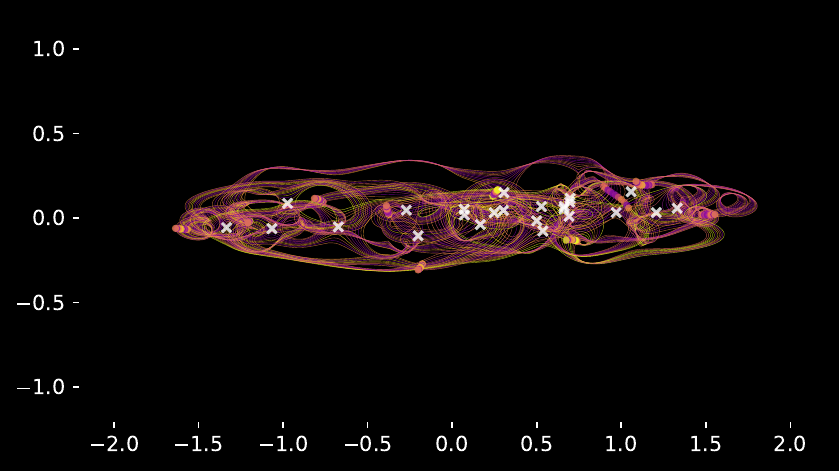}

		\includegraphics[width=0.45\textwidth]{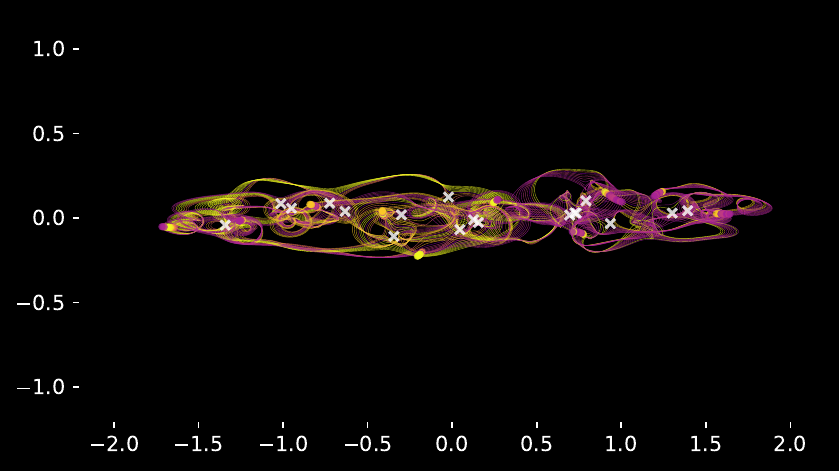}
		\includegraphics[width=0.45\textwidth]{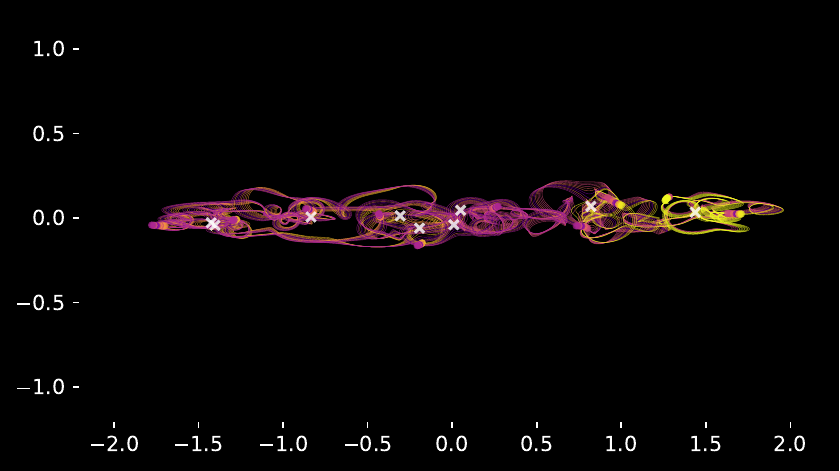}

		\caption{Eigenvalue trajectories for $s= 0.50, 0.51, \dots , 0.59, 0.60$ }
		\label{fig:pdf_image}
	\end{figure}

	\newpage

	\begin{figure}[htbp]
		\centering
		\includegraphics[width=0.45\textwidth]{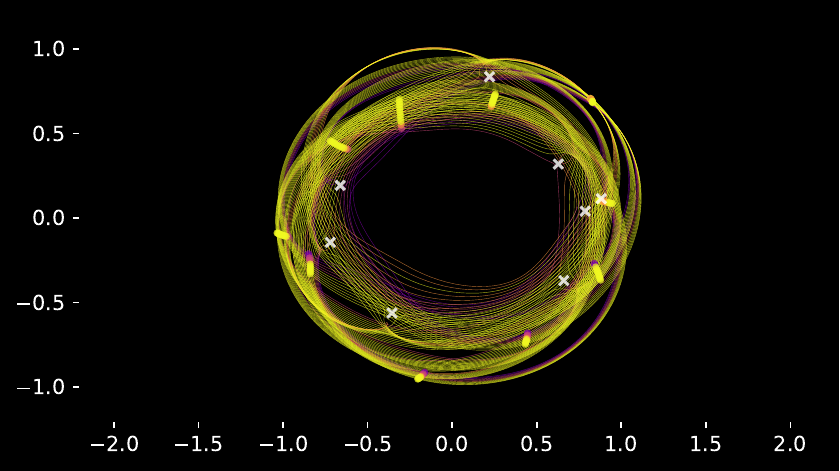}
		\includegraphics[width=0.45\textwidth]{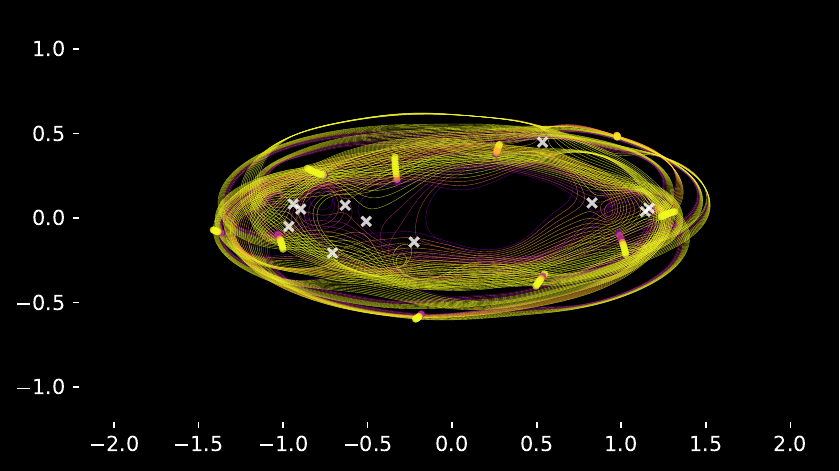}

		\includegraphics[width=0.45\textwidth]{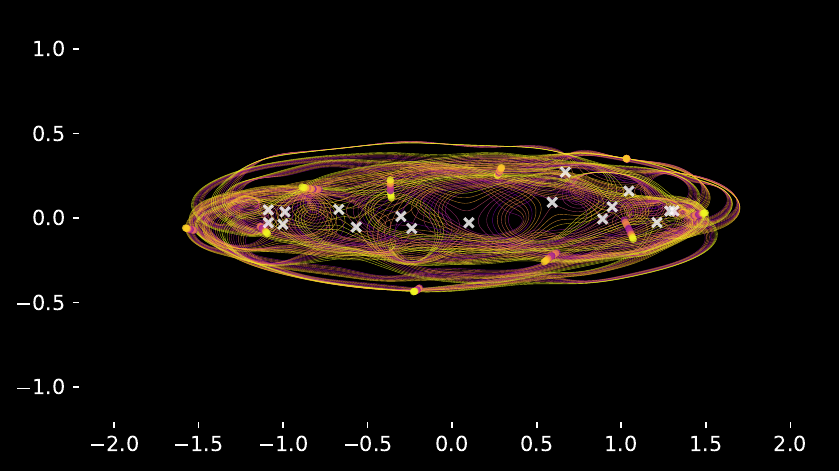}
		\includegraphics[width=0.45\textwidth]{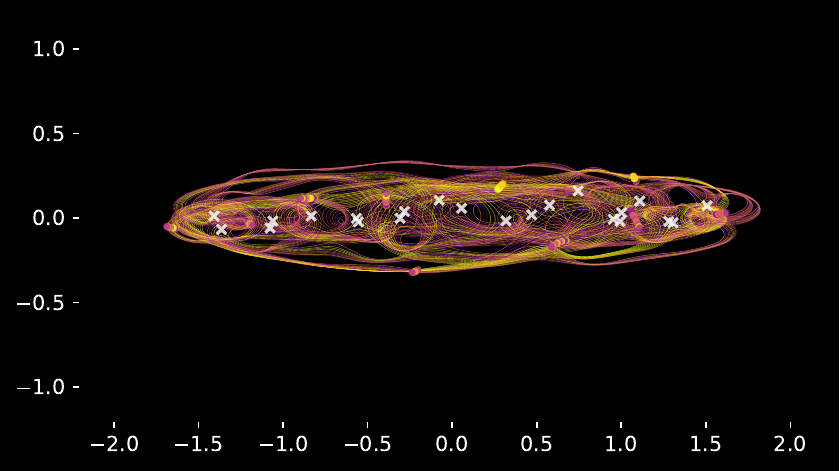}

		\includegraphics[width=0.45\textwidth]{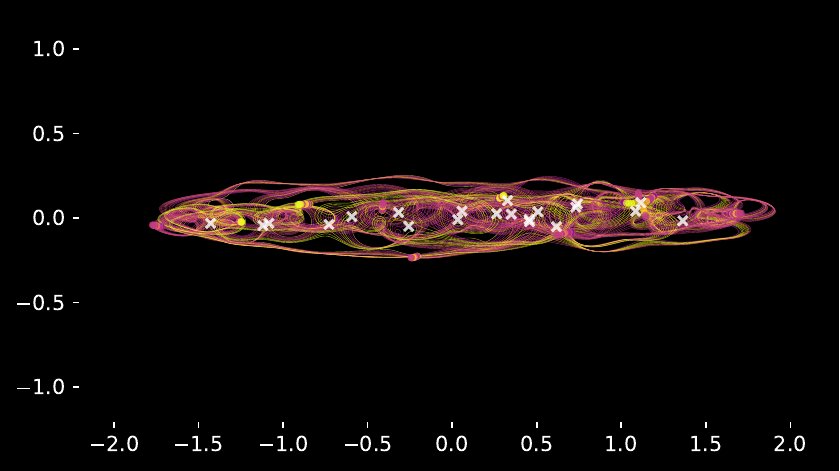}
		\includegraphics[width=0.45\textwidth]{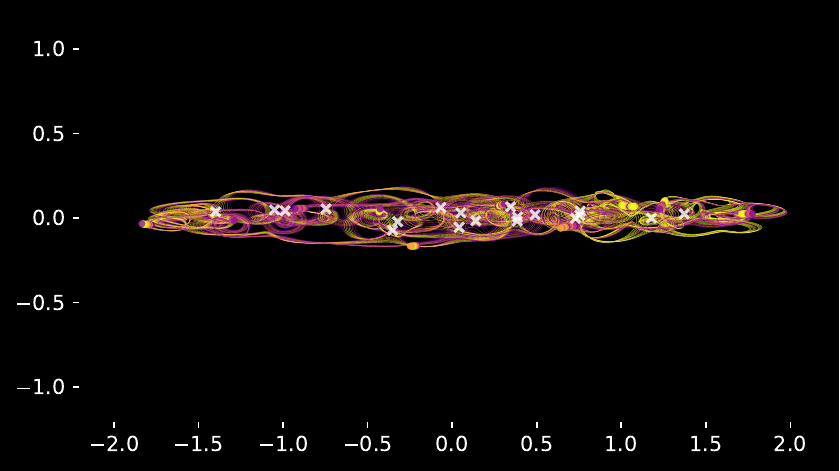}

		\caption{Eigenvalue trajectories for $s= 0.60, 0.61, \dots , 0.69, 0.70$ }
		\label{fig:pdf_image}
	\end{figure}

	\newpage

	\begin{figure}[htbp]
		\centering
		\includegraphics[width=0.45\textwidth]{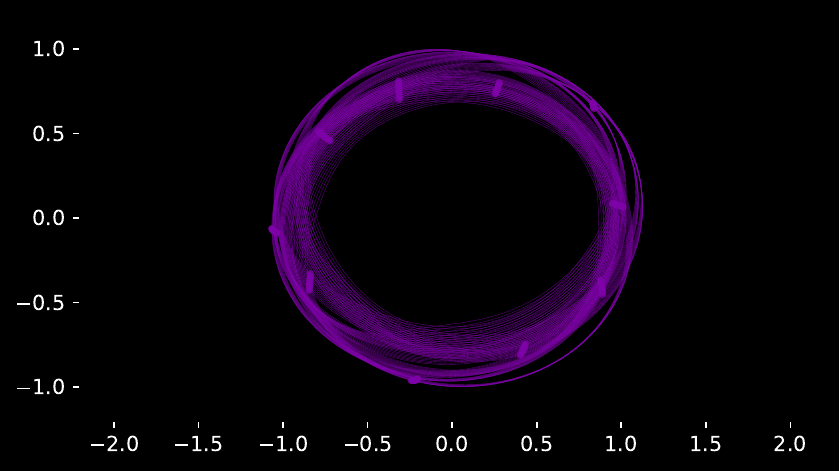}
		\includegraphics[width=0.45\textwidth]{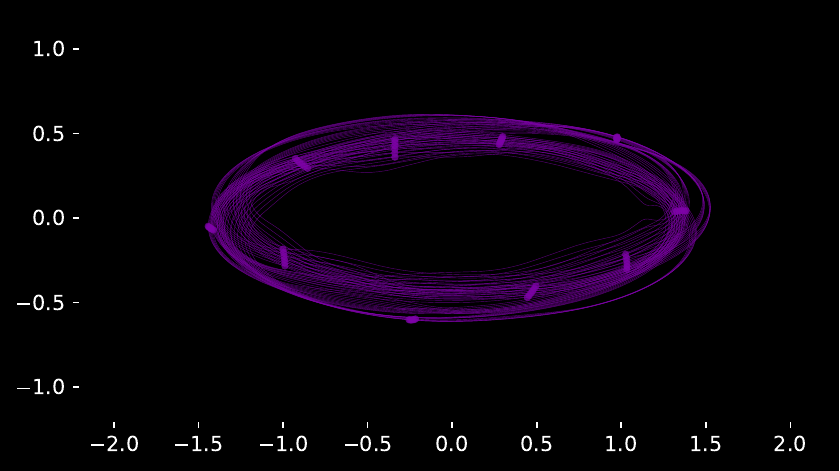}

		\includegraphics[width=0.45\textwidth]{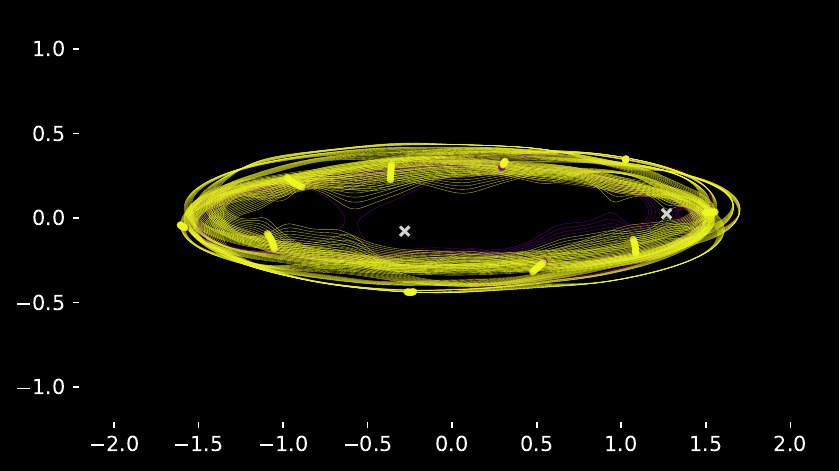}
		\includegraphics[width=0.45\textwidth]{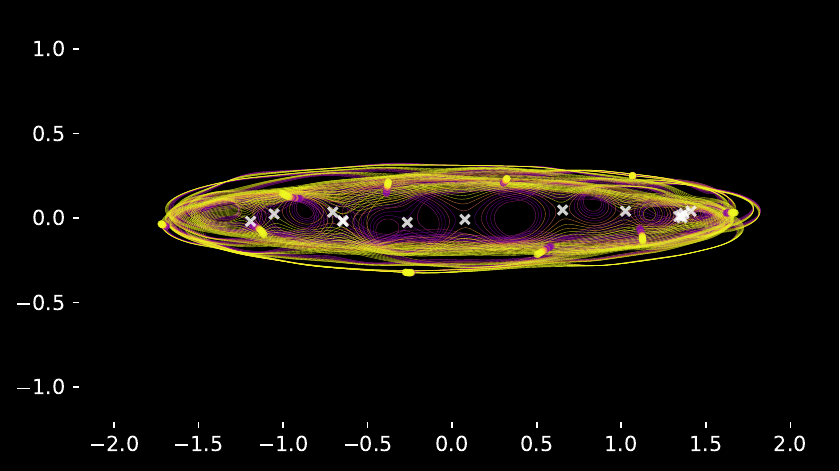}

		\includegraphics[width=0.45\textwidth]{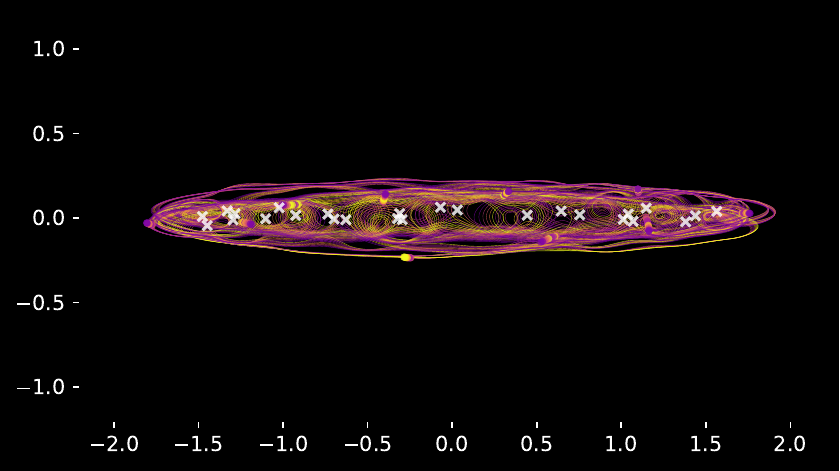}
		\includegraphics[width=0.45\textwidth]{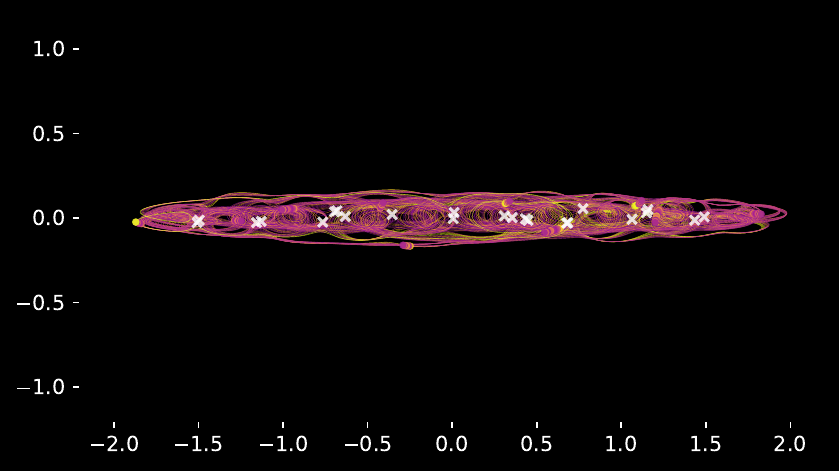}

		\caption{Eigenvalue trajectories for $s= 0.70, 0.71, \dots , 0.78, 0.80$ }
		\label{fig:pdf_image}
	\end{figure}

	\newpage

	\begin{figure}[htbp]
		\centering
		\includegraphics[width=0.45\textwidth]{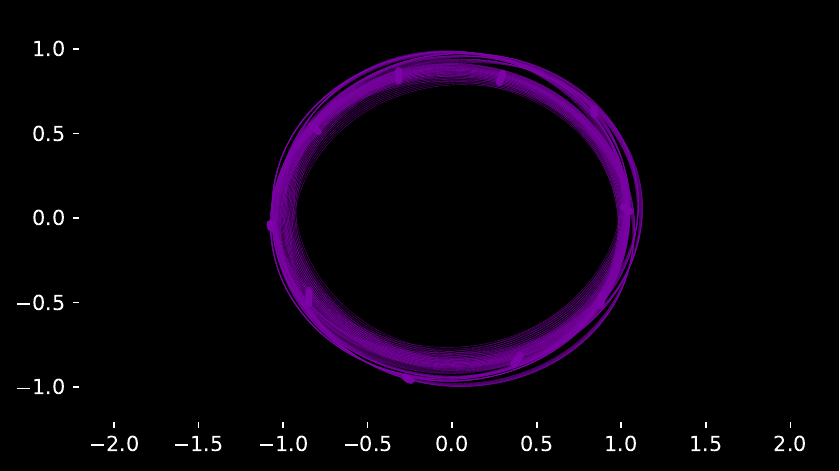}
		\includegraphics[width=0.45\textwidth]{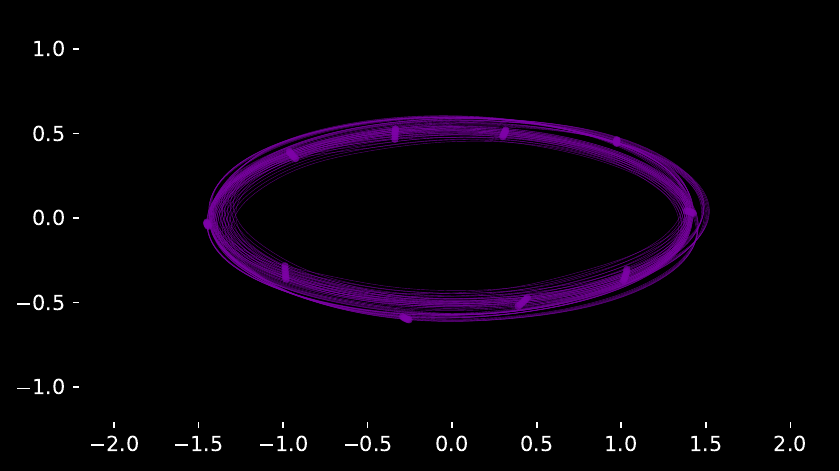}

		\includegraphics[width=0.45\textwidth]{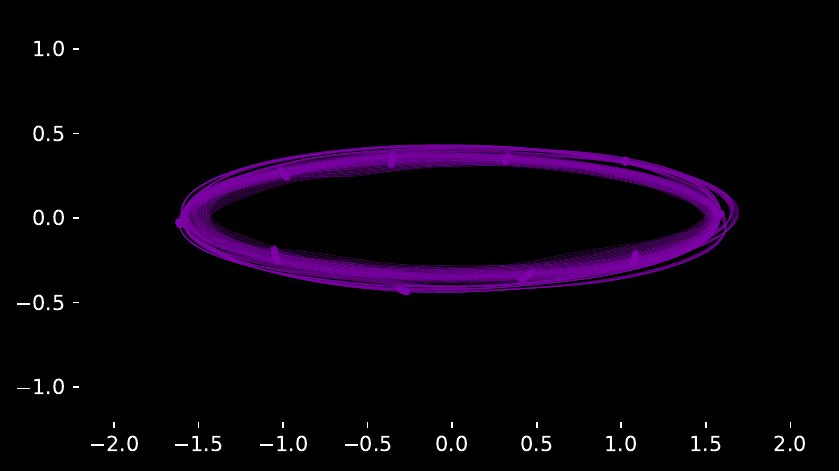}
		\includegraphics[width=0.45\textwidth]{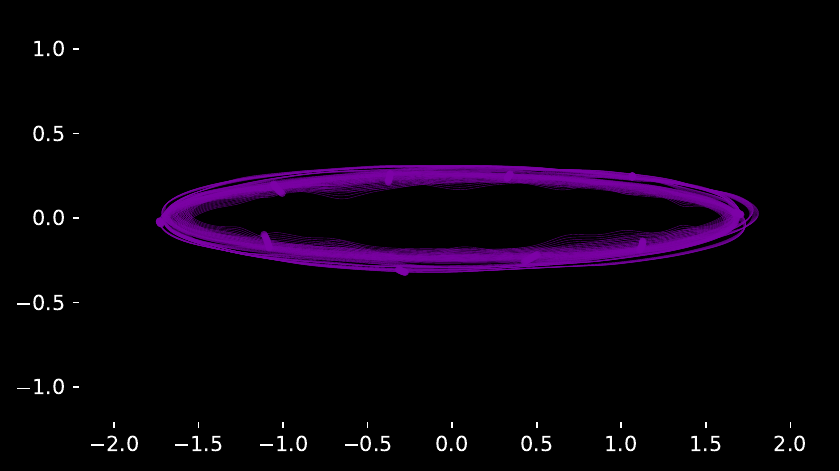}

		\includegraphics[width=0.45\textwidth]{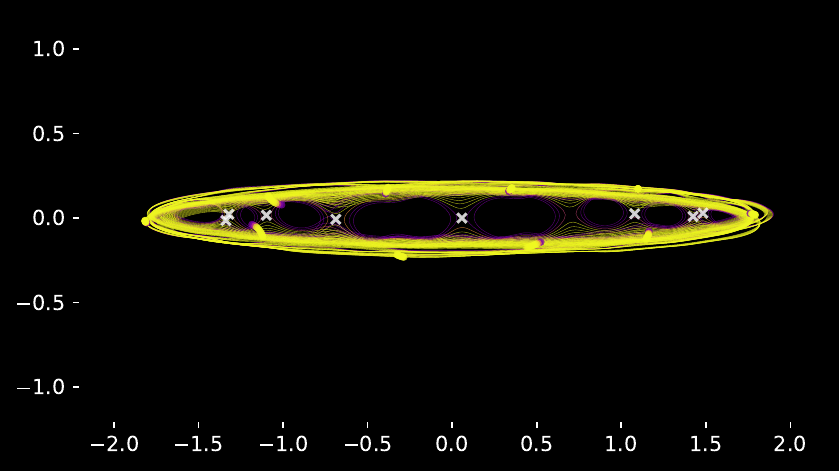}
		\includegraphics[width=0.45\textwidth]{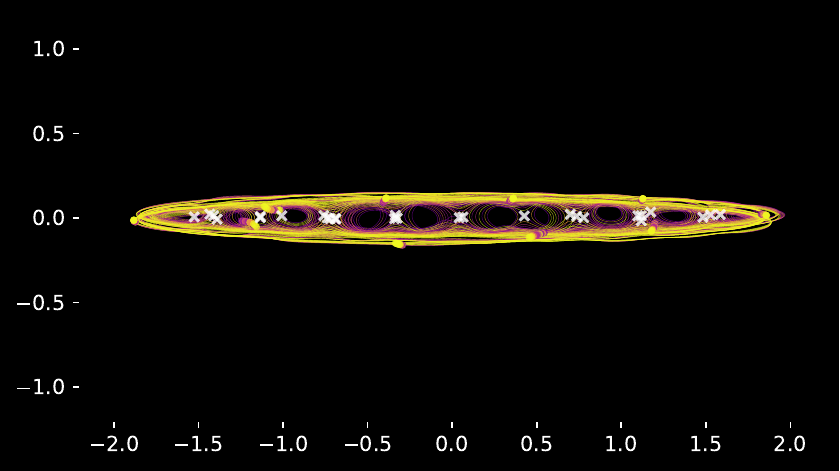}

		\caption{Eigenvalue trajectories for $s= 0.80, 0.81, \dots , 0.89, 0.90$ }
		\label{fig:pdf_image}
	\end{figure}

	\newpage

	\begin{figure}[htbp]
		\centering
		\includegraphics[width=0.45\textwidth]{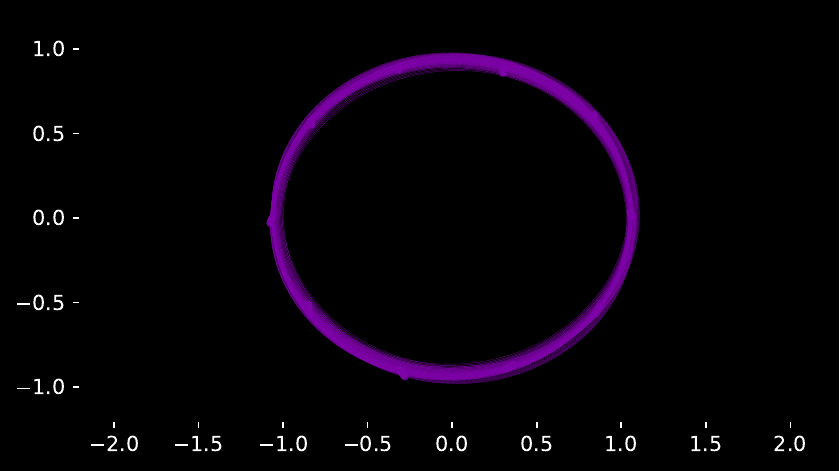}
		\includegraphics[width=0.45\textwidth]{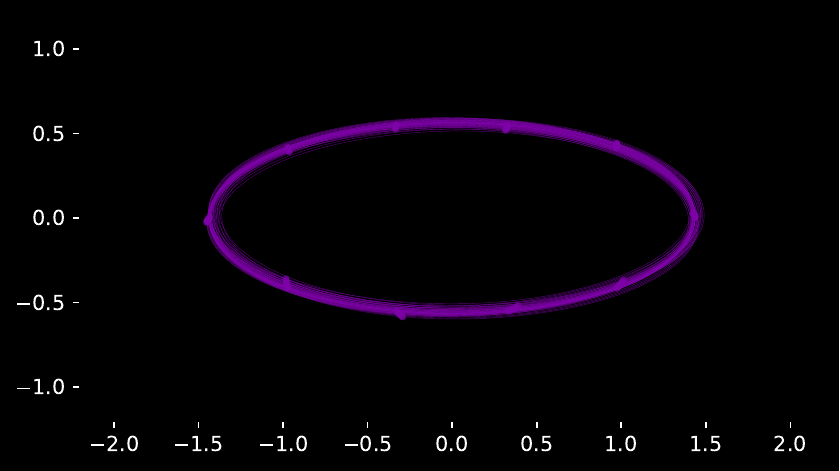}

		\includegraphics[width=0.45\textwidth]{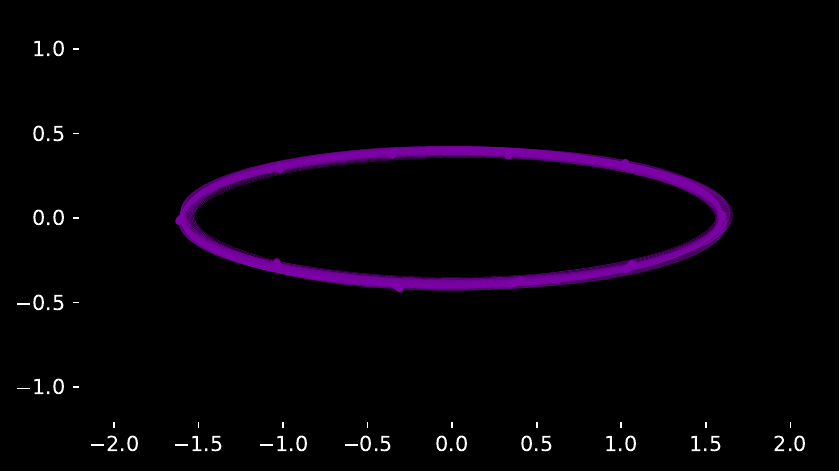}
		\includegraphics[width=0.45\textwidth]{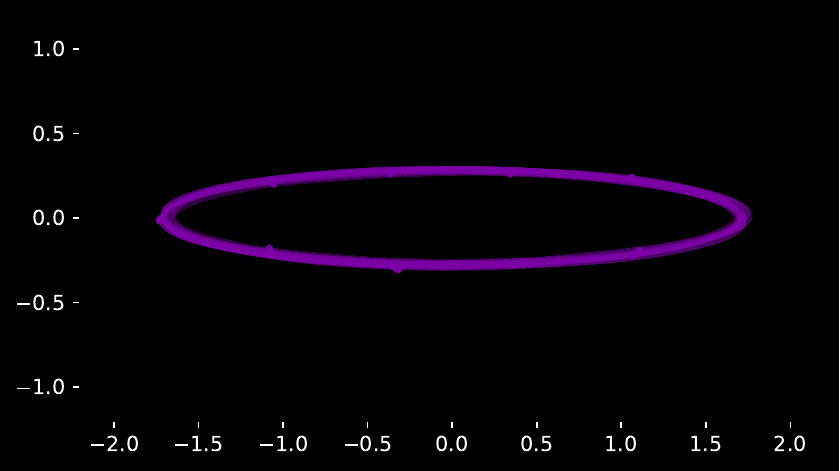}

		\includegraphics[width=0.45\textwidth]{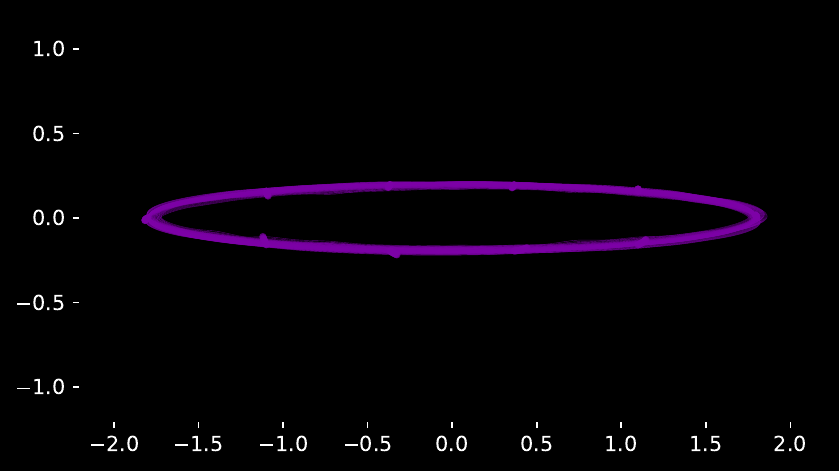}
		\includegraphics[width=0.45\textwidth]{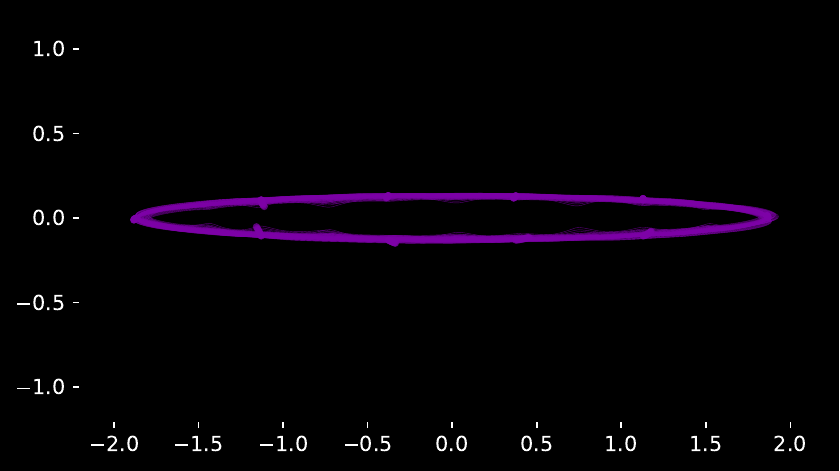}

		\caption{Eigenvalue trajectories for $s= 0.91, 0.92, \dots , 0.99, 1.0$ }
		\label{fig:pdf_image}
	\end{figure}

	\newpage

	\section{Initial tracks for cases with multiple initial eigenvalues} \label{appendix:repeated-eigenvalues}

	We include the eigenvalue tracks for the first 
	$s$-window ($s= 0.001, 0.01, 0.02, 0.03, \dots, 0.09. 0.10$) 
	of all fifteen traceless Bernoulli trials (out of $100$, from seeds $1000$ to $1099$) 
	with repeated initial eigenvalues. The seeds can be used to reproduce all results in this paper, 
	using the package in the repository \url{https://github.com/carlos-vargas-math/eigenvalue-collisions}.

	The cases with fewer distinct eigenvalues are seed $1001$, with a triple eigenvalue 
	and seeds $1020$, $1070$ 
	with two pairs of repeated eigenvalues.
	The triple eigenvalue from seed $1001$ produces a strong shockwave with three helices.
	Notice that sometimes the repeated eigenvalues clash immediately 
	and in some other cases they synchronize harmoniously without colliding. 
	
	A notable case of this is seed $1020$
	with two repeated eigenvalues, 
	one pair exploding immediately and the other synchronizing without colliding. 
	On the other hand, seed $1070$ has two exploding eigenvalues.
	From the rest of the seeds, with a single repeated eigenvalue,
	$1005$, $1009$, $1044$, $1047$, have exploding repeated eigenvalues, 
	while and $1017$, $1038$, $1043$, $1048$, $1072$, $1096$ have synchronized repeated eigenvalues.

	Our method had some trouble with the initial tracks for seeds $1001$, $1035$ and $1094$.
	The first track displayed in the figures for seeds $1001$, $1035$ and $1094$ 
	are for the values $s= 0.004$, $s= 0.017$, and $s=0.08$ respectively. 
	All other seeds tracks are indeed for $s=0.001, 0.01, \dots, 0.10$.
	 
	\newpage

	\begin{figure}[htbp]
		\centering
		\includegraphics[width=0.45\textwidth]{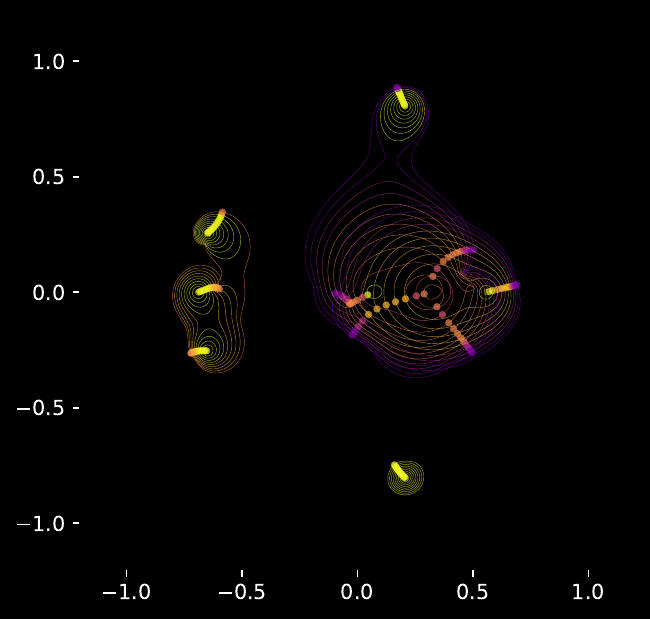}
		\includegraphics[width=0.45\textwidth]{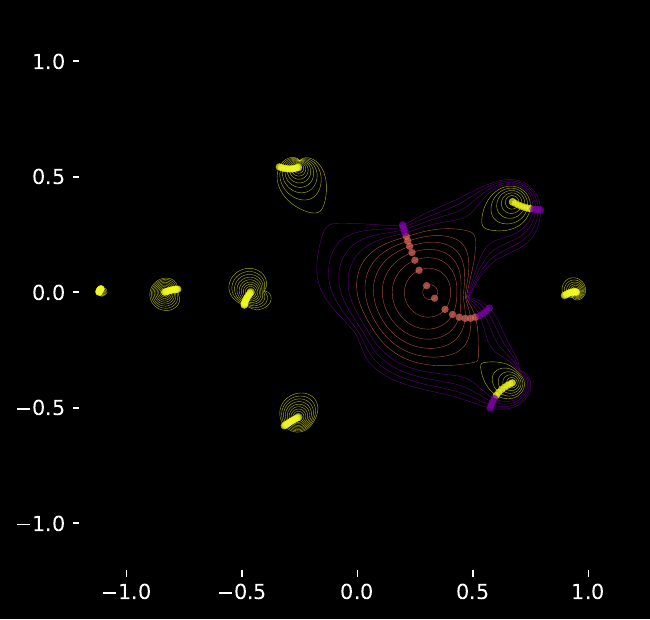}

		\includegraphics[width=0.45\textwidth]{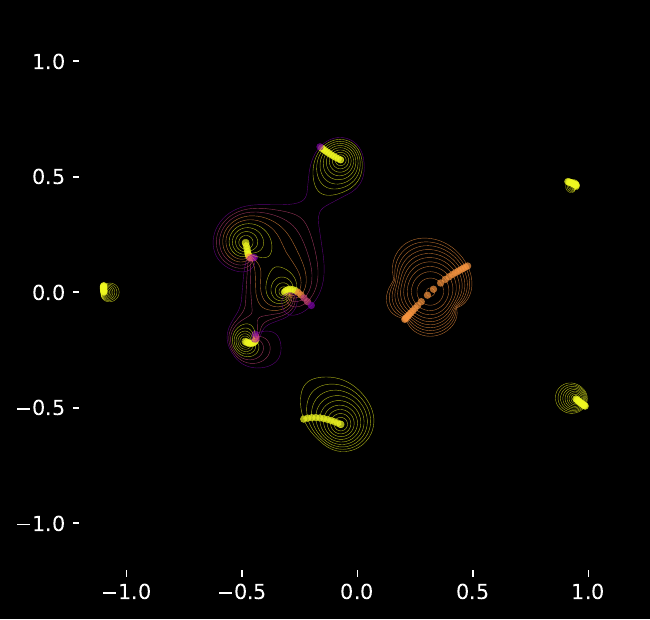}
		\includegraphics[width=0.45\textwidth]{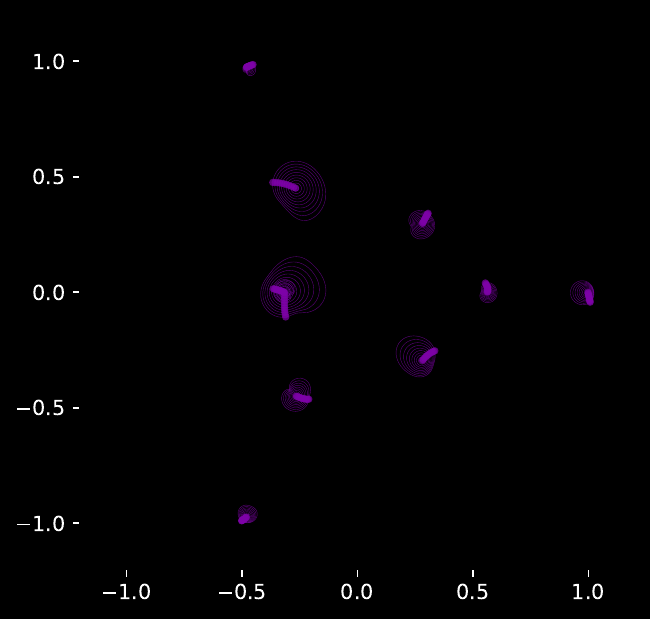}
		\caption{Eigenvalue trajectories for $s= 0.01, 0.02, \dots , 0.09, 0.10$. Seeds 1001, 1005, 1009, 1017}
		\label{fig:pdf_image}
	\end{figure}

	\newpage

	\begin{figure}[htbp]
		\centering
		\includegraphics[width=0.45\textwidth]{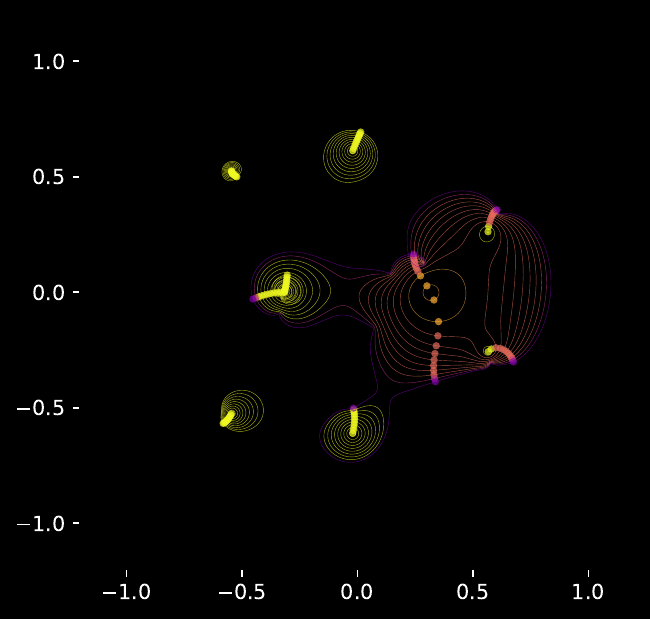}
		\includegraphics[width=0.45\textwidth]{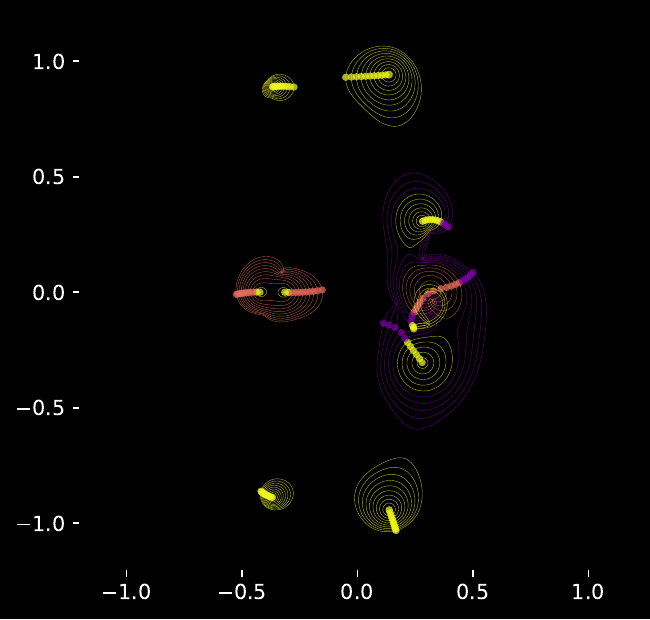}

		\includegraphics[width=0.45\textwidth]{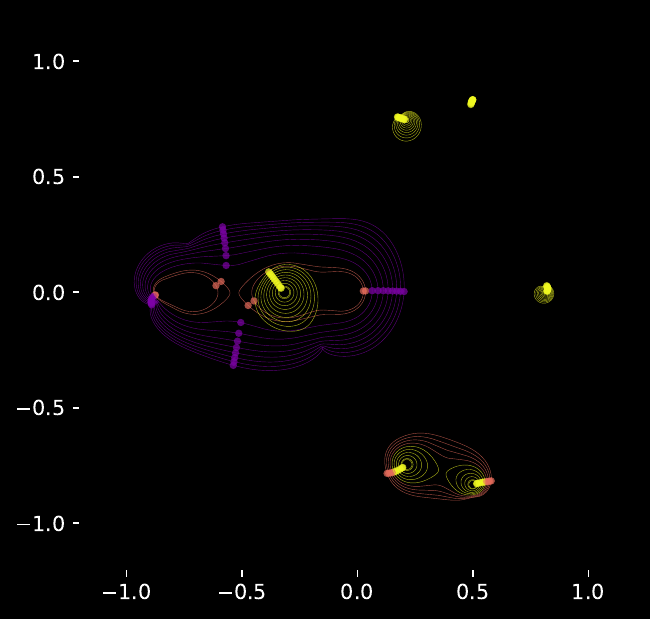}
		\includegraphics[width=0.45\textwidth]{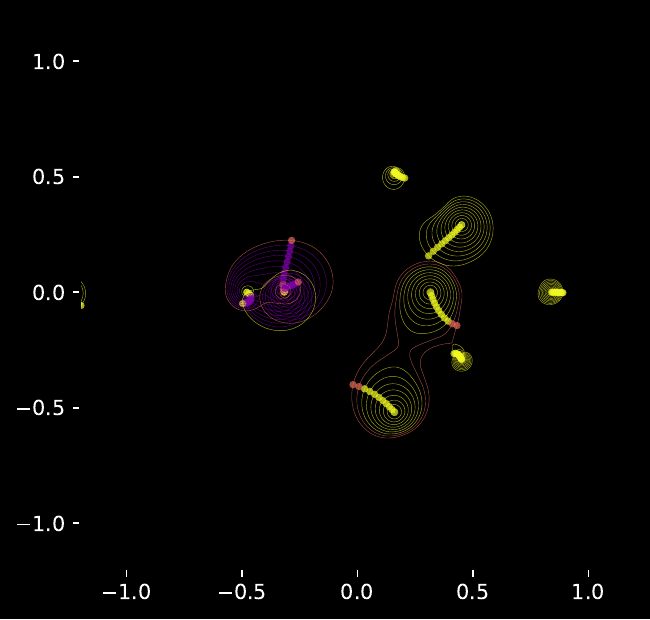}
		\caption{Eigenvalue trajectories for $s= 0.01, 0.02, \dots , 0.09, 0.10$. Seeds 1020, 1022, 1035, 1038}

		\label{fig:pdf_image}
 	\end{figure}

	 \newpage

	 \begin{figure}[htbp]
		 \centering
		 \includegraphics[width=0.45\textwidth]{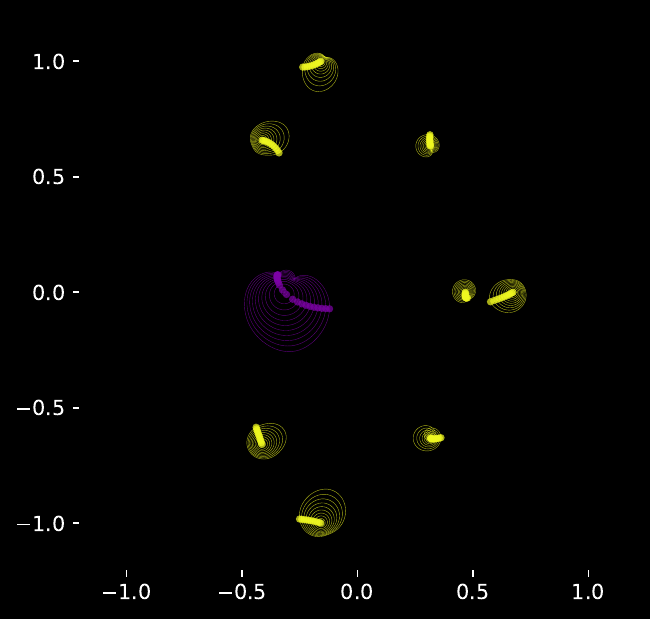}
		 \includegraphics[width=0.45\textwidth]{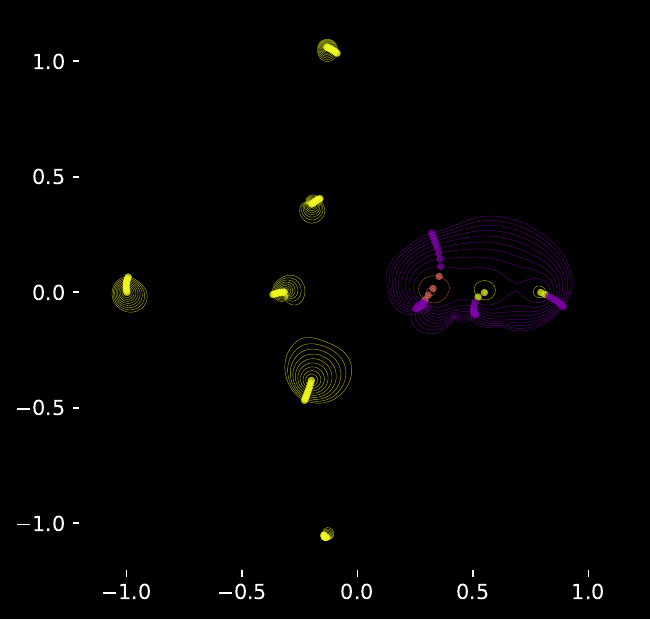}
 
		 \includegraphics[width=0.45\textwidth]{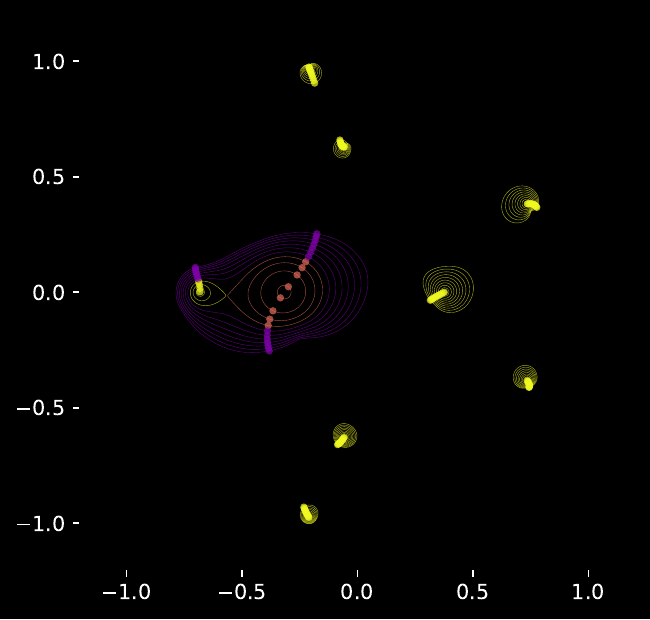}
		 \includegraphics[width=0.45\textwidth]{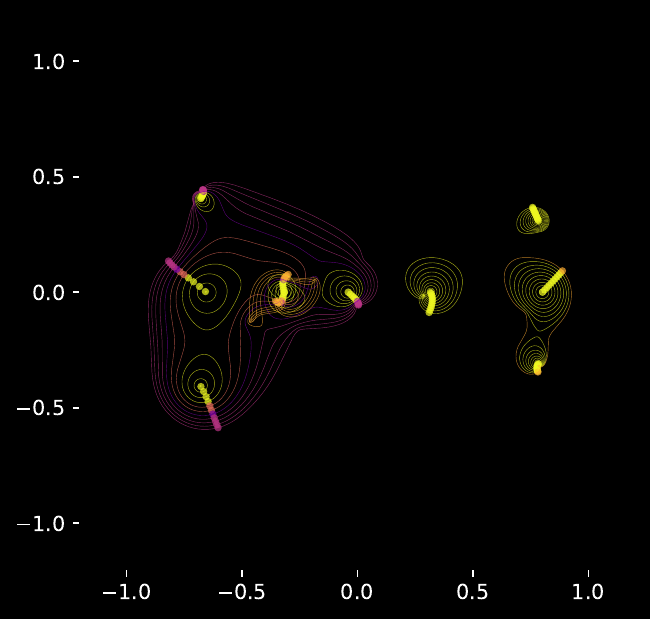}
		 \caption{Eigenvalue trajectories for $s= 0.01, 0.02, \dots , 0.09, 0.10$. Seeds 1043, 1044, 1047, 1048}
		 \label{fig:pdf_image}
	  \end{figure}

	  \newpage

	  \begin{figure}[htbp]
		  \centering
		  \includegraphics[width=0.45\textwidth]{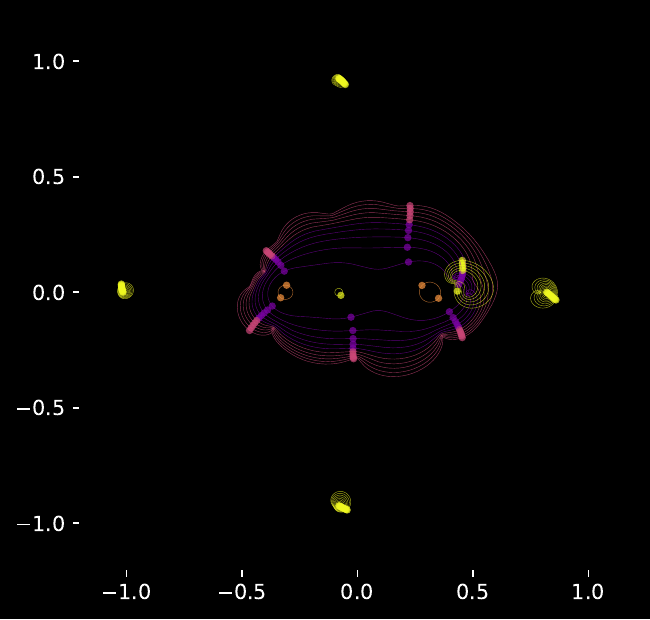}
		  \includegraphics[width=0.45\textwidth]{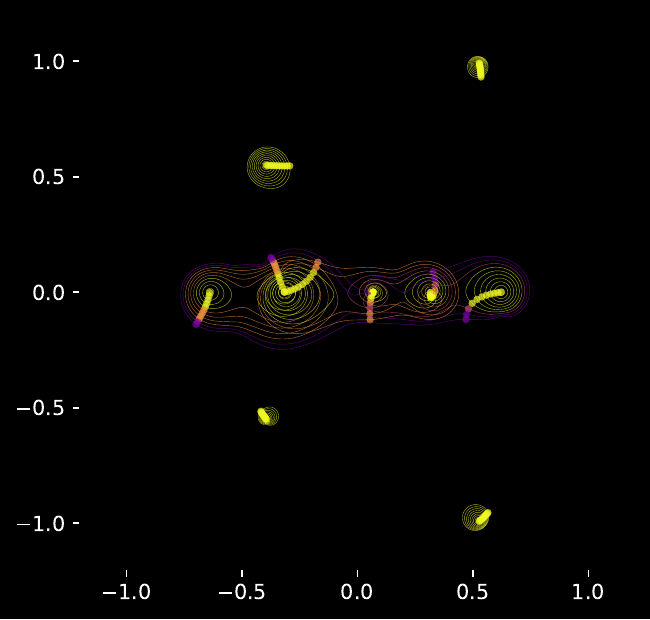}
  
		  \includegraphics[width=0.45\textwidth]{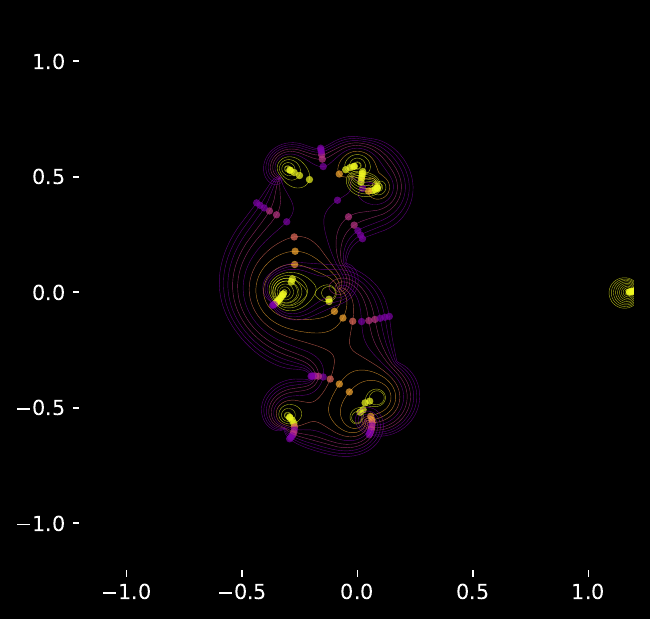}
		  \caption{Eigenvalue trajectories for $s= 0.01, 0.02, \dots , 0.09, 0.10$. Seeds 1070, 1072, 1094}
		  \label{fig:pdf_image}
	   \end{figure}

	   \newpage

	   \begin{figure}[htbp]
		   \centering
		   \includegraphics[width=0.8\textwidth]{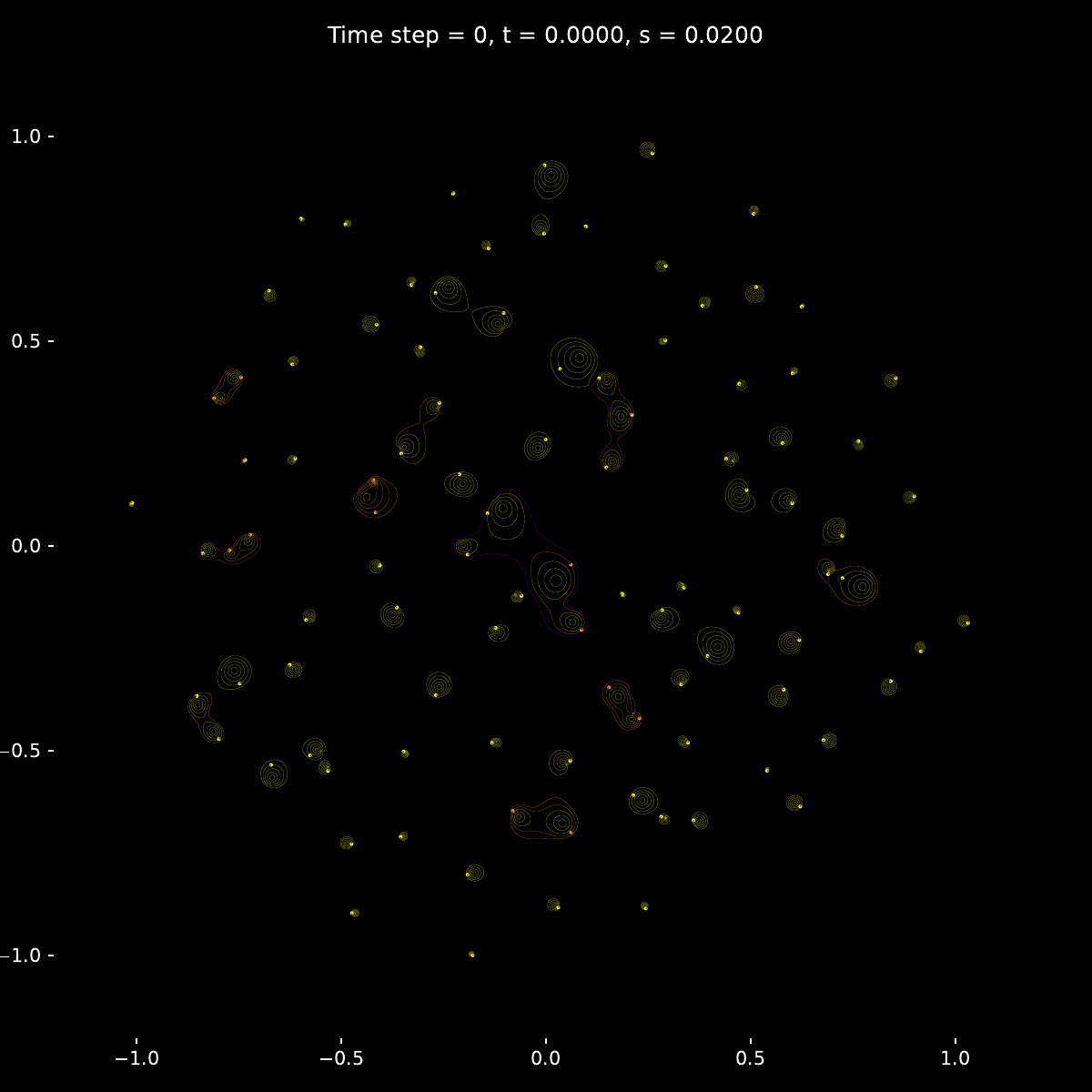}
		   \caption{Eigenvalue trajectories for $s= 0.025, 0.030, 0.035, 0.040, 0.045, 0.050$ }
		   \label{fig:pdf_image}
	   \end{figure}	 

	\section{Eigenvalue tracks for N=100, curve = circle} \label{appendix:N=100}

	We include some high resolution pictures for N=100 of the eigenvalue tracks.

	\newpage

	\begin{figure}[htbp]
		\centering
		\includegraphics[width=0.8\textwidth]{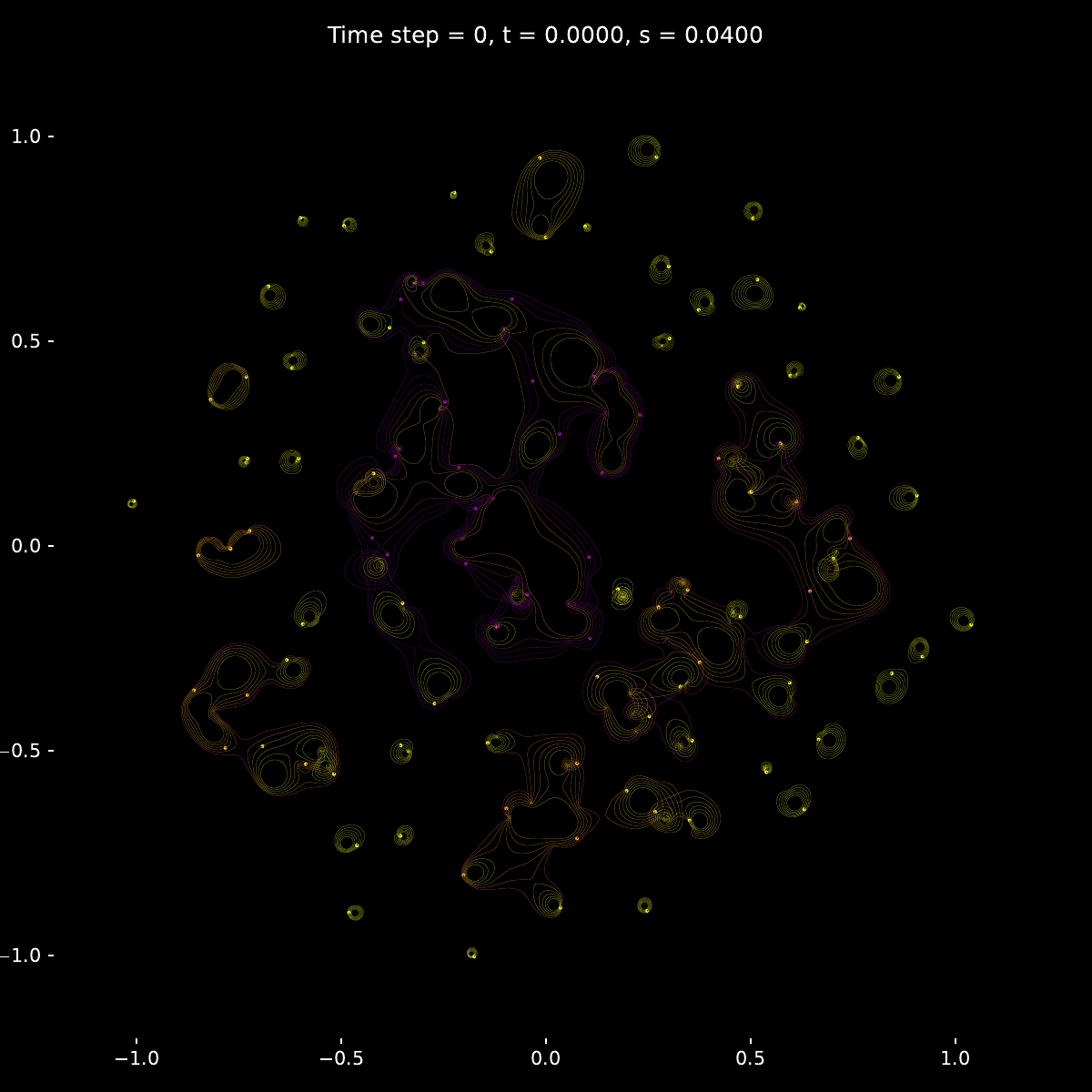}
		\caption{Eigenvalue trajectories for $s= 0.025, 0.030, 0.035, 0.040, 0.045, 0.050$ }
		\label{fig:pdf_image}
	\end{figure}
	
	\newpage

	\begin{figure}[htbp]
		\centering
		\includegraphics[width=0.8\textwidth]{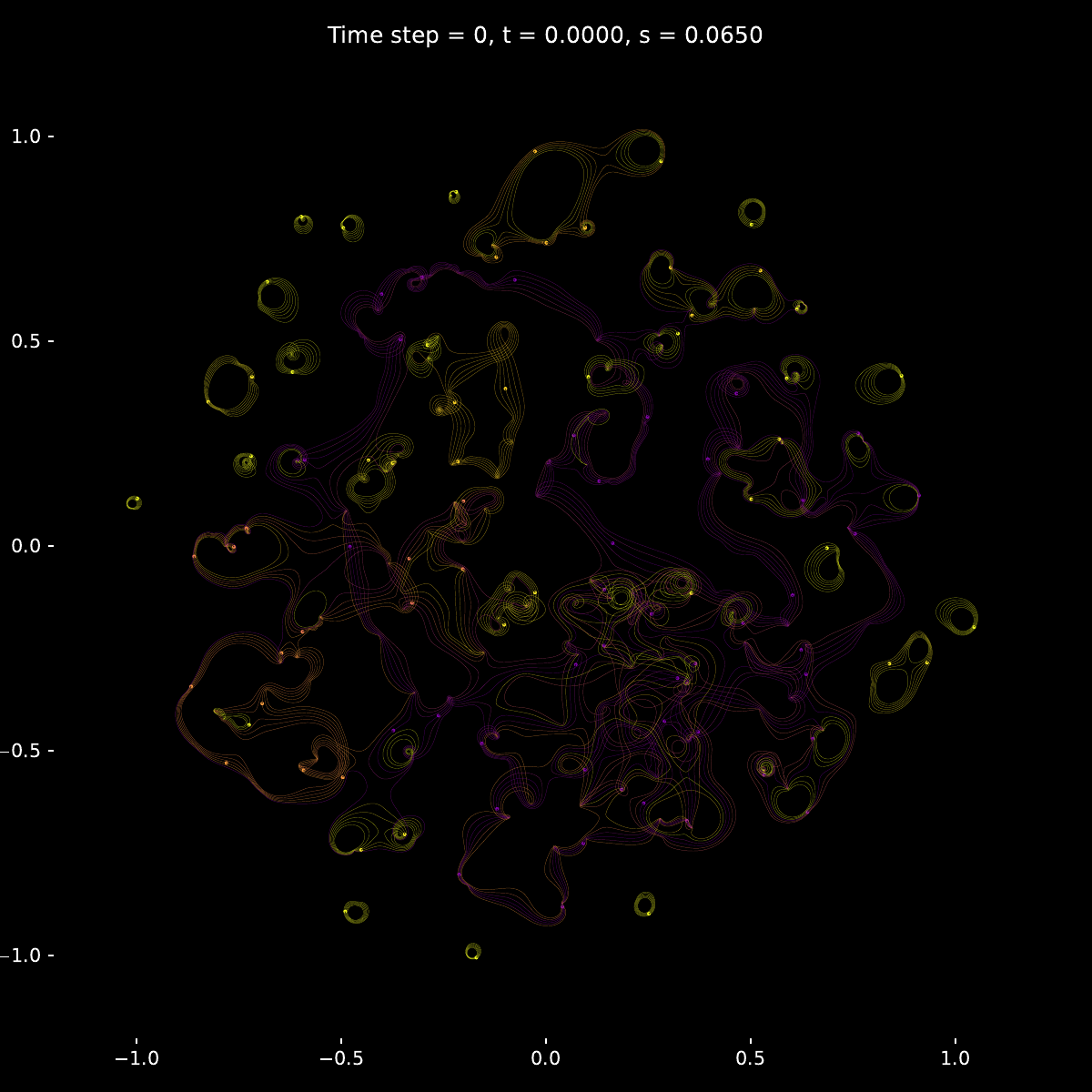}
		\caption{Eigenvalue trajectories for $s= 0.050, 0.055, 0.060, 0.065, 0.070, 0.075$ }
		\label{fig:pdf_image}
	\end{figure}

	\newpage

	\begin{figure}[htbp]
		\centering
		\includegraphics[width=0.8\textwidth]{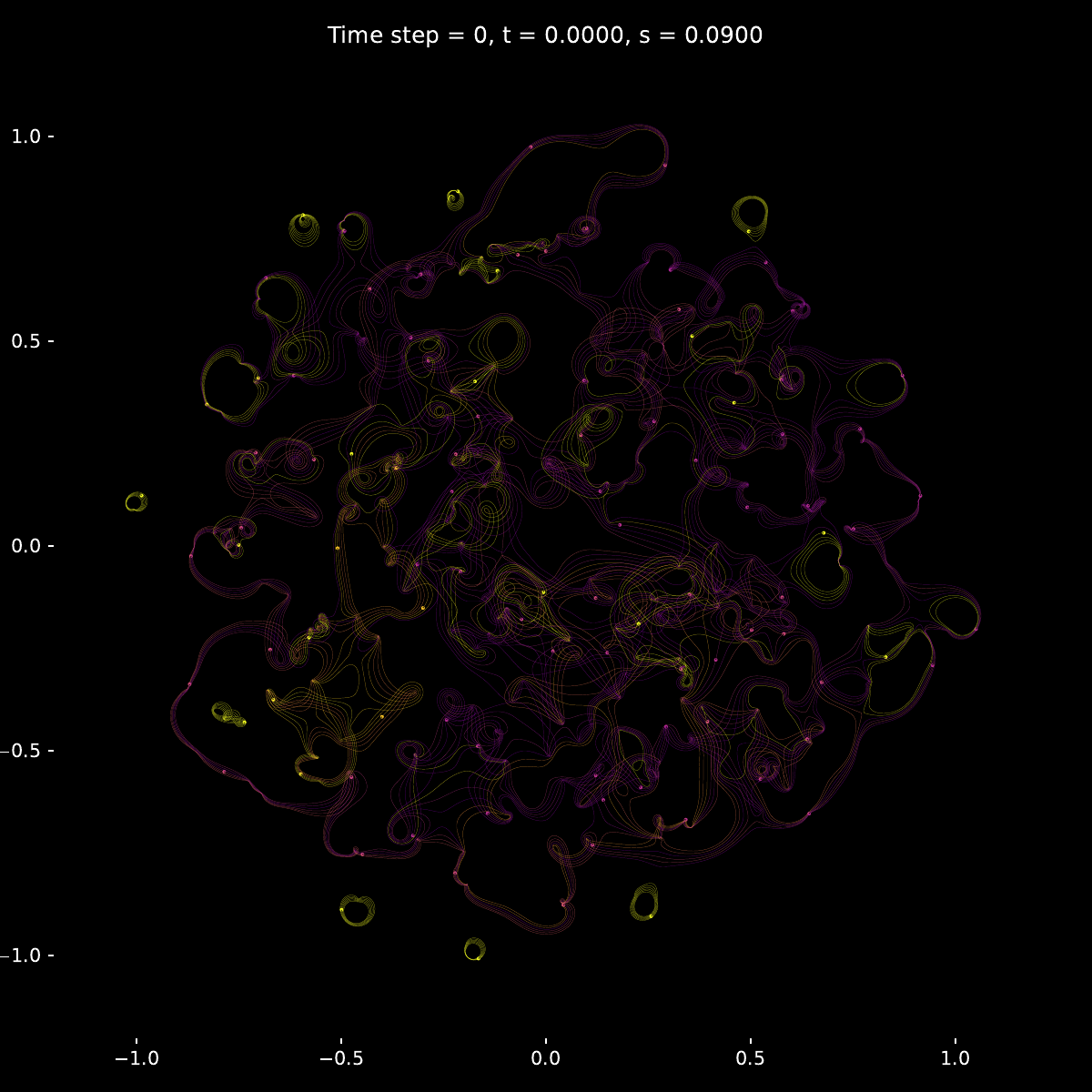}
		\caption{Eigenvalue trajectories for $s= 0.075, 0.080, 0.085, 0.090, 0.095, 0.100$ }
		\label{fig:pdf_image}
	\end{figure}

	\newpage

	\begin{figure}[htbp]
		\centering
		\includegraphics[width=0.8\textwidth]{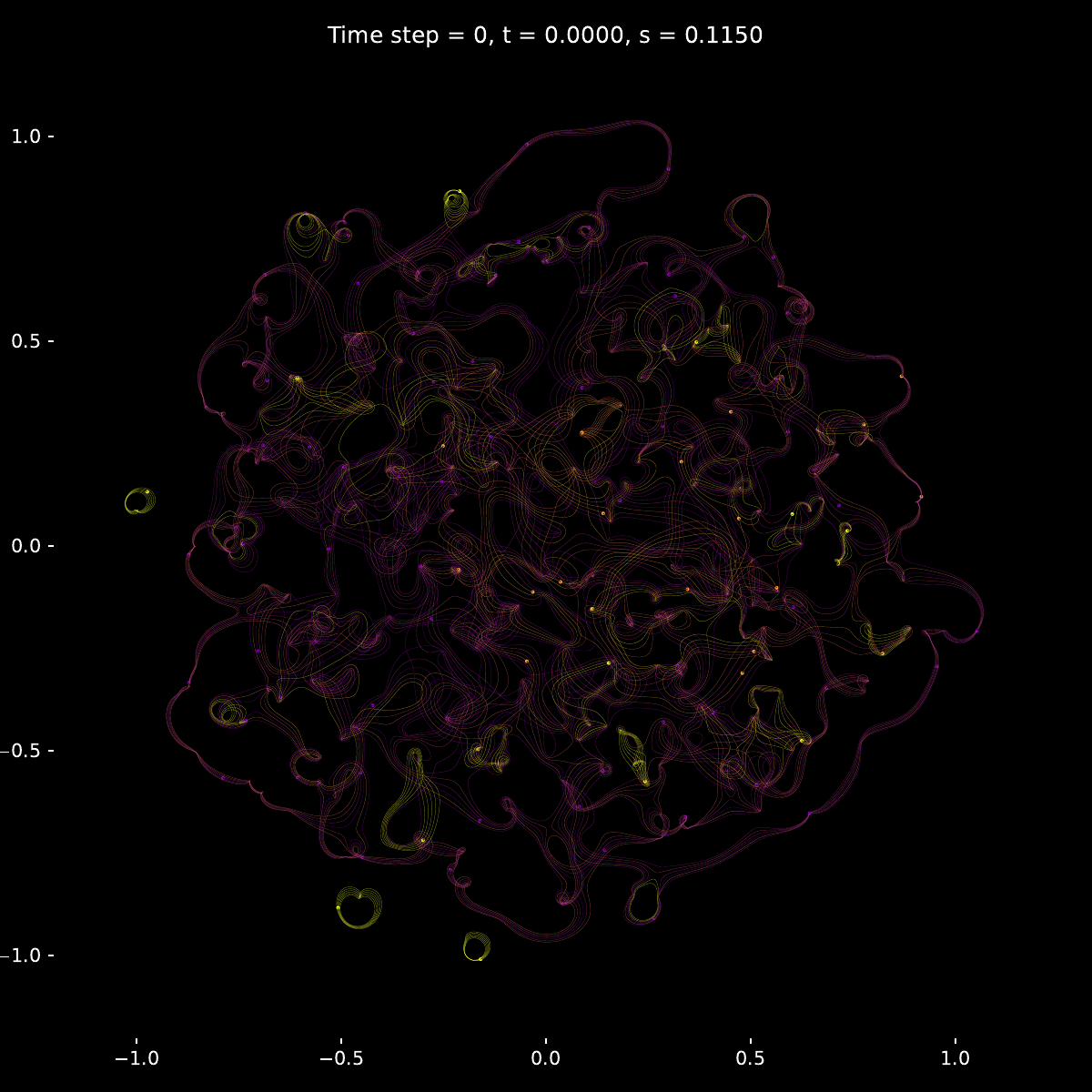}
		\caption{Eigenvalue trajectories for $s= 0.100, 0.105, 0.110, 0.115, 0.120, 0.125$ }
		\label{fig:pdf_image}
	\end{figure}

	\newpage

	\begin{figure}[htbp]
		\centering
		\includegraphics[width=0.8\textwidth]{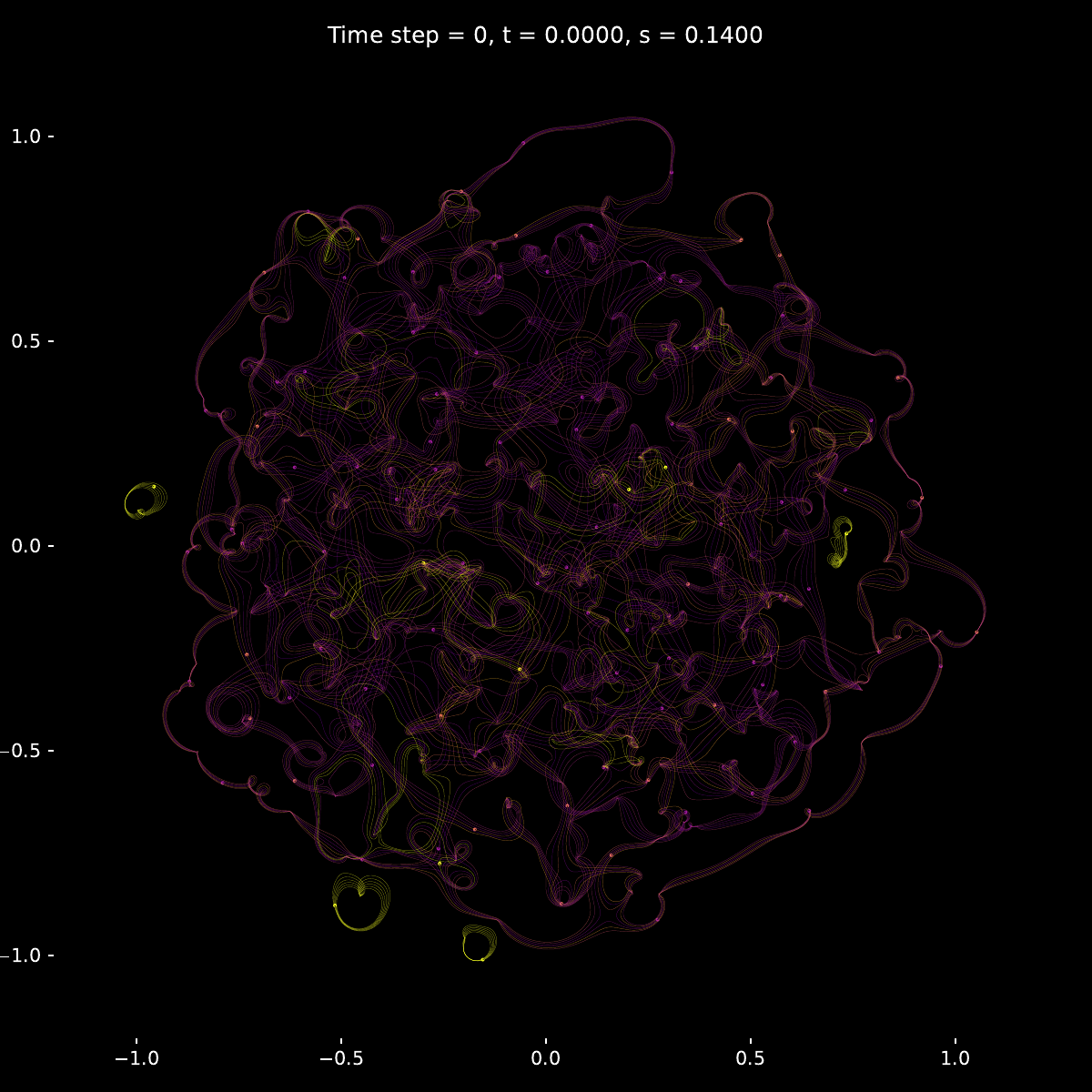}
		\caption{Eigenvalue trajectories for $s= 0.125, 0.130, 0.135, 0.140, 0.145, 0.150$ }
		\label{fig:pdf_image}
	\end{figure}

\end{document}